\documentclass[preprint,12pt]{elsarticle}

\usepackage{amssymb, bm}
\usepackage{amsthm, amsmath}

\usepackage{graphicx,color}
\usepackage{caption}
\usepackage{subcaption}
\usepackage{multirow}
\usepackage{mathtools,xparse}
\usepackage{algorithm}
\usepackage{algorithmic}

\usepackage{listings}

\makeatletter
\newenvironment{smallermatrix}[1][c]
{\null\,\vcenter\bgroup
  \Let@\restore@math@cr\default@tag
  \baselineskip0pt \lineskip0.4pt \lineskiplimit0pt
  \ialign\bgroup\if#1l\else\hfil\fi$\m@th\scriptstyle##$\if#1r\else\hfil\fi&&\thickspace\hfil
  $\m@th\scriptstyle##$\hfil\crcr
}{%
  \crcr\egroup\egroup\,%
}
\makeatother

\NewDocumentCommand{\ts}{O{c} e{^?_}}{
  \begin{smallermatrix}[#1]
  \mathstrut\IfValueT{#2}{#2} \\
  \mathstrut\IfValueT{#3}{#3} \\
  \mathstrut\IfValueT{#4}{#4}
  \end{smallermatrix}%
}

\journal{Arxiv}

\begin{document}

\begin{frontmatter}



\title{Generative Reduced Basis Method}




\author[inst2]{Ngoc Cuong Nguyen}

\affiliation[inst2]{organization={Center for Computational Engineering, Department of Aeronautics and Astronautics, Massachusetts Institute of Technology},
            addressline={77 Massachusetts
Avenue}, 
            city={Cambridge},
            state={MA},
            postcode={02139}, 
            country={USA}}
            

\begin{abstract}
We present a generative reduced basis (RB) approach to construct reduced order models for parametrized partial differential equations. Central to this approach is the construction of generative RB spaces that provide rapidly convergent approximations of the solution manifold. We introduce a generative snapshot method to generate significantly larger sets of snapshots from a small initial set of solution snapshots. This method leverages multivariate nonlinear transformations to enrich the RB spaces, allowing for a more accurate approximation of the solution manifold than commonly used techniques such as proper orthogonal decomposition and greedy sampling. The key components of our approach include ({\em i}) a Galerkin projection of the full order model onto the generative RB space to form the reduced order model; ({\em ii}) {\em a posteriori} error estimates to certify the accuracy of the reduced order model; and ({\em iii}) an offline-online decomposition to separate the computationally intensive model construction, performed once during the offline stage, from the real-time model evaluations performed many times during the online stage. The error estimates allow us to efficiently explore the parameter space and select parameter points that maximize the accuracy of the reduced order model. Through numerical experiments, we demonstrate that the generative RB method not only improves the accuracy of the reduced order model but also provides tight error estimates.

\end{abstract}







\end{frontmatter}


\section{Introduction}
\label{sec:intro}


Reduced basis (RB) methods are widely used for solving parametrized partial differential equations (PDEs)  \cite{Chen2009, Chen2010a, deparis07, EftangPateraRonquist10,Eftang2012a, Huynh2007a, Huynh2010, PhuongHuynh2013, Karcher2018, Knezevic2011, Knezevic2010,  nguyen04:_handb_mater_model, Nguyen2007, Nguyen_SantaFE08, Calcolo, ARCME, Rozza05_apnum, Sen2006b, veroy04:_certif_navier_stokes, Veroy2002, Vidal-Codina2014, Vidal-Codina2018a}. One of the fundamental concepts underlying the effectiveness of RB methods is the Kolmogorov $n$-width, which measures how well a function space can be approximated by an $n$-dimensional subspace. The Kolmogorov $n$-width quantifies the smallest possible error when approximating the solution manifold of a parametrized PDE with an $n$-dimensional subspace \cite{Maday2008c}. It offers a theoretical framework for understanding the potential efficiency of a reduced basis: the lower the $n$-width for a given $n$, the more effectively the reduced basis can approximate the solution manifold. 

The Kolmogorov $n$-width theory is related to Proper Orthogonal Decomposition (POD) \cite{sirovich87:_turbul_dynam_coher_struc_part,Kunisch2001,Nguyen2008a}. POD provides a practical method for constructing low-dimensional approximation spaces by extracting the most energetic modes from a set of solution snapshots.  For problems with smooth solution manifolds, where the Kolmogorov $n$-width decays rapidly, POD offers an effective low-dimensional approximation. However, POD faces limitations when the solution manifold contains localized, non-smooth, highly oscillatory features. In such cases, POD requires a large number of snapshots to achieve acceptable accuracy because the $n$-width decays slowly.

 
Greedy algorithms have been introduced to address some of the limitations of POD \cite{tanbui07_greedy,Buffa2012a,bui07:_goal_pod,Hesthaven2014,DeVore2013,Haasdonk2013,Hoang2013,Binev2011,Nguyen2016,Lieberman2010}. The greedy algorithm, combined with a posteriori error estimates, provides an efficient way to explore the parameter space and select the most informative parameter points. Both POD and the greedy algorithm encounter difficulties when applied to problems with slowly decaying Kolmogorov $n$-width or high-dimensional parameter spaces. In such cases, the solution manifold may not be adequately captured by a small set of snapshots, limiting the accuracy of these methods.

 
To improve the accuracy and efficiency of ROMs, several extensions to traditional RB methods have been proposed. Adaptive RB approach partitions the parameter domain and constructs an individualized ROM for each partition \cite{EftangPateraRonquist10,Haasdonk2011,Hess2019,Maday2013a}.  Other adaptive RB techniques have been developed to update the RB space during the online phase according to various criteria associated with changes of the system dynamics in parameters and time \cite{Hesthaven2022a,Peherstorfer2015,Peherstorfer2015a,Sapsis2009,Zimmermann2018}. Multiscale RB method \cite{BoyavalSIAM08,Nguyen2008} and static condensation RB element method \cite{Eftang2013,PhuongHuynh2013,Smetana2015} have been developed to handle problems in which the  coefficients of the differential operators are characterized by a large number of independent parameters. Boyaval et al. \cite{Boyaval2010,Boyaval2009b} develop a RB approach for solving variational problems with stochastic parameters  and reducing the variance in the Monte-Carlo simulation of a stochastic differential equation \cite{Vidal-Codina2014,Vidal-Codina2016}.

RB methods applied to convection-dominated problems often struggle due to slowly decaying Kolmogorov $n$ widths caused by moving shocks and discontinuities. To address this issue, various techniques recast the problem in a more suitable coordinate frame using parameter-dependent maps. Approaches include POD-Galerkin methods with shifted reference frames \cite{Rowley2003}, optimal transport methods \cite{Iollo2014,Heyningen2023,nguyen2023optimal,Nguyen2024,Nguyen2024d}, and phase decomposition. Recent methods like shifted POD \cite{Reiss2018}, transport reversal \cite{Rim2018}, and transport snapshot \cite{Nair2019} use time-dependent shifts, while registration \cite{Taddei2020} and shock-fitting methods  \cite{Zahr2018,Zahr2020} minimize residuals to determine the map for better low-dimensional representations. 

Model reduction on nonlinear manifolds was originated with the work \cite{Rutzmoser2017} on reduced-order modeling of nonlinear structural dynamics using quadratic manifolds. This approach was further extended through the use of deep convolutional autoencoders by Lee and Carlberg \cite{Lee2020}. Recently, quadratic manifolds are further developed to address the challenges posed by the Kolmogorov barrier in  model order reduction \cite{Geelen2023,Barnett2022}. Another approach is model reduction through lifting or variable transformations, as demonstrated by Kramer and Willcox \cite{Kramer2019}, where nonlinear systems are reformulated in a quadratic framework, making the reduced model more tractable. This method was expanded in the ``Lift \& Learn'' framework \cite{Qian2020} and operator inference techniques \cite{McQuarrie2023,Kramer2024}, leveraging physics-informed machine learning for large-scale nonlinear dynamical systems.

Another advancement aimed at improving the efficiency of ROMs is the development of empirical interpolation methods (EIM). EIM was introduced to deal with non-affine parameter dependencies and nonlinearities  \cite{Barrault2004, Grepl2007}.  EIM constructs an approximation of nonlinear or non-affine terms by selecting interpolation points and functions to form an interpolation model. EIM has been widely used to construct efficient ROMs for nonaffine and nonlinear PDEs \cite{Grepl2007,Nguyen2007,Galbally2010,Eftang2010b,Drohmann2012,Manzoni2012,Hesthaven2014,Kramer2019,Hesthaven2022,Chen2021}. Several attempts have been made to extend the EIM in various ways. The best-points interpolation method (BPIM) \cite{Nguyen2008a,Nguyen2008d} employs POD to generate the basis set and least squares to compute the interpolation point set.  Generalized empirical interpolation method (GEIM) \cite{Maday2015a,Maday2013} generalizes EIM by replacing the pointwise function evaluations by more general measures defined as linear functionals.

The first-order empirical interpolation method (FOEIM), introduced in \cite{Nguyen2023d,Nguyen2024,nguyen2024firstorderempiricalinterpolationmethod}, extends the standard EIM by using partial derivatives of parametrized nonlinear functions to generate additional basis functions and interpolation points. FOEIM has been further extended to a class of high-order empirical interpolation methods, which use higher-order partial derivatives to enrich the basis functions and interpolation points \cite{nguyen2024highorderempiricalinterpolationmethods}. High-order EIM offers more flexible and accurate approximations compared to standard EIM, significantly improving the accuracy of hyper-reduced ROMs by capturing non-affine and nonlinear behaviors while maintaining computational efficiency. 

In this paper, we introduce a generative reduced basis (RB) approach to constructing reduced-order models. Central to our approach is the construction of generative RB spaces, which offer a rapidly converging approximation of the solution manifold. We introduce the generative snapshot method to generate significantly larger sets of snapshots from a small initial set of solution snapshots. This  method enriches the RB spaces without requiring additional FOM solutions, thereby enabling more accurate and efficient approximations of the solution manifold. The method complements existing model reduction techniques such as POD, greedy sampling, adaptive RB methods, optimal transport methods,  nonlinear manifold methods, and nonintrusive and data-driven ROMs.

The proposed approach integrates several key components to enhance the accuracy and efficiency of ROMs. First, a Galerkin projection of the full-order model onto the generative RB spaces forms the reduced-order model. Second, we develop a posteriori error estimates that are inexpensive and tight. These error estimates are instrumental in guiding a greedy sampling algorithm, which effectively explores the parameter space and selects parameter points to optimize the accuracy of ROMs. Moreover, the approach leverages an offline-online decomposition to ensure computational efficiency. 

Through a series of numerical experiments, we show that the generative snapshot method produces enriched RB spaces approximating the solution manifold far more accurately than the standard RB space. Furthermore, generative RB methods significantly improve both the accuracy and reliability of reduced-order models (ROM). Our results show that the generative RB approach not only accelerates the convergence of ROMs but also provides tight error bounds, making it particularly effective for high-dimensional and complex parametric problems.

The paper is organized as follows. In Section 2, we provide an overview of dimensionality reduction methods for constructing linear approximation spaces. In Section 3, we introduce the generative snapshot method and present numerical results demonstrating its effectiveness. We propose a generative reduced basis method for affine linear PDEs in Section 4 and present numerical results to demonstrate the method in Section 5. Finally, in Section 6, we conclude with remarks on future work.

\section{Linear Approximation Spaces}



\subsection{Parametric Solution Manifold}

Let $\Omega \subset  \mathbb{R}^D$ represent the physical domain, where $\bm x \in \Omega$  denotes the spatial coordinate. This physical domain may represent, for example, the geometry of a region in which a physical process (e.g., fluid flow, heat transfer, or structural deformation) takes place. let $\mathcal{D} \subset \mathbb{R}^P$ denote the parameter domain, where  $\bm \mu$ represents a set of  $P$-dimensional parameters. These parameters encapsulate a wide range of physical, material, or geometrical properties of the system under consideration.


The function $u(\bm x, \bm \mu)$ is the solution of a parametrized partial differential equation (PDE), which could be linear or nonlinear, time-dependent or time-independent. We define the solution manifold $\mathcal{M}$ as the set of all possible solutions corresponding to every admissible parameter configuration in $\mathcal{D}$:
\begin{equation}
\mathcal{M} = \{u(\bm  \mu), \ \forall \bm \mu \in \mathcal{D}\} ,   
\end{equation}
where each element $u(\bm  \mu) \in X$ represents the spatial solution $u(\bm x, \bm  \mu)$ for a particular parameter point $\bm \mu$. Here $X$ is an appropriate Hilbert space equipped with an inner product $(\cdot, \cdot)_X$ and norm $\|v\|_X = \sqrt{(v,v)_X}$. The manifold $\mathcal{M}$ can be viewed as a high-dimensional object that captures the full variability of the solution across the parameter domain. This high-dimensional nature presents significant computational challenges, as computing the full solution manifold requires solving the PDE for many different parameter values, which can be computationally prohibitive in real-time or many-query contexts, such as optimization or uncertainty quantification.

Let $\Xi_K = \{\widehat{\bm \mu}_k \in \mathcal{D}\}_{k=1}^K$ be a set of $K$ distinct parameter points sampled from the parameter domain. These points are used to compute corresponding solutions to the parametrized PDE, $u(\hat{\bm \mu}_k)$, which are used to define the following finite-dimensional space
\begin{equation}
\mathcal{M}_{K} = \mbox{span} \{u(\hat{\bm \mu}_k), 1 \le k \le K\} .   
\end{equation}
We assume that $K$ is a large number such that the finite-dimensional space $\mathcal{M}_{K}$ can approximate well any function in the solution manifold $\mathcal{M}$.  However, constructing the space $\mathcal{M}_{K}$ is computationally expensive, since it requires a large number of solutions. The primary purpose of the finite-dimensional space $\mathcal{M}_{K}$  is to assess the quality of a RB space. 

\subsection{Kolmogorov $n$-Width}

Let $V_M = \mbox{span} \{v_m, 1 \le m \le M\}$ be an approximation space of $M$ dimension. We consider linear approximation in which any function $v \in V_M$ can be expressed as a linear combination of its basis functions, $v(\bm x) = \sum_{m=1}^M \alpha_m v_m (\bm x)$, where $\alpha_m, 1 \le m \le M,$ are real coefficients. For any given function $w \in X$, its best approximation in the space $V_M$ is defined as
\begin{equation}
v^* = \arg \inf_{v \in V_M} \|v - w\|_X .    
\end{equation}
The best approximation can be computed as $v^* = \sum_{m=1}^M \alpha^*_m v_m$, where $\bm \alpha^* = \arg \inf_{\bm \alpha \in \mathbb{R}^M} \|\sum_{m=1}^M \alpha_m v_m - w\|_X$. To assess the accuracy of the space $V_M$ in approximating the solution space $\mathcal{M}_K$, we define the following error metric: 
\begin{equation}
\label{errormetric}
d(V_M, \mathcal{M}_K) =  \sup_{u \in \mathcal{M}_K} \inf_{v \in V_M} \|v - u\|_X   .  
\end{equation}
This quantity measures the worst-case approximation error by calculating the maximum deviation between the solutions in  $\mathcal{M}_K$  and their best approximations within the space $V_M$. The smaller the value of the error metric he more accurate the space $V_M$ is in approximating the functions in $\mathcal{M}_K$.

For a fixed dimension $M$, we wish to find the space $V_M^* \in X$ of dimension $M$ that minimizes the worst-case approximation error across the entire solution manifold $\mathcal{M}_K$.  Mathematically, this optimal approximation space is defined as
\begin{equation}
V^*_M = \arg \inf_{V_M} \sup_{u \in \mathcal{M}_K} \inf_{v \in V_M} \|v - u\|_X  .  
\end{equation}
The associated error metric $d(V_M^*, \mathcal{M}_K)$ is known as the Kolmogorov $n$-width. The Kolmogorov $n$-width characterizes the smallest possible worst-case error that can be achieved using any  $M$-dimensional approximation space. This provides a fundamental limit on the performance of any linear approximation method of fixed dimension $M$. If the Kolmogorov $n$-width converges rapidly to zero as $M$ increases, it implies that the solution manifold $\mathcal{M}_K$ is close to being representable by a low-dimensional linear subspace. Conversely, if the Kolmogorov $n$-width converges slowly, it indicates that solution manifold is inherently complex and high-dimensional. 

\subsection{Proper Orthogonal Decomposition} 

While the Kolmogorov $n$-width gives a powerful theoretical limit, finding the optimal approximation space that achieves this bound is generally not computable in practice. The concept of finding the best approximation space is closely related to Proper Orthogonal Decomposition (POD). POD provides a practical method for constructing an approximation space $V_M$ that approximates the optimal space $V_M^*$ by identifying the most energetic modes in the solution manifold. These modes correspond to the principal components of the solution manifold and can be obtained by applying singular value decomposition (SVD) to the snapshot data.

POD generates an orthogonal basis that minimizes the projection error in the $L_2$ norm. Mathematically, the POD basis functions solve the following minimization problem
\begin{equation}
\min_{V_M} \sum_{k=1}^K \left\| u(\bm \mu_k) - \sum_{m=1}^M (u(\bm \mu_k), v_m)_X \, v_m \right\|_X^2 
\end{equation}
where $V_M$ is constrained to be an orthogonal space of dimension $M$. This minimization problem amounts to solving the following eigenvalue problem
\begin{equation}
\bm C \bm a = \lambda \bm a
\end{equation}
where $\bm C \in \mathbb{R}^{K \times K}$ is the covariance matrix with entries $C_{kk'} = (u(\bm \mu_k), u(\bm \mu_{k'}))_X$, $1 \le k, k' \le K$. The eigenvectors  provide the coefficients that describe the linear combination of the snapshots that form the optimal POD basis. Once the eigenvalue problem has been solved, the POD basis functions can be computed as a linear combination of the original snapshots. The
$m$-th POD basis function $v_m$ is constructed as 
\begin{equation}
v_m = \sum_{k=1}^K a_{mk} u(\bm \mu_k), 
\end{equation} 
where $a_{mk}$ are components of the eigenvector associated with the $m$-th largest eigenvalue.





\subsection{Greedy Sampling Algorithm}

While POD provides an orthogonal basis that closely approximates the optimal space $V_M^*$, constructing this basis requires access to the entire solution manifold  $\mathcal{M}_K$. One of the primary challenges in constructing the POD basis is the high computational cost of generating the snapshots. One effective approach to reduce the number of required snapshots is the use of greedy algorithm. This method begins with a small set of snapshots and adds new ones in an iterative fashion by selecting those that maximally reduce the approximation error at each step. Greedy algorithm relies on inexpensive error estimates in order to effectively explore the parameter space.

The greedy algorithm starts with a small set of parameter points $\{\bm \mu_i\}_{i=1}^n$ randomly selected in the parameter domain. These initial parameter points are used to compute the corresponding solutions, $\{u(\bm \mu_n)\}_{i=1}^N$, which span the initial RB space $W_n = \mbox{span} \{u(\bm \mu_n), 1 \le i \le n\}$. It finds the next parameter point $\bm \mu_{n+1}$ that maximizes an error estimate $\varepsilon_n(\bm \mu)$ over a training set $\Xi_{\rm train}$: 
\begin{equation}
\bm \mu_{n+1} = \arg \max_{\bm \mu \in \Xi_{\rm train}} \varepsilon_n(\bm \mu) ,   
\end{equation}
where $\varepsilon_n(\bm \mu)$ is an estimate of the true error $\|u(\bm \mu) - u_n(\bm \mu) \|_X$, with $u_n(\bm \mu) \in W_n$ being an appropriate approximation to $u(\bm \mu)$. The training sample $\Xi_{\rm train}$ can be either a very fine grid in the parameter domain or an exhaustive list of randomly selected parameter points. Then the RB space $W_{n+1}$ is augmented by adding $u(\bm \mu_{n+1})$ to $W_n$. This process is repeated until $\varepsilon_N(\bm \mu)$ is less than a specified tolerance, where $N$  is the total number of basis functions in the final RB space $W_N = \mbox{span} \{u(\bm \mu_n), 1 \le i \le N\}$. 

Although the greedy algorithm is widely used for constructing RB spaces, it has several limitations. One of the key challenges lies in the computation of error estimates for numerous parameter points within the training set. For the algorithm to be effective, these error estimates must be inexpensive, rigorous, and tight. However, achieving all three properties simultaneously is often difficult—error estimates that are inexpensive may lack rigor or accuracy, while rigorous and tight estimates can be computationally expensive. Each time a new parameter point  is selected, the corresponding solution must be computed, which requires solving the parametrized PDE. This can be costly, especially when many iterations of the greedy algorithm are required. In high-dimensional parameter spaces, the size of the training set grows exponentially with the number of parameters. This makes it difficult to sample the parameter space adequately, leading to potential gaps in the training set. If the training set does not adequately cover the parameter space, the algorithm might fail to capture critical parameter regions where the solution exhibits significant variation. As a result, the reduced basis might not provide accurate approximations for certain parameter values. The convergence of the algorithm is strongly dependent on the smoothness and structure of the solution manifold. For certain problems with complex solution manifolds, the greedy algorithm may converge slowly, necessitating a large number of basis functions to achieve the desired accuracy.

\section{Generative Snapshot Method}

We introduce a method for generating an expanded set of snapshots from a small initial set of solutions by applying multivariate nonlinear transformations. These transformations enrich the original set of solutions, creating new snapshots that capture more complex patterns within the solution manifold. After generating these new snapshots, POD is applied to construct generative RB spaces. The key advantage of this method is that it significantly enlarges the snapshot space without requiring additional solutions. The method allows for more accurate approximations of the solution manifold by producing high-quality RB spaces. To demonstrate the effectiveness of the generative snapshot method, we present a series of numerical examples that compare the performance of generative RB spaces to that of the traditional RB space.

\subsection{Multivariate Nonlinear Transformations}


Let $\sigma(v)$ be a proper nonlinear function of a scalar variable $v$. We assume that $\sigma$ is smooth and at least twice differentiable with respect to $v$,  and that the inverse function $\sigma^{-1}(w)$ exists. We then introduce the following two-variable function
\begin{equation}
G(v_1, v_2) = \sigma(v_1) + \frac{\partial \sigma(v_1)}{\partial v}(v_2 - v_1) 
\end{equation}
and define the corresponding space
\begin{equation}
\mathcal{G}_{N^2} = \mbox{span} \{\rho_{mn} \equiv  \sigma^{-1}(G(\xi_n, \xi_{m})),  \ 1 \le n, m \le N \} .
\end{equation} 
Here $\{\xi_n \equiv u(\bm \mu_n)\}_{n=1}^N$ is the set of $N$ solutions defining the  RB space $W_N$:
\begin{equation}
W_N = \mbox{span} \{\xi_n \equiv u(\bm \mu_n), 1 \le i \le N\}  .  
\end{equation}
It can be easily shown that $W_N \subset \mathcal{G}_{N^2}$. The dimension of $\mathcal{G}_{N^2}$ may far exceed $N$, but it is less than or equal to $N^2$. 

We further introduce a three-variable function to capture higher-order nonlinear interactions:
\begin{equation}
H(v_1, v_2, v_3) = \sigma(v_1) + \frac{\partial \sigma(v_1)}{\partial v}(v_2 - v_1) + \frac{1}{2} \frac{\partial^2 \sigma(v_1)}{\partial v^2}(v_2 - v_1) (v_3 - v_1) ,
\end{equation}
and define the associated space
\begin{equation}
\mathcal{H}_{N^3} = \mbox{span} \{\varrho_{mnp} \equiv  \sigma^{-1}(H(\xi_n, \xi_{m}, \xi_p)),  \ 1 \le n, m, p \le N \} .
\end{equation}
It can be demonstrated that $W_N \subset \mathcal{G}_{N^2} \subset \mathcal{H}_{N^3}$. The space $\mathcal{G}_{N^2}$ contains $N^2$ functions, while $\mathcal{H}_{N^3}$ extends this by incorporating $N^3$ functions, offering a richer and more accurate representation of the solution manifold.

This hierarchical structure of spaces $W_N \subset \mathcal{G}_{N^2} \subset \mathcal{H}_{N^3}$ enables a flexible approach for generating snapshots, allowing for progressively more detailed approximations of the solution space and construction of efficient and accurate ROMs. In addition to enhancing ROM accuracy, this approach also enables the development of efficient error estimation techniques. The hierarchical structure facilitates the generation of tight error bounds by comparing the approximations across the nested spaces.

\subsection{Generative Reduced Basis Spaces}

Due to the potentially large number of snapshots in $\mathcal{G}_{N^2}$ and $\mathcal{H}_{N^3}$, generating basis functions via POD may become computationally intensive. We can reduce the number of snapshots by using the nearest parameter points as follows. We compute the distances $d_{nn'} = \|\bm \mu_n - \bm \mu_{n'}\|, 1 \le n,n' \le N,$ between each pair of parameter points $\bm \mu_{n}$  and  $\bm \mu_{n'}$ in $S_N$. For any given $\bm \mu_n \in S_N$, we define $S_L(\bm \mu_n)$ a set of exactly $L$ parameter points that are closest to $\bm \mu_n$. For any solution $\xi_n \in W_N$, we define $W_L(\xi_n) = \{u(\bm \mu_l), \, \forall \bm \mu_l \in S_L(\bm \mu_l) \}$, which contains  exactly $L$ solutions. We then introduce the following sets of snapshots
\begin{equation}
\begin{split}
\mathcal{G}_{NL} & = \mbox{span} \{\rho_{nl} \equiv  \sigma^{-1}(G(\xi_n, \xi_{l})),  \forall \xi_l \in W_L(\xi_n), \ 1 \le n\le N \}, \\
\mathcal{H}_{NL^2} & = \mbox{span} \{\varrho_{nll'} \equiv  \sigma^{-1}(H(\xi_n, \xi_{l}, , \xi_{l'})),  \forall \xi_l, \xi_{l'} \in W_L(\xi_n), \ 1 \le n\le N \} .
\end{split}
\end{equation} 
These snapshot sets satisfy $W_N \subset \mathcal{G}_{NL} \subset \mathcal{H}_{NL^2} \subset \mathcal{H}_{N^3}$ and contain $N L$ and $N L^2$ functions, respectively.

We apply POD to both $\mathcal{G}_{NL}$ and $\mathcal{H}_{NL^2}$ in order to compress the functions within these spaces. Thus, we obtain two sets of orthogonal basis functions that efficiently approximate the solution manifold. These basis functions provide a more compact representation for the reduced basis approximation while preserving accuracy. Specifically, the compressed spaces corresponding to $\mathcal{G}_{NL}$ and $\mathcal{H}_{NL^2}$ are given by
\begin{equation}
\Phi_N^{M_1} = \mbox{span} \{\phi_m, 1 \le m \le M_1 \}, \quad \Psi_N^{M_2} = \mbox{span} \{\psi_m, 1 \le m \le M_2 \} ,
\end{equation}
where the superscripts $M_1$ and $M_2$ denote the dimensions of the two RB spaces, and the subscript $N$ indicates the number of solutions used to construct these spaces.  The dimension $M_1$ (or $M_2$) is chosen such that the sum of the eigenvalues corresponding to the first $M_1$ (or $M_2$) POD modes captures most of the total sum of all the eigenvalues. Both $M_1$ and $M_2$ can be far greater than $N$. Typically, $M_2$ is  greater than $M_1$ due to the fact $\mathcal{G}_{NL} \subset \mathcal{H}_{NL^2}$.

Due to $W_N \subset \mathcal{G}_{NL} \subset \mathcal{H}_{NL^2}$, we can expect that the error metrics for the three different RB spaces satisfy the below relationship  
\begin{equation}
    d(\Psi_N^{M_2}, \mathcal{M}_K) \le d(\Phi_N^{M_1}, \mathcal{M}_K) \le d(W_N, \mathcal{M}_K) .
\end{equation} 
Hence, the generative RB spaces should provide better approximations than the standard RB space. The convergence rates of these error metrics as a function of $N$ reflect the performance of the generative RB spaces, $\Phi_N^{M_1}$ and $\Psi_N^{M_2}$, relative to the standard RB space $W_N$.



\subsection{Nonlinear Activation Functions}

The nonlinear function $\sigma$ plays a crucial role in the generative snapshot method, as it facilitates the transformation of a small set of original solutions into a significantly larger set of enriched snapshots through multivariate nonlinear operations. By introducing nonlinearity into the generated snapshots, $\sigma$ enables the new set to capture more complex and intricate relationships within the solution manifold that are not easily represented by the original solutions alone. This nonlinearity is essential for expanding the expressiveness of the reduced basis (RB) space, as it allows the method to approximate the solution manifold more effectively.

The function $\sigma$ must satisfy two key requirements to ensure the method's effectiveness. First, $\sigma$ must be at least twice differentiable to ensure that the transformations and expansions generated by $\sigma$ are smooth and stable, allowing for accurate computation of higher-order interactions when constructing enriched snapshots. Second, the inverse function $\sigma^{-1}$ must exist and be well-defined, ensuring that the nonlinear transformations can be consistently inverted when required. The existence of an inverse is crucial for defining the transformed snapshot spaces and maintaining consistency in the generation of new basis functions.

The function $\sigma$ not only transforms individual solutions but also enables the generation of snapshots that incorporate higher-order correlations between multiple solutions. For example, in the generative snapshot method, the nonlinear function is applied to pairs or triples of solutions to create spaces of $N^2$ or $N^3$ snapshots. This multivariate expansion enhances the richness of the snapshot space and improves the approximation accuracy of ROMs by capturing interactions that would otherwise be missed in a purely linear approximation approach. It allows the generative snapshot method to construct high-quality RB spaces that significantly improve the accuracy and generalization capabilities of ROMs.

Beyond these fundamental requirements, the choice of $\sigma$ can significantly impact the effectiveness of the method. A carefully chosen nonlinear function can balance complexity and computational efficiency, ensuring that the enriched snapshot set provides substantial improvements in accuracy without introducing excessive computational overhead. Common candidates for $\sigma$ include smooth polynomial functions, sigmoidal functions, expotential functions, or even more complex functions derived from machine learning models, depending on the nature of the solution space. 



In this paper, we explore several options for the nonlinear function $\sigma$, which play a crucial role in the generative snapshot method by introducing the necessary nonlinearity to capture the solution manifold. As shown in Table \ref{tab1}, we consider a variety of commonly used activation functions, including hyperbolic tangent, sigmoid, arctangent, softplus,  exponential, and quadratic functions. These functions are smooth and continuously differentiable, making them well-suited for mathematical models that require  gradient calculations and smooth transformations. They are widely employed in neural networks due to their ability to model nonlinear relationships, enabling the network to learn complex data distributions. In the context of the generative snapshot method, these nonlinear functions are natural candidates for transforming the original set of solutions into a larger, enriched snapshot set that accurately captures the solution manifold.


\begin{table}[h!]
\centering
\renewcommand{\arraystretch}{2}
\begin{tabular}{|c|c|c|c|c|}
\hline
$\sigma(x)$ & $\partial \sigma(x)/\partial x$  & $\partial^2 \sigma(x)/\partial x^2$ & $\sigma^{-1}(x)$\\ \hline

 $\tanh(x)$ 
& $1 - \sigma^2(x)$ 
& $2 \sigma(x)(\sigma^2(x) - 1)$ 
& $ \frac{1}{2} \ln \left( \frac{1+x}{1-x} \right)$ \\ \hline

 $\frac{1}{1 + e^{-x}}$
& $\sigma(x)(1 - \sigma(x))$ 
& $\sigma(x)(1 - \sigma(x))(1 - 2\sigma(x))$ 
& $ \ln\left( \frac{x}{1 - x} \right)$ \\ \hline

 $\arctan(x)$
& $\frac{1}{1 + x^2}$ 
& $\frac{-2x}{(1 + x^2)^2}$ 
& $ \tan(x)$ \\ \hline

 $\ln(1 + e^x)$
& $\frac{1}{1 + e^{-x}}$ 
& $\frac{e^x}{(1 + e^x)^2}$ 
& $ \ln(e^x - 1)$ \\ \hline

 $e^x$
& $e^x$
& $e^x$
& $ \ln(x)$ \\ \hline

 $x^2$
& $2x$
& $2$
& $ \sqrt{x}$ \\ \hline

\end{tabular}
\caption{Nonlinear activation functions, their first and second-order derivatives, and inverse functions.}
\label{tab1}
\end{table}

\subsection{One-Dimensional Parametrized Function}

\begin{figure}[htbp]
	\centering
 \includegraphics[width=0.55\textwidth]{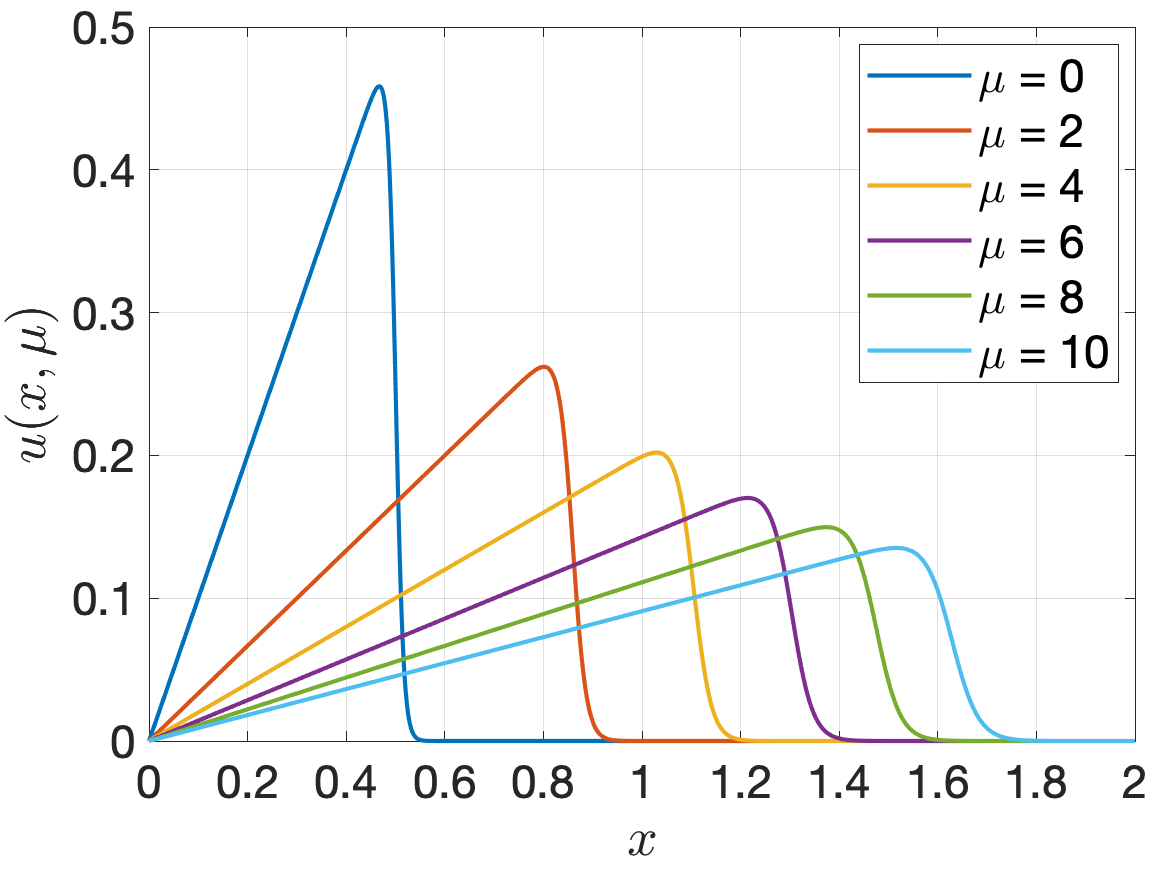}
 \caption{Plots of $u(x,\mu)$ defined in (\ref{ex1u}) as a function of $x$ for  different values of $\mu$.}
	\label{ex1fig1}
\end{figure}

We present numerical results from a simple test case to assess the performance of the generative snapshot method. The test case involves the following parametrized function 
\begin{equation}
\label{ex1u}
u(x, \mu) = \frac{x}{(\mu + 1) \left(1 + \sqrt{\frac{\mu+1}{\exp(62.5)}} \exp \left(\frac{125 x^2}{\mu+1} \right) \right) }    
\end{equation}
in a physical domain $\Omega = [0,2]$ and parameter domain $\mathcal{D} = [0, 10]$. We choose the function space $X = L^2(\Omega)$ with the inner product $(w, v)_X = \int_{\Omega} w v d x$. Figure \ref{ex1fig1} shows the plots of $u$ for various values of $\mu$. The solution has a sharp shock profile at $x =0.5$ for $\mu = 0$. As $\mu$ increases, this shock profile moves to the right boundary and smears out. 

The sample $S_N$ consists of $N$ parameter points from an extended Chebyshev distribution over the parameter domain $\mathcal{D}$. Thus, the standard RB space $W_N$ consists of $N$ snapshots that correspond to these extended Chebyshev points. The generative RB spaces, $\Phi_N^{M_1}$ and $\Psi_N^{M_2}$, are also constructed from these $N$ solutions by using the generative snapshot method for two different values of $L$, which are $L = \min(4,N)$ and $L = \min(8,N)$. We assess the performance of these RB spaces by studying the convergence of the error metric defined in Equation (\ref{errormetric}) as a function of $N$. To this end, we employ $K=200$ parameter points sampled uniformly over $\mathcal{D}$. Consequently, the discrete solution manifold $\mathcal{M}_K$ is composed of 200 snapshots. 

\begin{figure}[htbp]
	\centering
	\begin{subfigure}[b]{0.32\textwidth}
		\centering		\includegraphics[width=\textwidth]{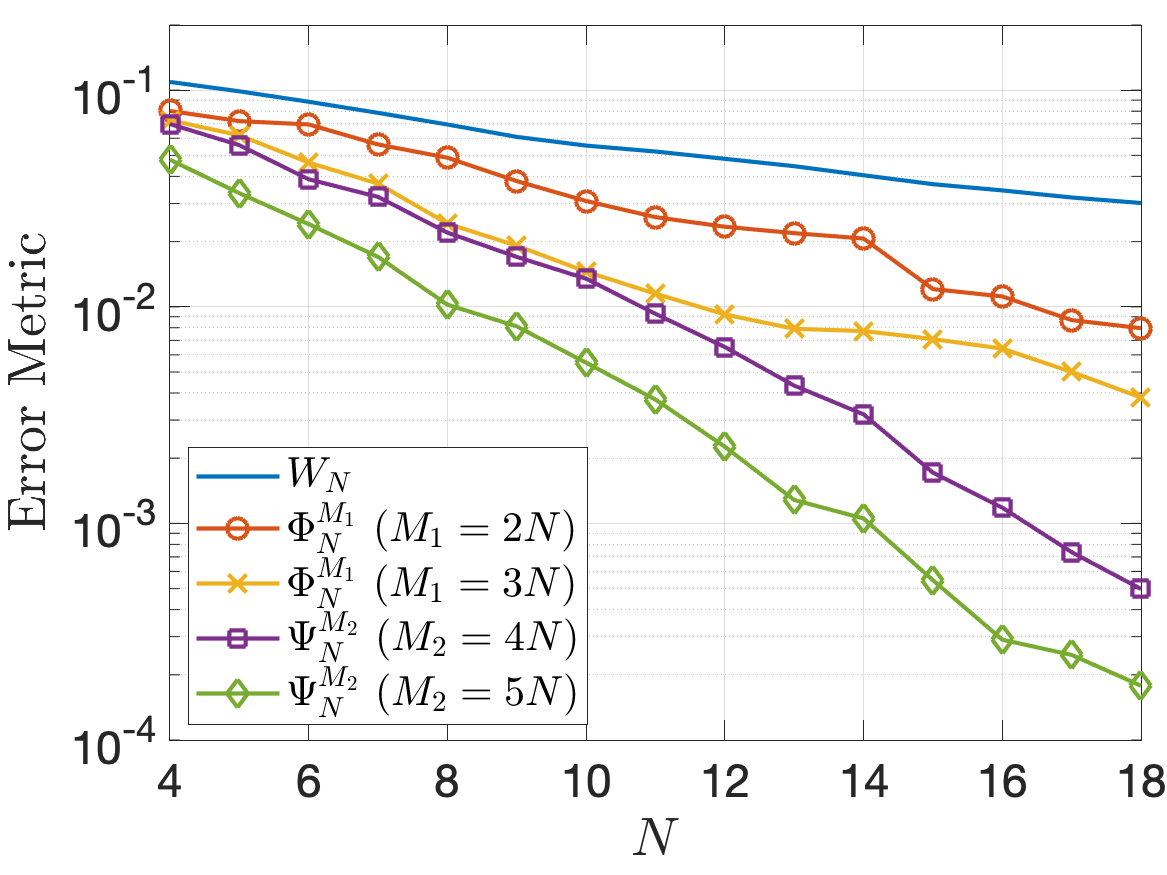}
		\caption{Hyperbolic tangent function}
	\end{subfigure}	
	\begin{subfigure}[b]{0.32\textwidth}
		\centering		\includegraphics[width=\textwidth]{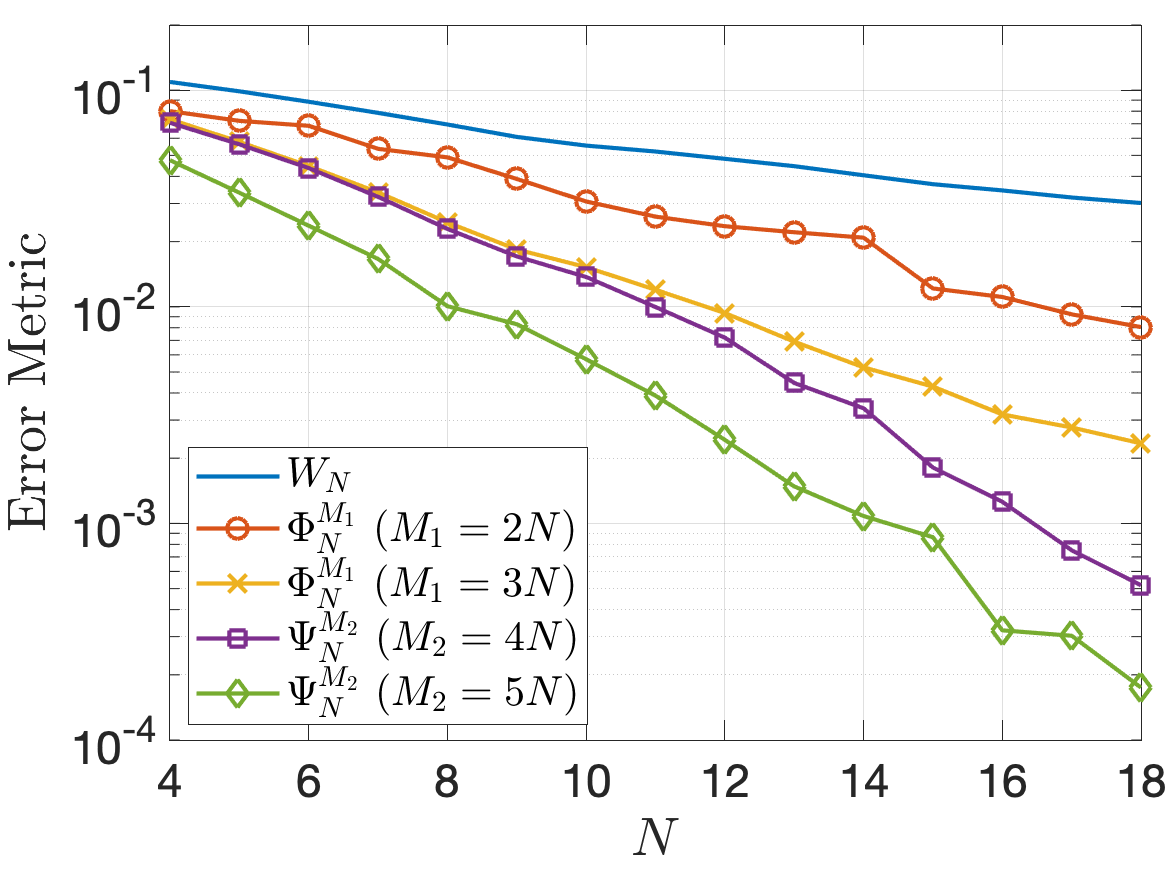}
		\caption{Sigmoid function}
	\end{subfigure}
 	\begin{subfigure}[b]{0.32\textwidth}
		\centering		\includegraphics[width=\textwidth]{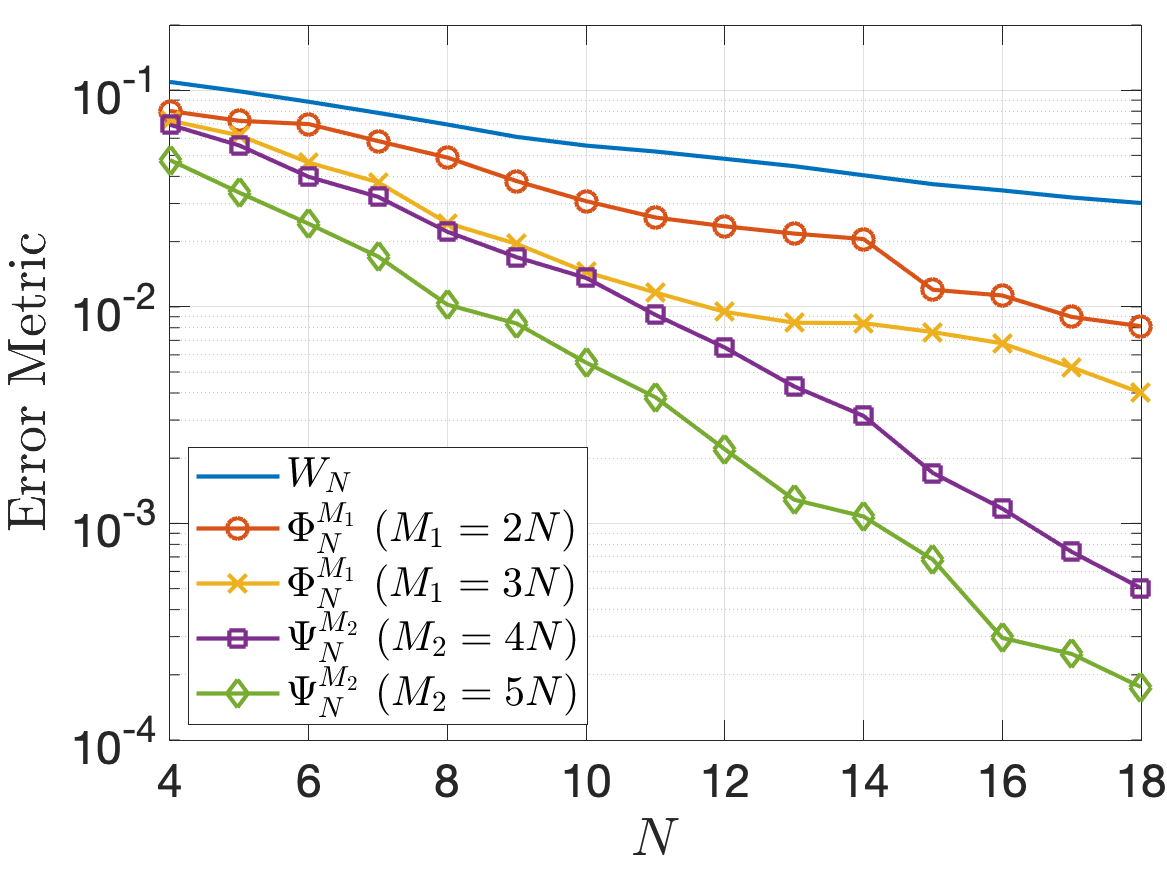}
		\caption{Arctangent function}
	\end{subfigure}	
	\begin{subfigure}[b]{0.32\textwidth}
		\centering		\includegraphics[width=\textwidth]{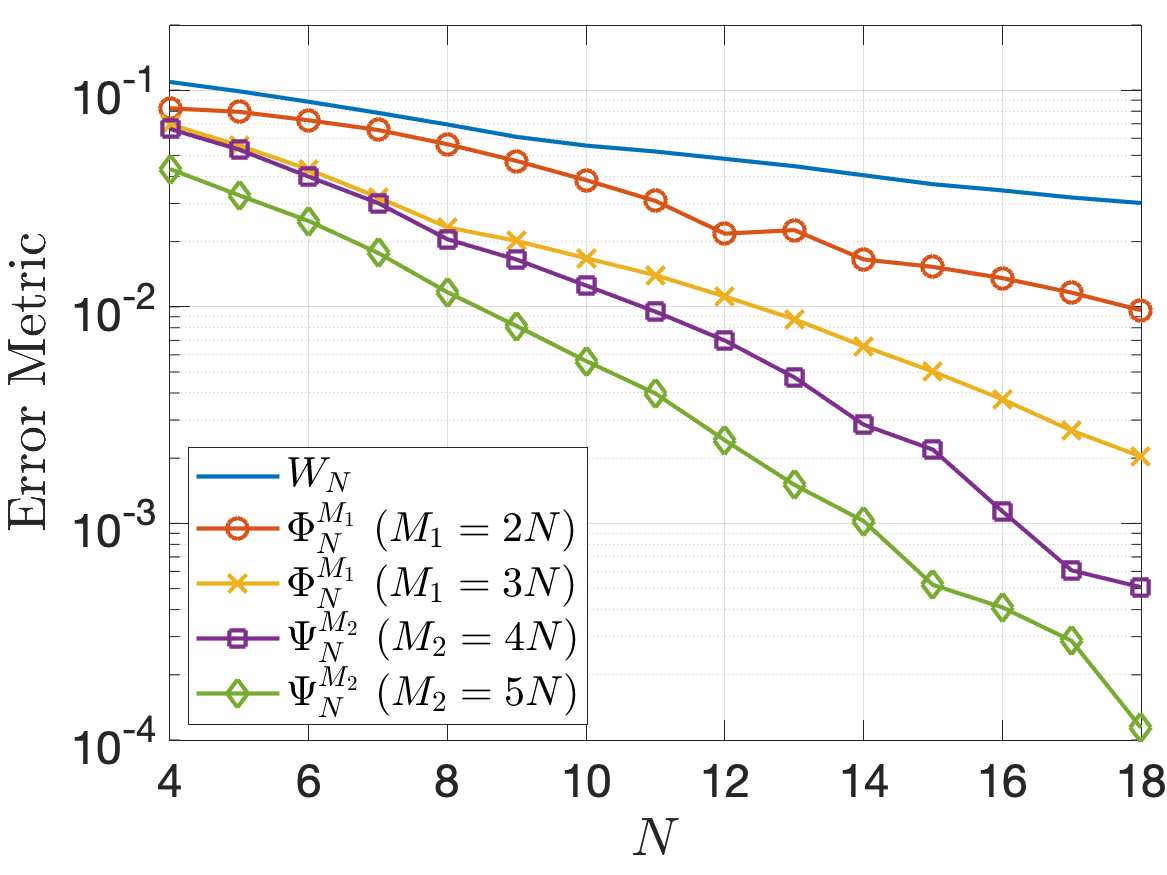}
		\caption{Softplus function}
	\end{subfigure}
  	\begin{subfigure}[b]{0.32\textwidth}
		\centering		\includegraphics[width=\textwidth]{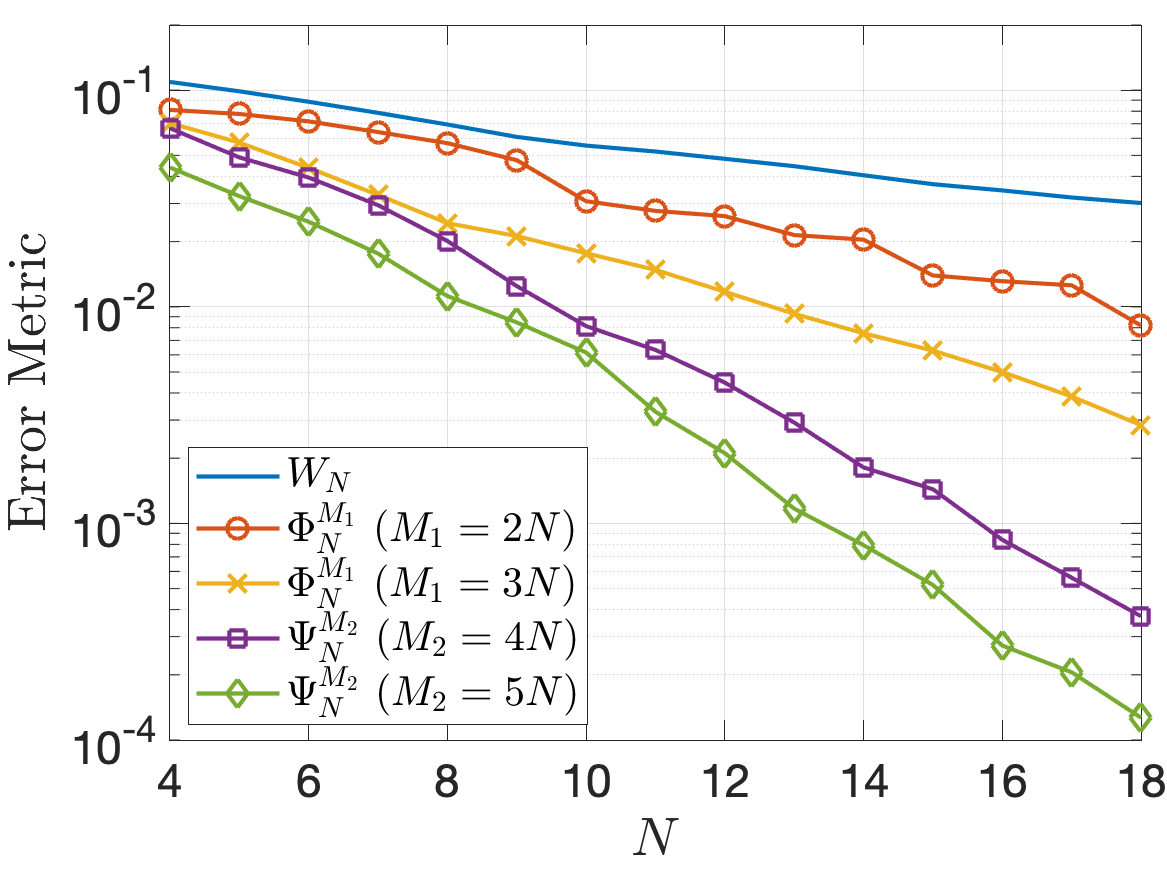}
		\caption{Exponential function}
	\end{subfigure}	
	\begin{subfigure}[b]{0.32\textwidth}
		\centering		\includegraphics[width=\textwidth]{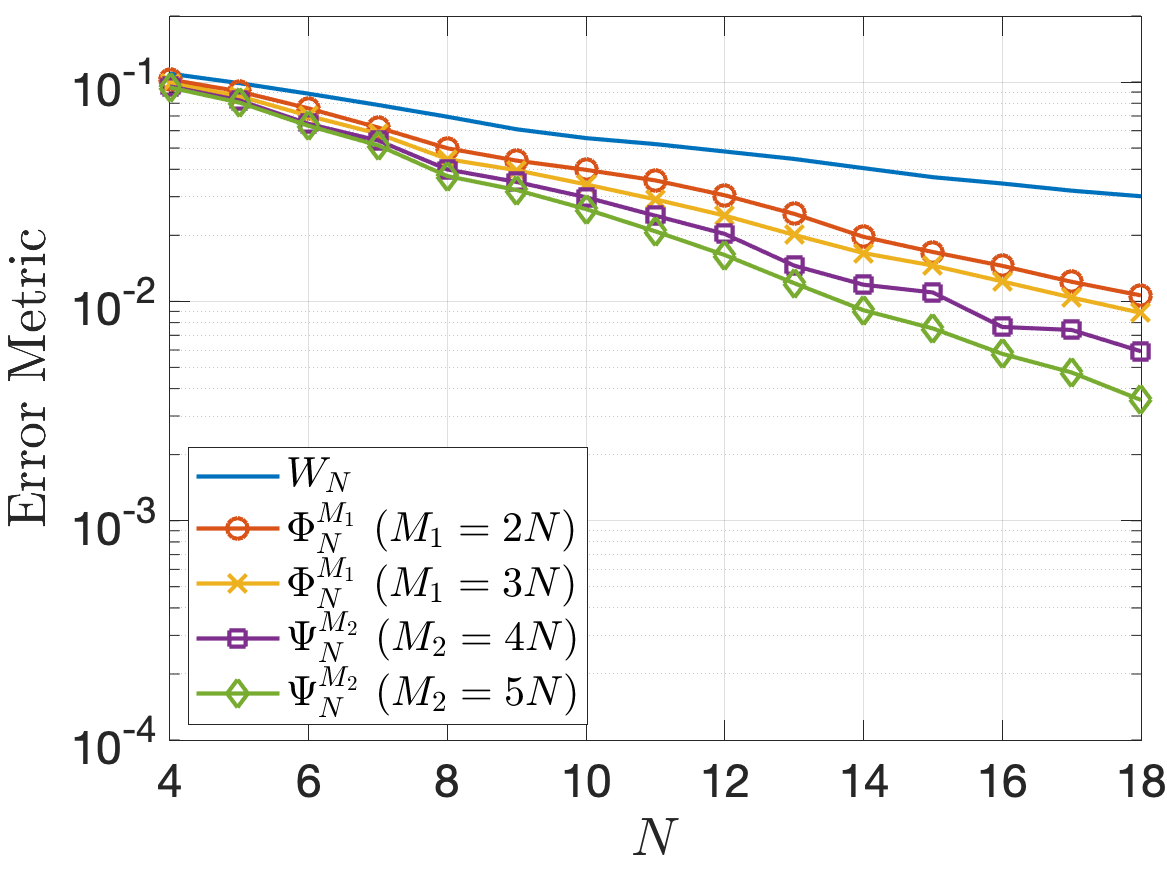}
		\caption{Quadratic function}
	\end{subfigure}
	\caption{The error metric for the standard RB space $W_N$ and the generative RB spaces $\Phi_N^{M_1}$ and $\Psi_N^{M_2}$ as a function of $N$ for $L= \min(4, N)$.}
	\label{ex1fig2}
\end{figure}

\begin{figure}[htbp]
	\centering
	\begin{subfigure}[b]{0.32\textwidth}
		\centering		\includegraphics[width=\textwidth]{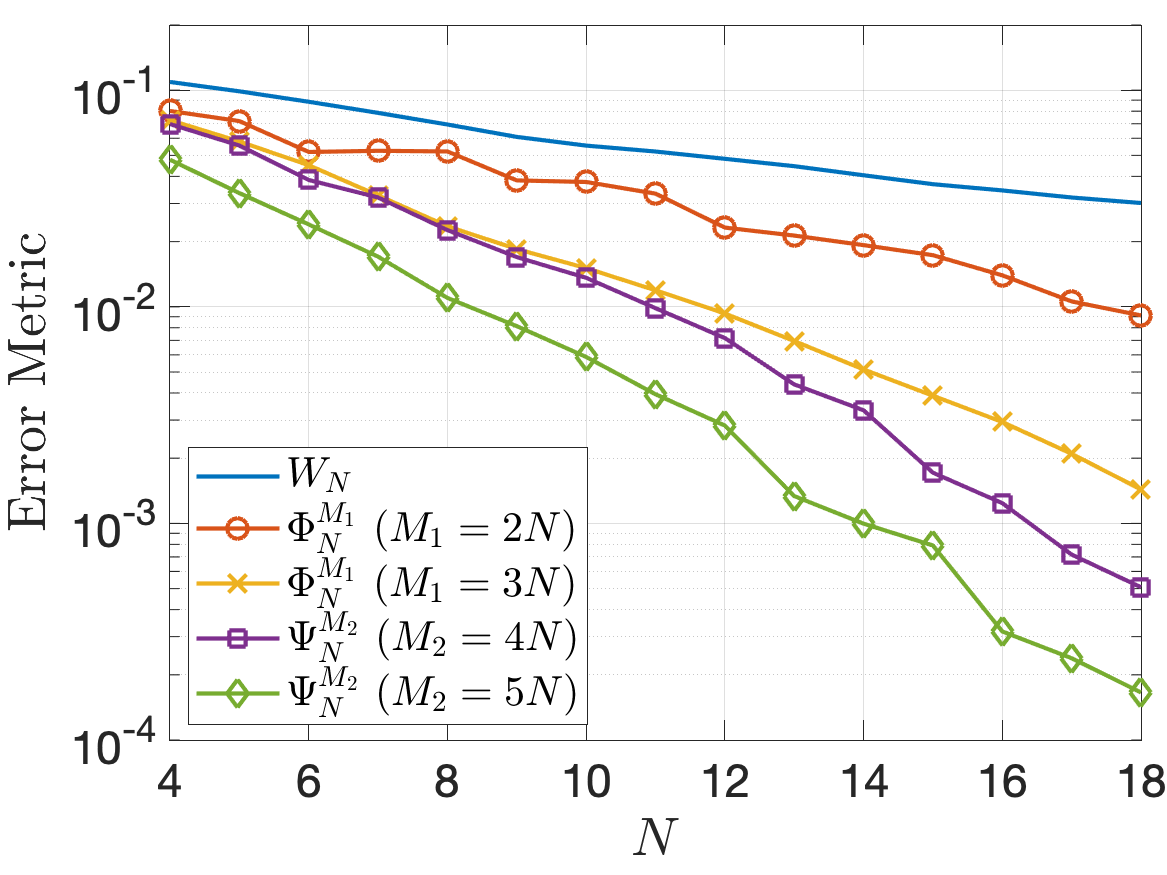}
		\caption{Hyperbolic tangent function}
	\end{subfigure}	
	\begin{subfigure}[b]{0.32\textwidth}
		\centering		\includegraphics[width=\textwidth]{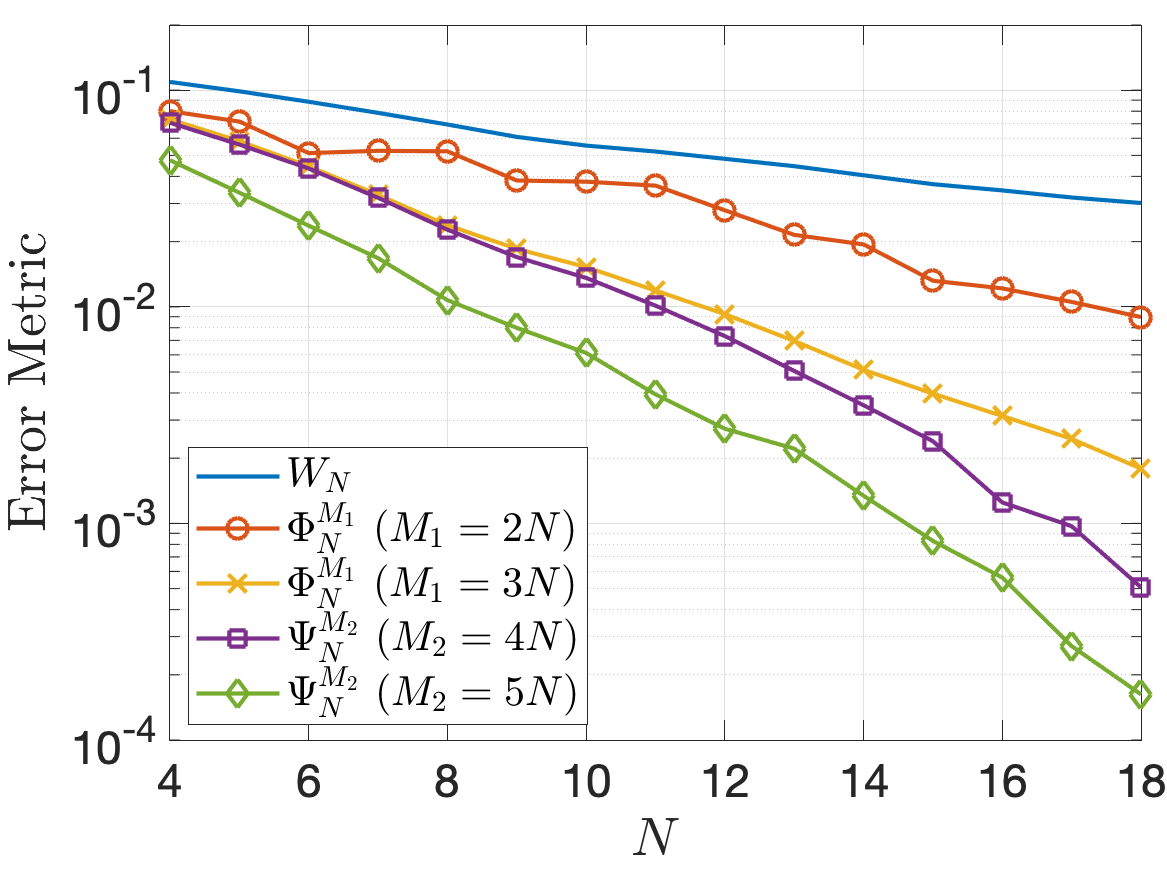}
		\caption{Sigmoid function}
	\end{subfigure}
 	\begin{subfigure}[b]{0.32\textwidth}
		\centering		\includegraphics[width=\textwidth]{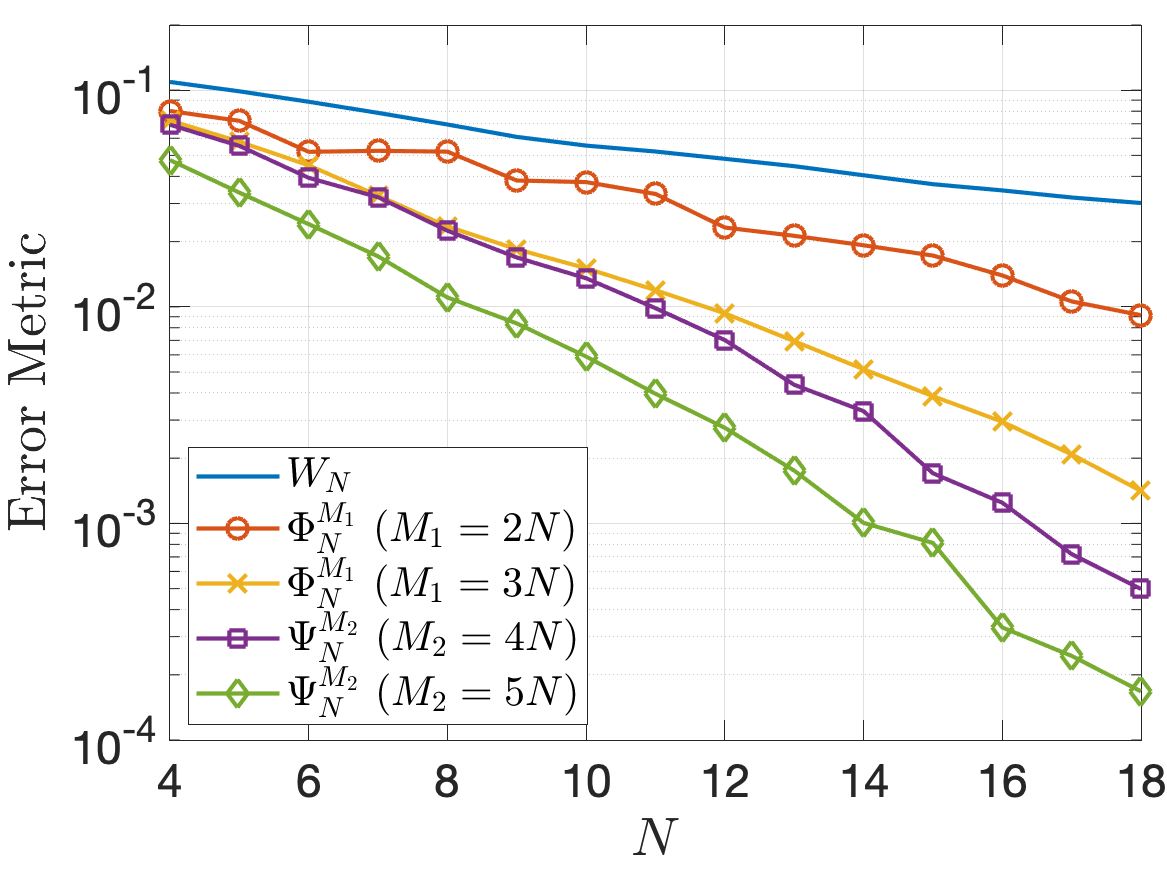}
		\caption{Arctangent function}
	\end{subfigure}	
	\begin{subfigure}[b]{0.32\textwidth}
		\centering		\includegraphics[width=\textwidth]{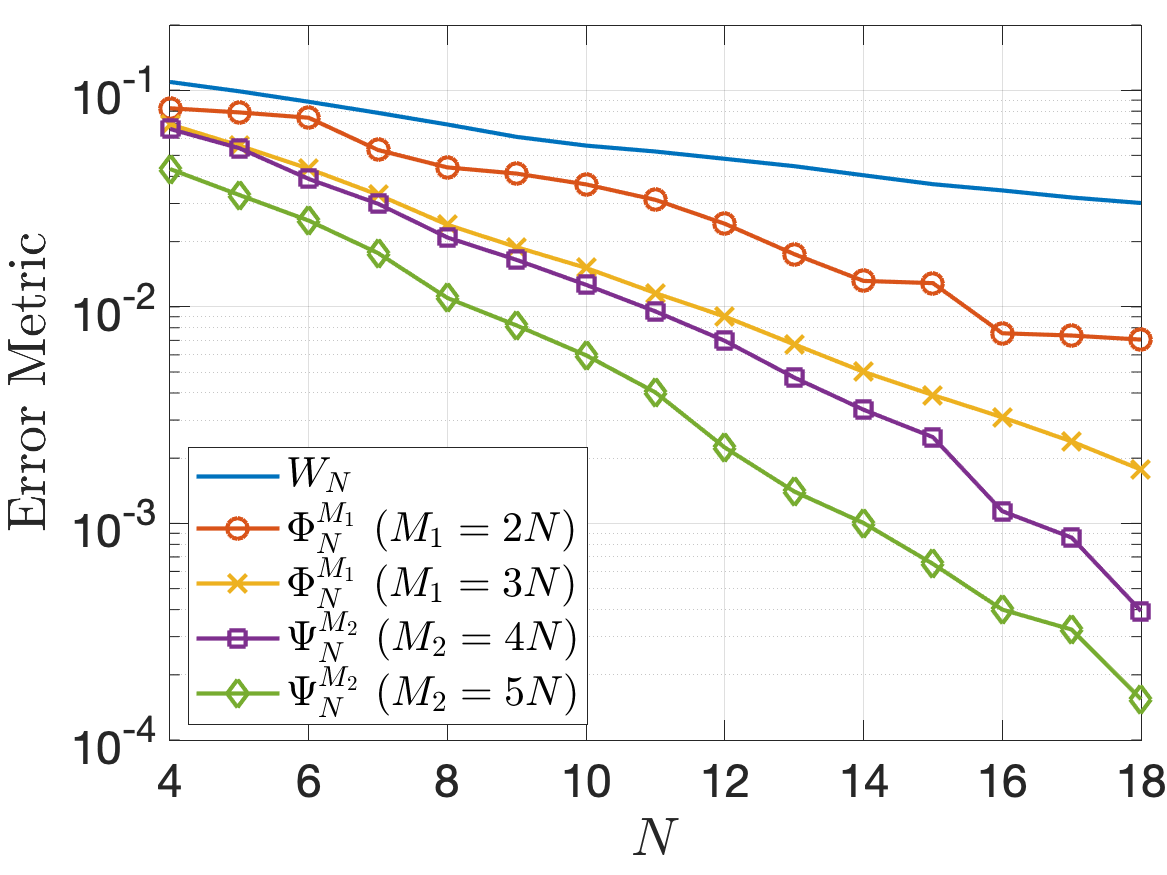}
		\caption{Softplus function}
	\end{subfigure}
  	\begin{subfigure}[b]{0.32\textwidth}
		\centering		\includegraphics[width=\textwidth]{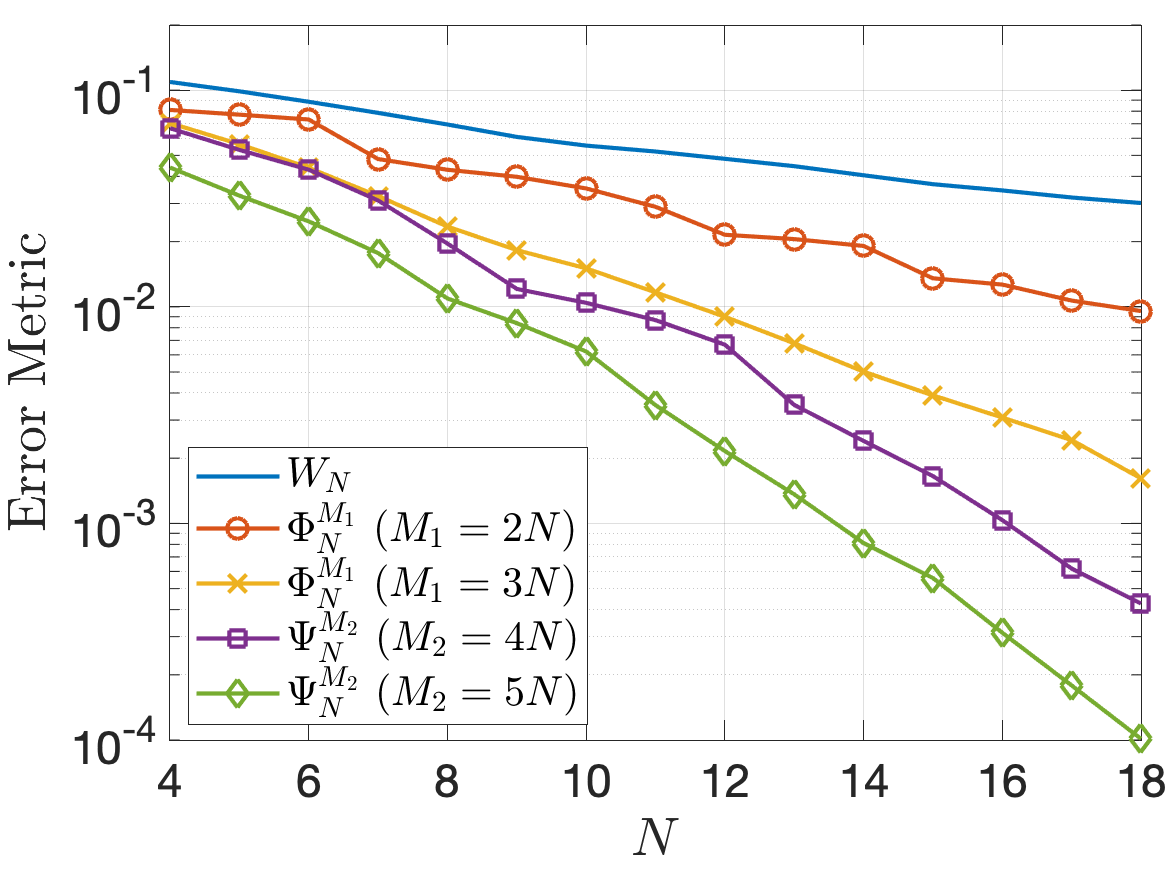}
		\caption{Exponential function}
	\end{subfigure}	
	\begin{subfigure}[b]{0.32\textwidth}
		\centering		\includegraphics[width=\textwidth]{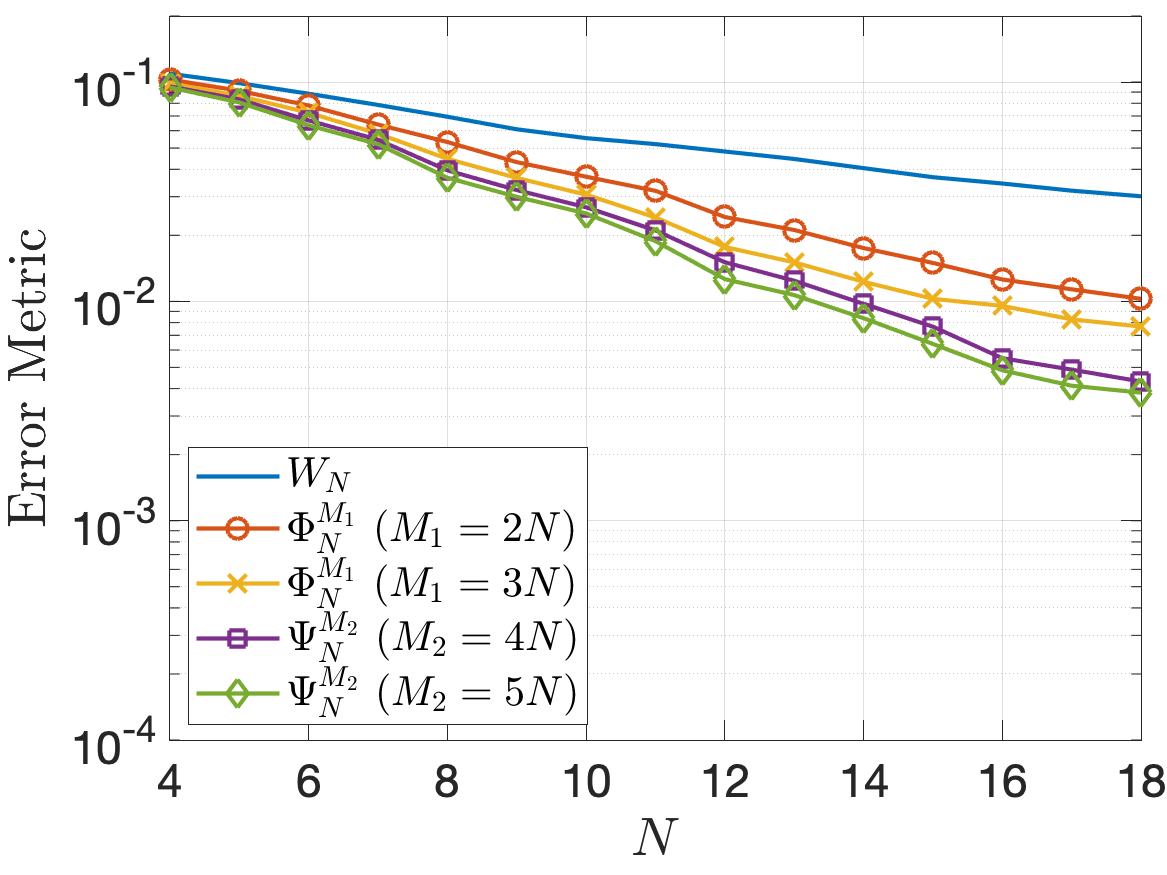}
		\caption{Quadratic function}
	\end{subfigure}
	\caption{The error metric for the standard RB space $W_N$ and the generative RB spaces $\Phi_N^{M_1}$ and $\Psi_N^{M_2}$ as a function of $N$  for $L= \min (8, N)$.}
	\label{ex1fig3}
\end{figure}

Figures \ref{ex1fig2} and \ref{ex1fig3} display the error metric for the standard and generative RB spaces as a function of $N$ for $L = \min(4,N)$ and $L = \min(8,N)$, respectively. As expected, increasing
$N$ reduces the error in all RB spaces. However, the generative RB spaces show significantly lower errors compared to the standard RB space for all choices of the nonlinear function, with the exception of the quadratic function. For the hyperbolic tangent, sigmoid, Softplus, and arctangent, and expotential functions, the error decreases more rapidly in the generative RB spaces as $N$ increases. This demonstrates the ability of these transformations to enrich the RB spaces. In contrast, the quadratic function does not offer the same degree of error reduction. Furthermore, increasing $L$ leads to reduction in the error for the space $\Phi_N^{M_1}$, but it has very little impact on the space $\Phi_N^{M_2}$. This indicates that the higher-order correlations are sufficient for improving accuracy, and further increasing the number of snapshots offers diminishing returns.


In summary, the results show that the generative RB spaces based on higher-order transformations (such as $\Psi_N^{M_2}$) consistently outperform both the standard RB space and lower-order generative spaces (such as $\Phi_N^{M_1}$). This performance gap becomes more pronounced as $N$ increases.  This highlights the importance of selecting appropriate nonlinear functions and transformation orders when constructing generative RB spaces.

\subsection{Two-Dimensional Parametrized Function}

\begin{figure}[h!]
	\centering
	\begin{subfigure}[b]{0.49\textwidth}
		\centering		\includegraphics[width=\textwidth]{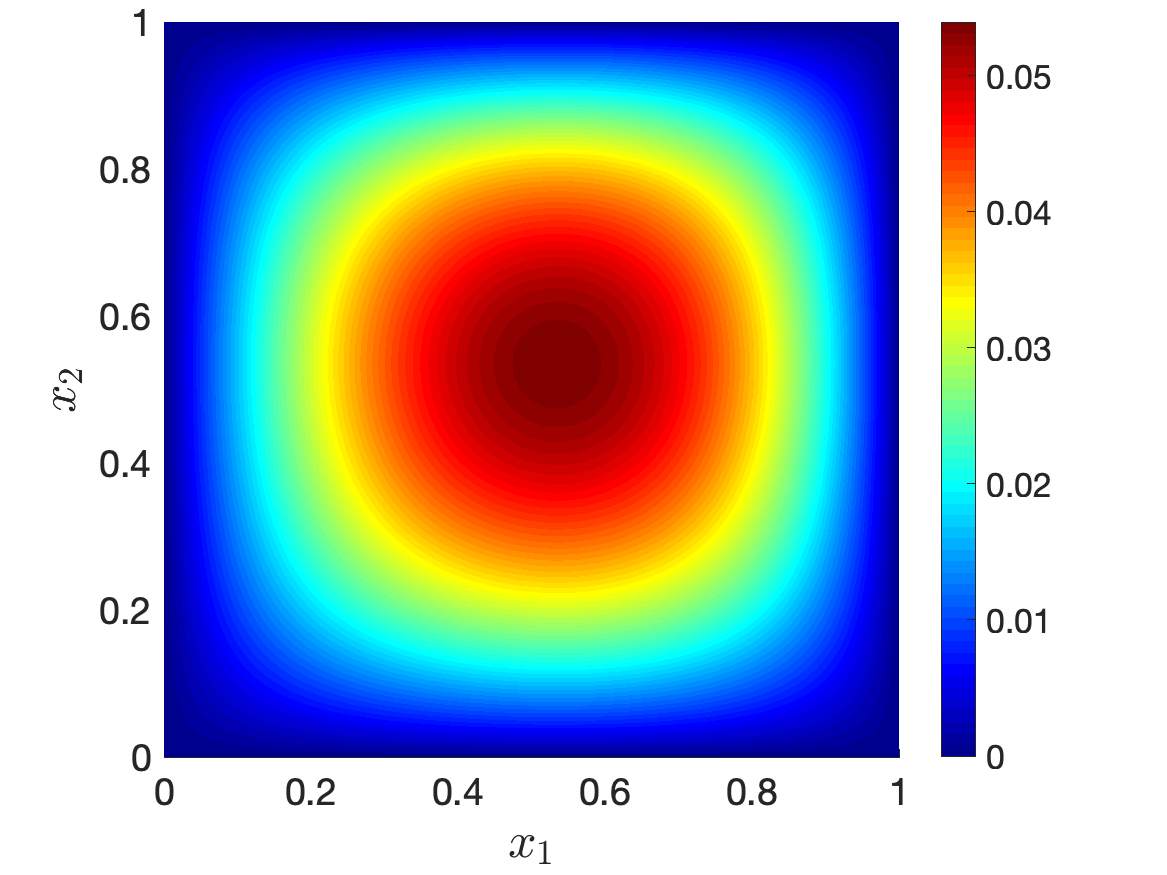}
		\caption{$\bm \mu = (1,1)$.}
	\end{subfigure}
	\hfill
	\begin{subfigure}[b]{0.49\textwidth}
		\centering		\includegraphics[width=\textwidth]{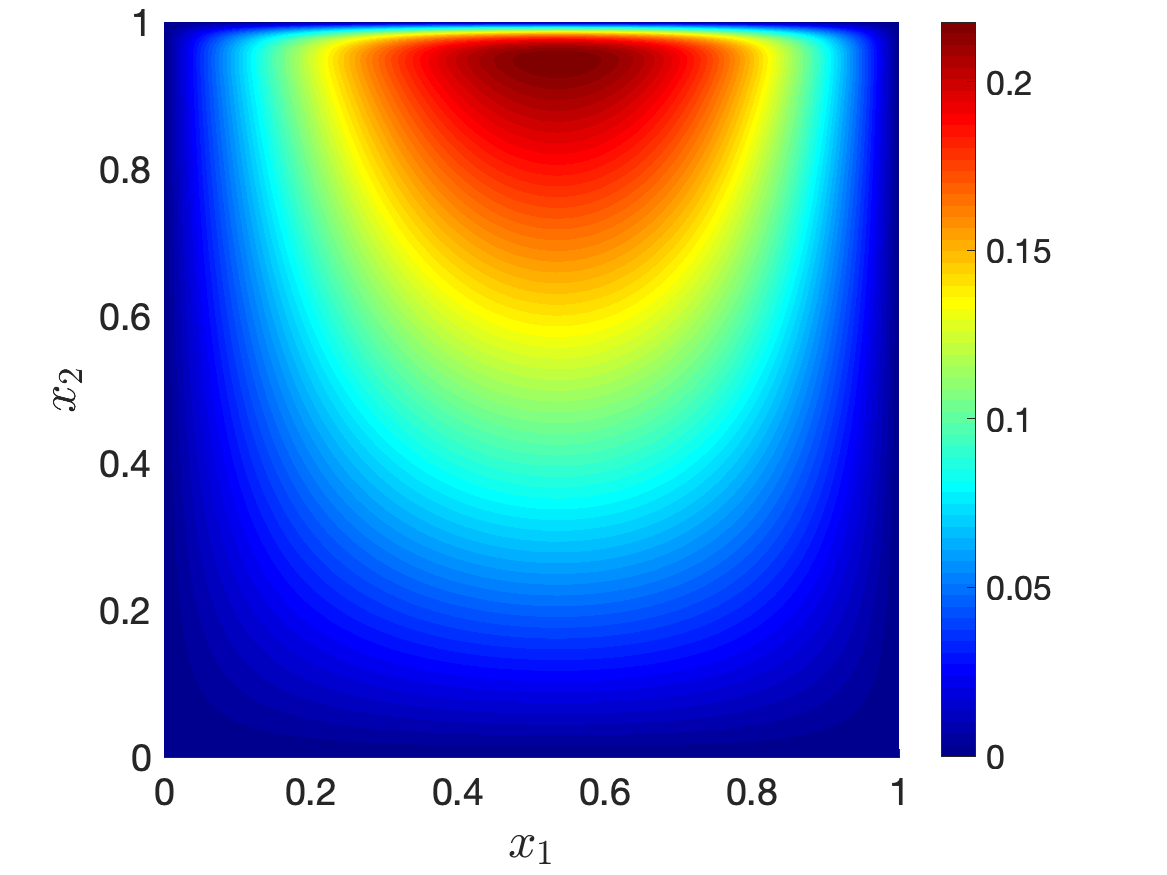}
		\caption{$\bm \mu = (50,1)$.}
	\end{subfigure} \\
        \begin{subfigure}[b]{0.49\textwidth}
		\centering		\includegraphics[width=\textwidth]{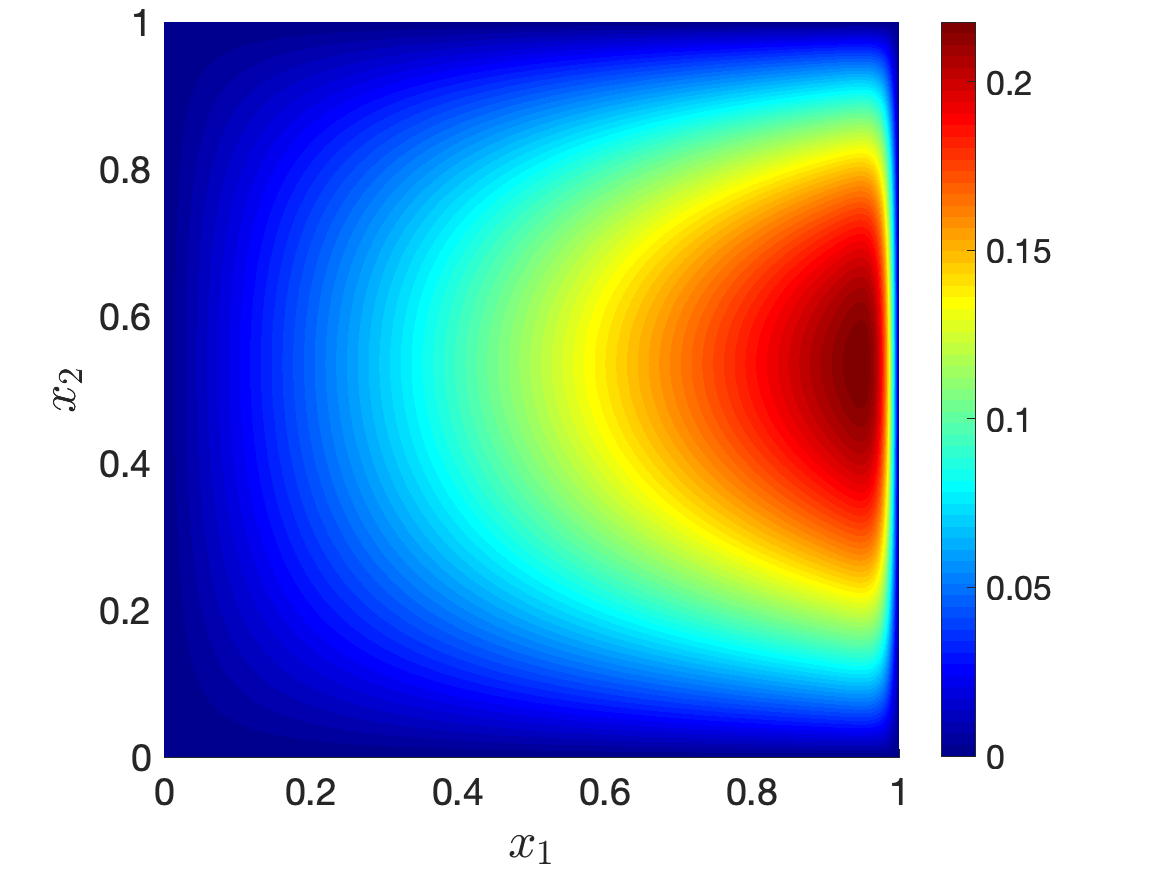}
		\caption{$\bm \mu = (1,50)$.}
	\end{subfigure}
	\hfill
	\begin{subfigure}[b]{0.49\textwidth}
		\centering		\includegraphics[width=\textwidth]{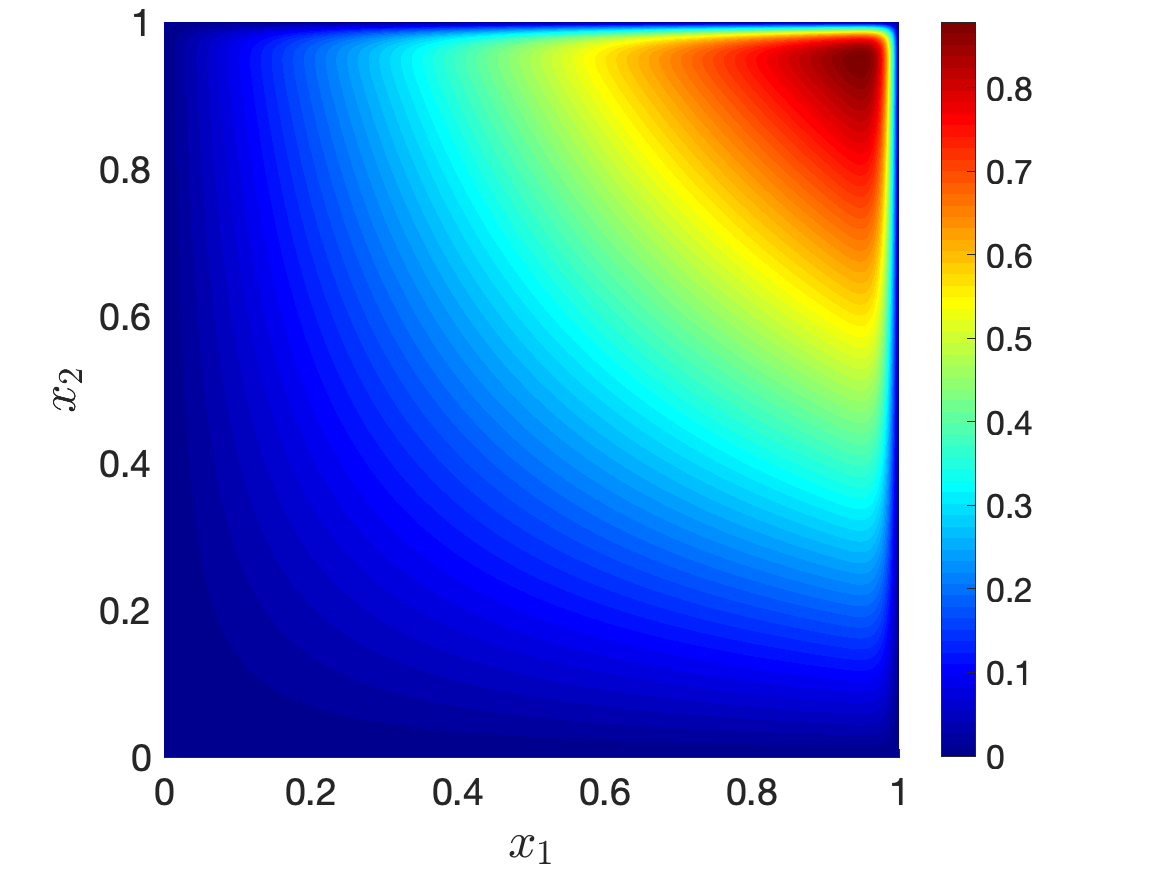}
		\caption{$\bm \mu = (50,50)$.}
	\end{subfigure}
	\caption{Plots of $u(\bm x, \bm \mu)$ defined in (\ref{ex2u}) for four different values of $\bm \mu$.}
	\label{ex2fig1}
\end{figure}

We present numerical results from a two-dimensional test case to assess the performance of the generative snapshot method. The test case involves the following parametrized function 
\begin{equation}
\label{ex2u}
u(\bm x,\bm \mu) = x_1 x_2 \tanh ((1-x_1) \mu_1) \tanh ((1-x_2) \mu_2)
\end{equation}
in a physical domain $\Omega = [0,1]^2$ and parameter domain $\mathcal{D} = [1, 50]^2$. We choose the function space $X = L^2(\Omega)$. Figure \ref{ex2fig1} shows four instances of $u(\bm x, \bm \mu)$ corresponding to the four corners of the parameter domain.  As $\bm \mu$ increases, the solution develops sharp boundary layers.

We choose for $S_N$ a set of $\sqrt{N} \times \sqrt{N}$ parameter points from the extended Chebyshev distribution over the parameter domain $\mathcal{D}$. The generative RB spaces, $\Phi_N^{M_1}$ and $\Psi_N^{M_2}$, are constructed from these $N$ solutions by using the generative snapshot method for  $L = \min(5,N)$. We employ $40 \times 40$ parameter points sampled uniformly over $\mathcal{D}$ to construct the discrete solution manifold $\mathcal{M}_K$ of $K = 1600$ snapshots. 

\begin{figure}[htbp]
	\centering
	\begin{subfigure}[b]{0.32\textwidth}
		\centering		\includegraphics[width=\textwidth]{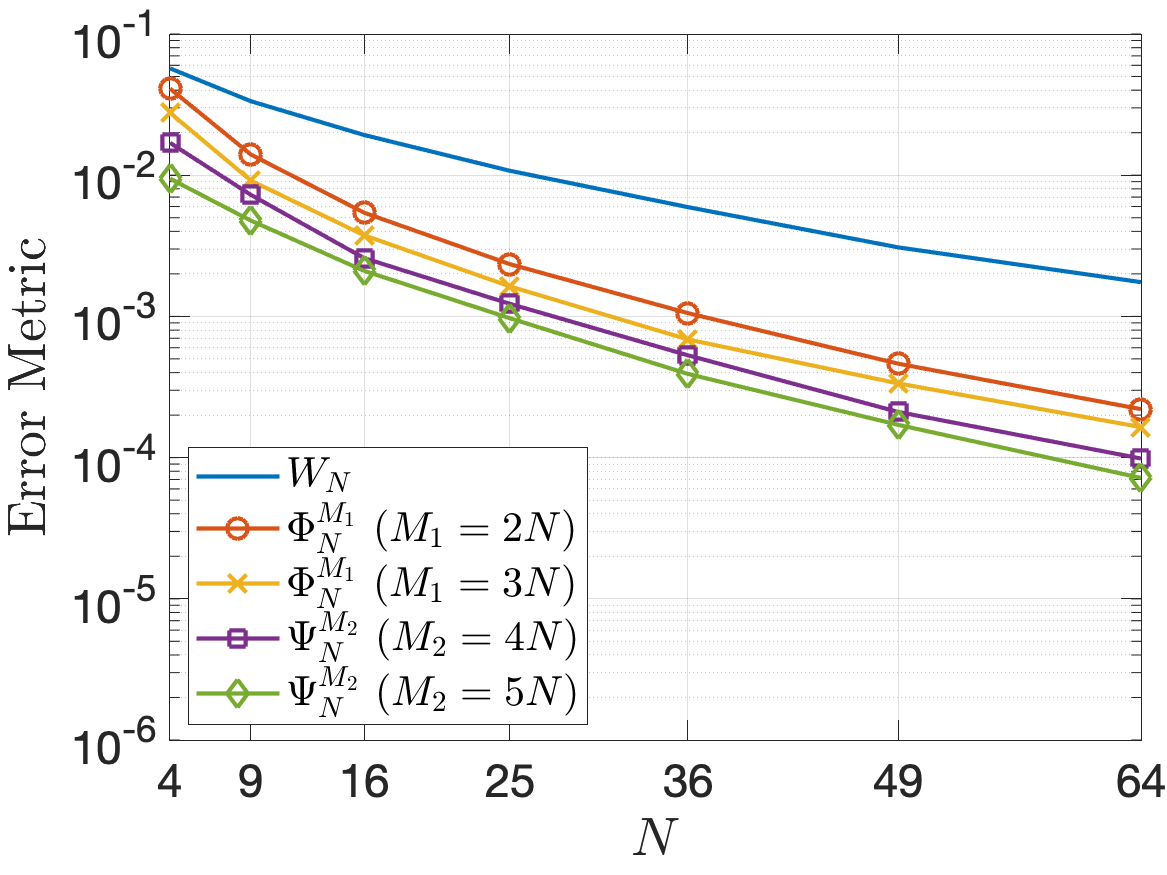}
		\caption{Hyperbolic tangent function}
	\end{subfigure}
	\begin{subfigure}[b]{0.32\textwidth}
		\centering		\includegraphics[width=\textwidth]{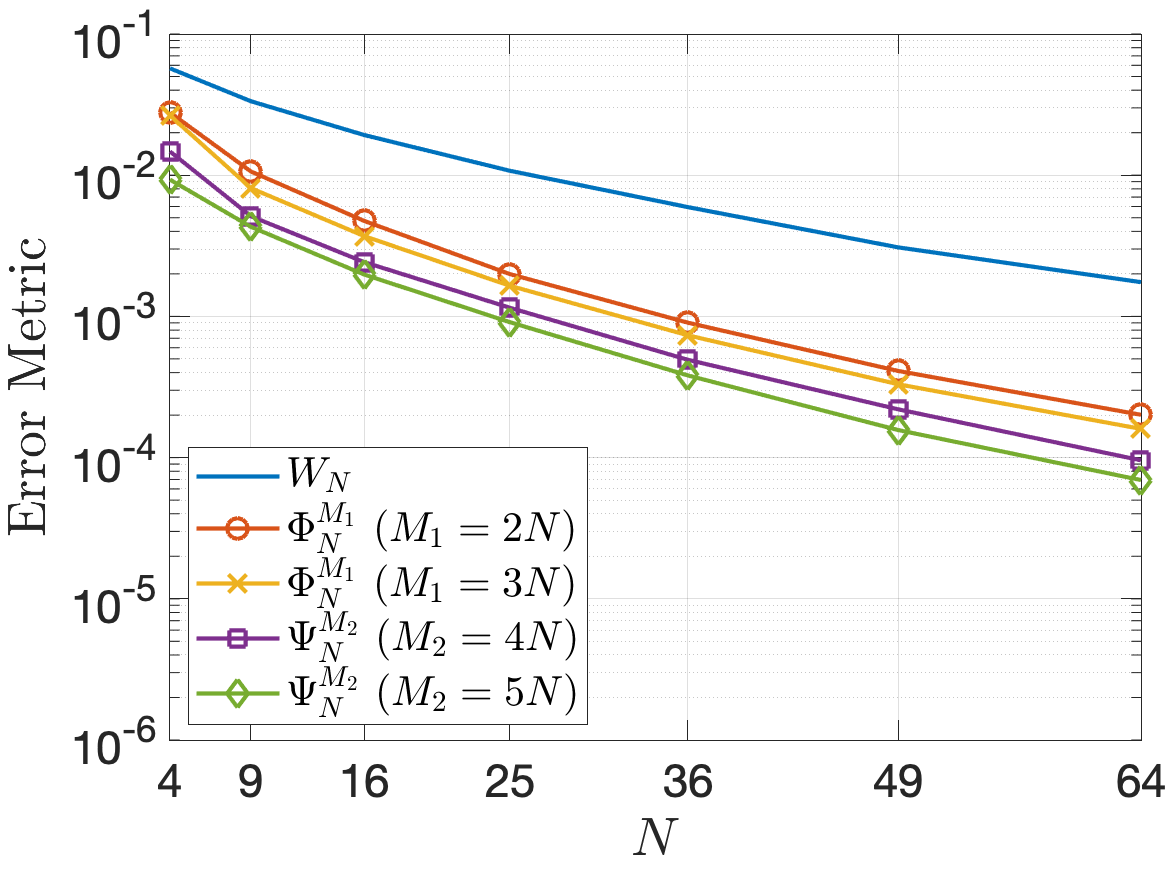}
		\caption{Sigmoid function}
	\end{subfigure}
 	\begin{subfigure}[b]{0.32\textwidth}
		\centering		\includegraphics[width=\textwidth]{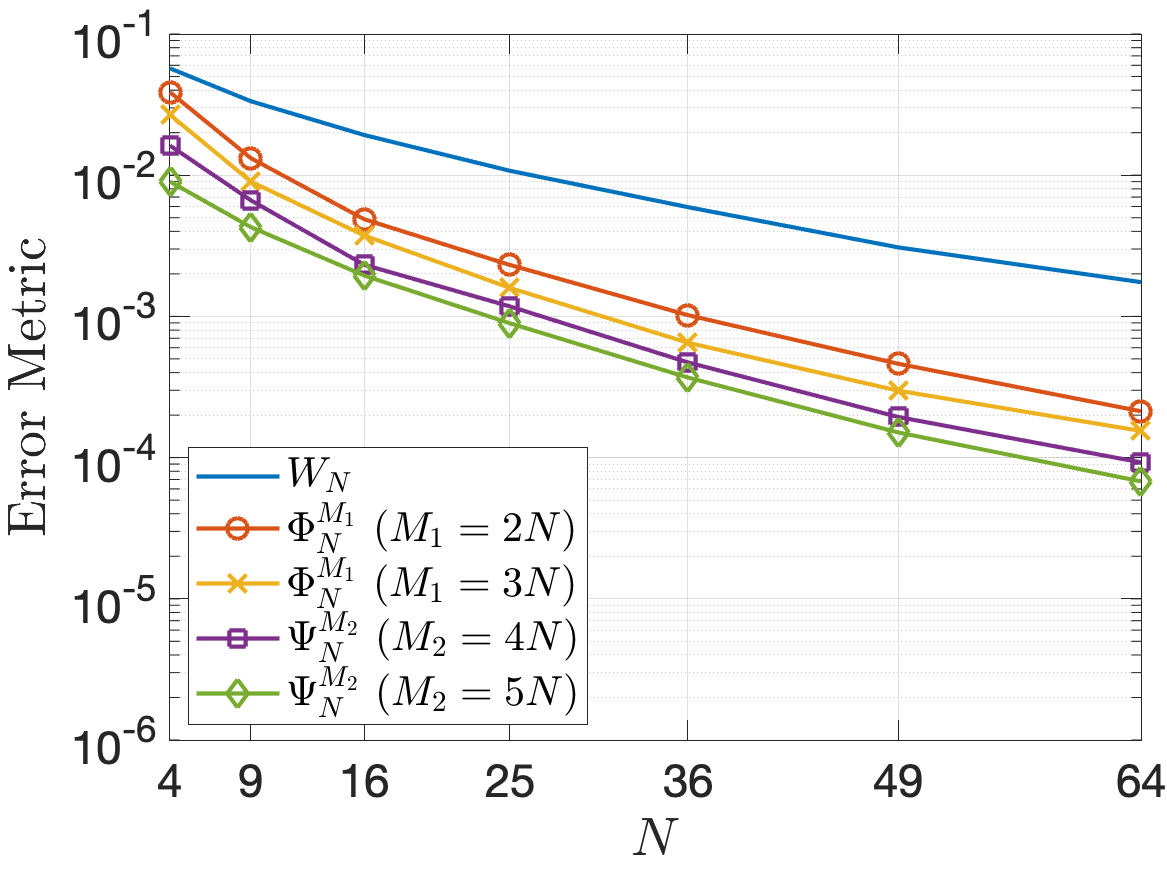}
		\caption{Arctangent function}
	\end{subfigure}
	\begin{subfigure}[b]{0.32\textwidth}
		\centering		\includegraphics[width=\textwidth]{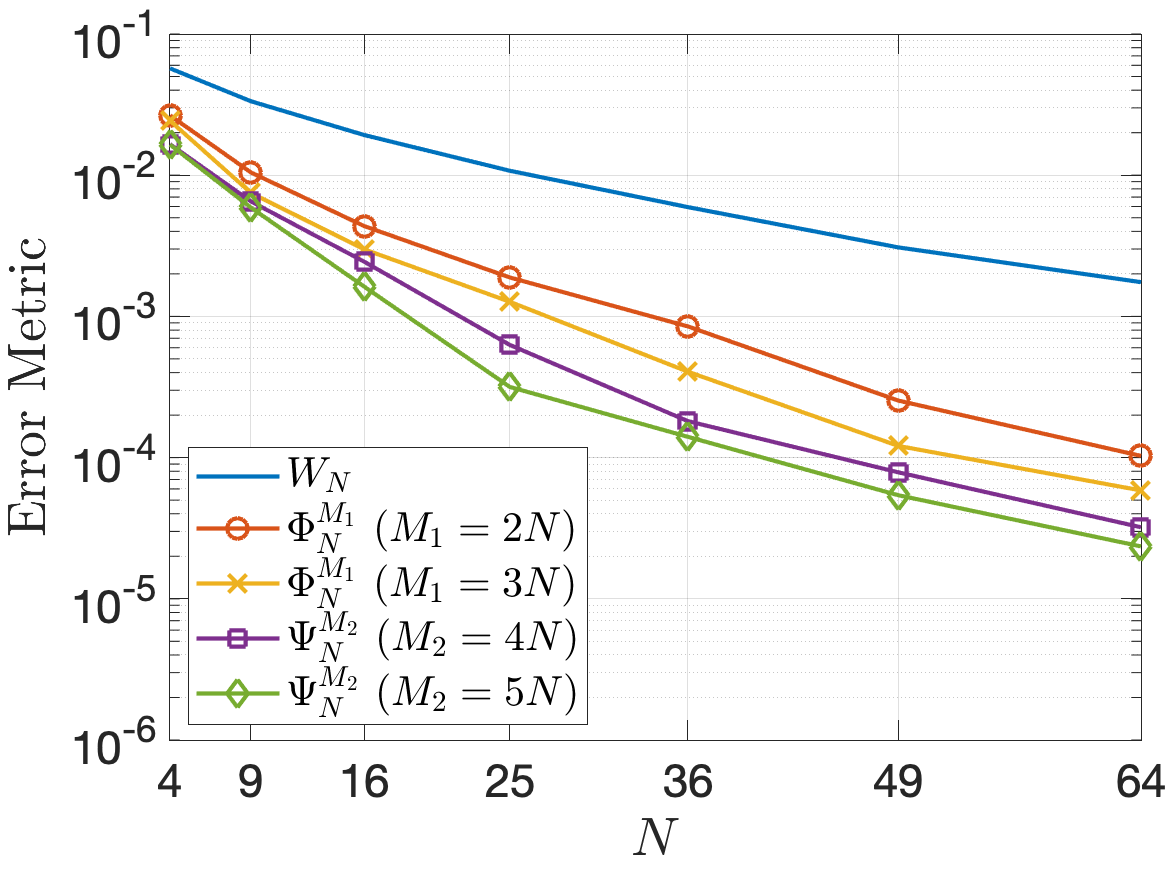}
		\caption{Softplus function}
	\end{subfigure}
  	\begin{subfigure}[b]{0.32\textwidth}
		\centering		\includegraphics[width=\textwidth]{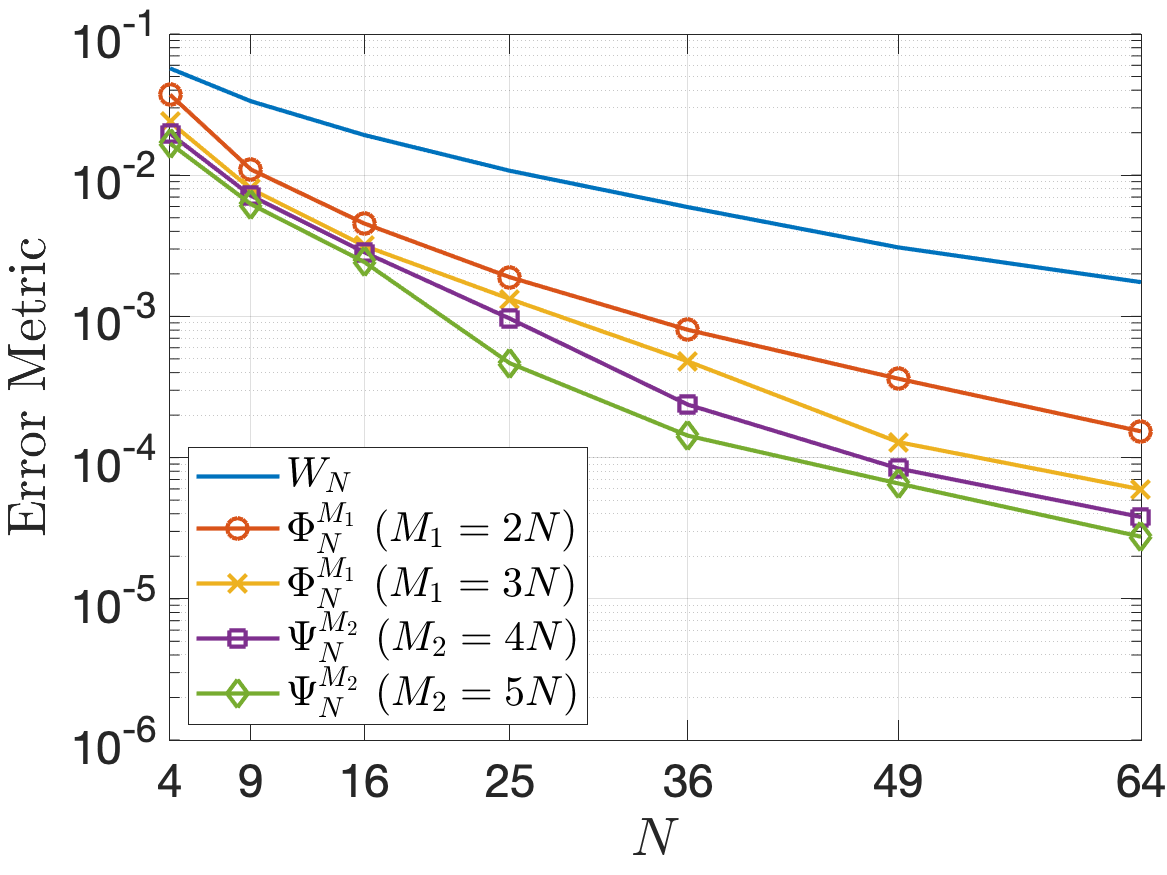}
		\caption{Exponential function}
	\end{subfigure}
	\begin{subfigure}[b]{0.32\textwidth}
		\centering		\includegraphics[width=\textwidth]{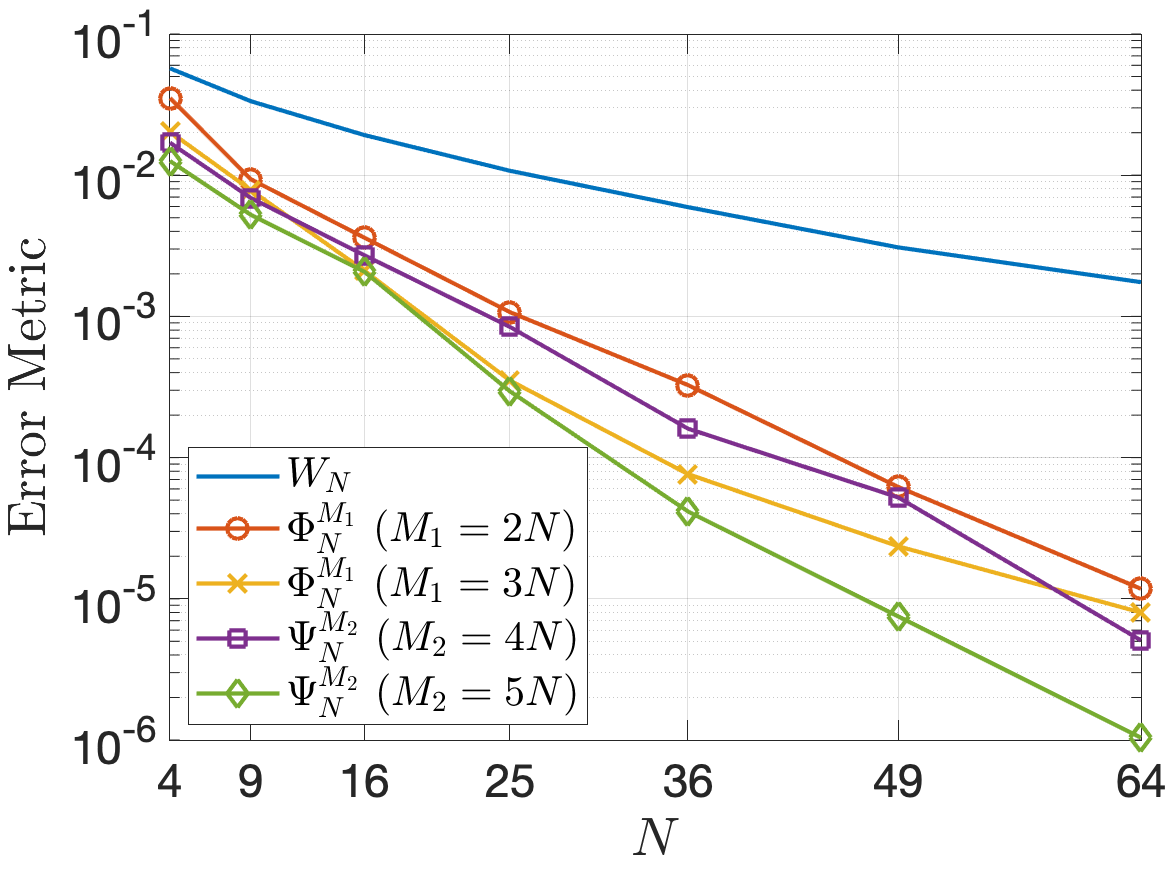}
		\caption{Quadratic function}
	\end{subfigure}
	\caption{The error metric for the standard RB space $W_N$ and the generative RB spaces $\Phi_N^{M_1}$ and $\Psi_N^{M_2}$ as a function of $N$ for $L= \min(5, N)$.}
	\label{ex2fig2}
\end{figure}

We assess the performance of the RB spaces by showing the convergence of their error metric as a function of $N$ in Figure \ref{ex2fig2}. As expected, increasing 
$N$ consistently reduces the error across all RB spaces. However, the generative RB spaces show significantly lower errors compared to the standard RB space for all choices of the nonlinear function. Among the tested nonlinear activation functions, the hyperbolic tangent, sigmoid, and arctangent show very similar convergence behaviors for the error metric. These functions show a steady reduction in error as $N$ increases, indicating their effectiveness in capturing non-linearity within the solution space. The softplus and exponential activation functions, on the other hand, result in even faster convergence, yielding lower errors compared to the hyperbolic tangent, sigmoid, and arctangent functions.

Interestingly, in contrast to the previous example, the quadratic function exhibits the fastest convergence and the smallest error in this case. This indicates that while the performance of the quadratic function may vary depending on the characteristics of the underlying solution manifold, it can be highly effective in specific contexts where polynomial-based transformations are sufficient to capture the relevant features. For $N=64$, the error for $\Psi_{N}^{M_2}$ with $M_2 = 5N$ is 1000 times smaller than the error for the standard RB space $W_N$. This highlights the exceptional performance of the generative RB space with quadratic transformations. The significant error reduction suggests that these higher-order transformations allow the generative RB space to capture the complex structure of the solution manifold more accurately.

Overall, these results confirm the advantages of using nonlinear transformations to expand the snapshot set and construct generative RB spaces. The ability of  softplus and exponential functions to accelerate error convergence, combined with the strong performance of the quadratic function in this particular case, demonstrates the flexibility of the generative approach. The results also underscore the importance of selecting the appropriate non-linear transformation based on the specific characteristics of the problem to be addressed, as the choice of $\sigma$ can significantly influence both the convergence rate and the accuracy of the RB spaces.



\subsection{Three-Dimensional Parametrized Function}

In the last example, we consider the  parametrized spherical Bessel function \cite{Nguyen2023,Nguyen2023b,Nguyen2024b,nguyen2024c}: 
\begin{equation}
\label{ex3u}
u(\bm x, \mu)  =  \frac{\sin (\mu \pi r) }{\mu \pi r} \exp \left(1 -\frac{1}{\sqrt{\left(1 - r^3 \right)^2 + 10^{-6}}} \right), 
\end{equation}
where $r = \sqrt{x_1^2 +x_2^2 +x_3^2}$ is the radial distance from the origin in the unit sphere domain $\Omega$, and $\mu$ is the parameter defined in the range $\mathcal{D} = [1, 20]$. Spherical Bessel functions are commonly used in modeling wave propagation, oscillatory phenomena, and scattering problems, making this an ideal test case to evaluate the performance of RB spaces in handling highly oscillatory solutions.

\begin{figure}[htbp]
	\centering
 \includegraphics[width=0.8\textwidth]{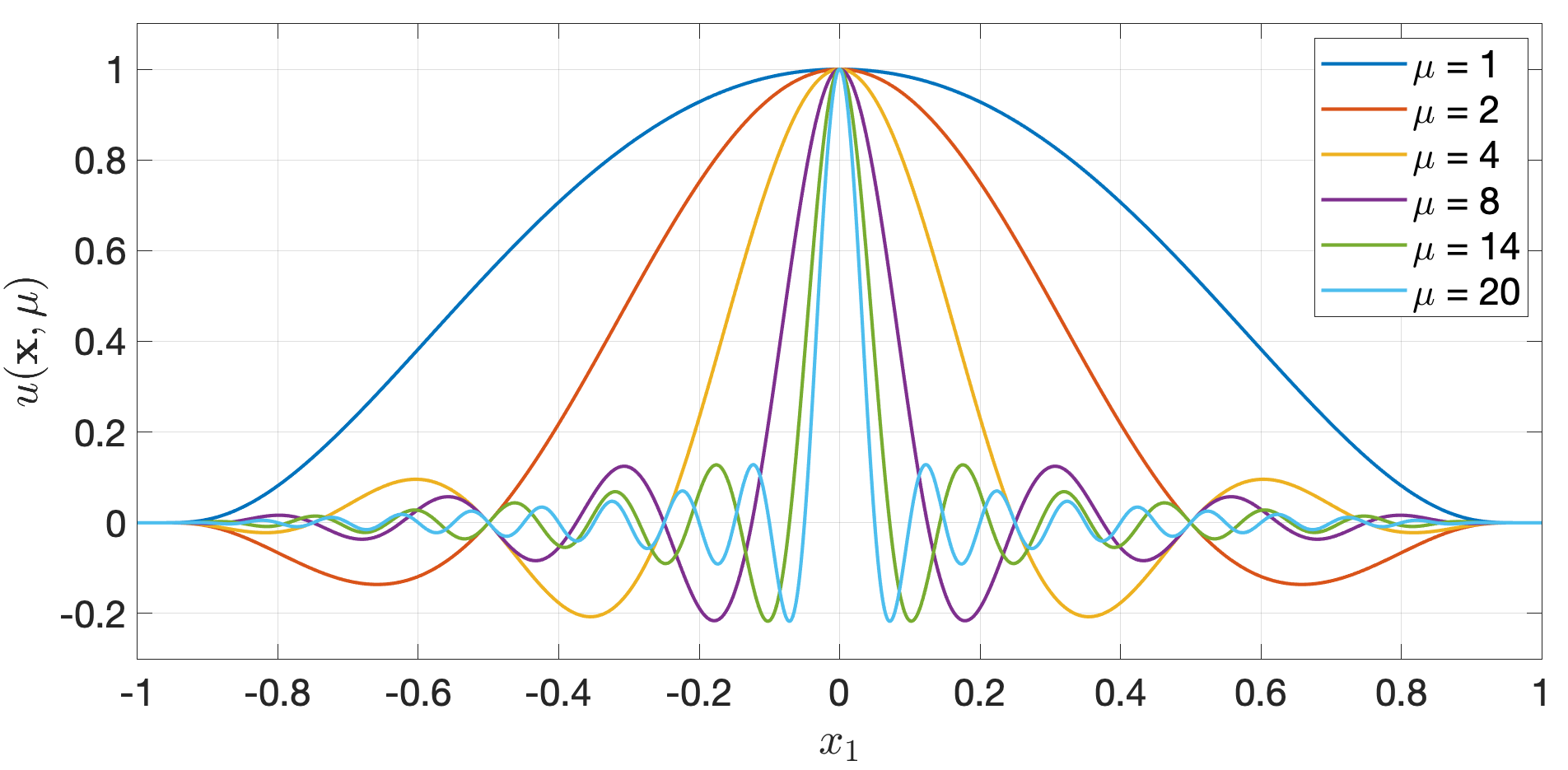}
 \caption{Plots of $u(\bm x,\mu)$  in (\ref{ex3u}) along the line $x_2 = x_3 = 0$ for  different values of $\mu$.}
	\label{ex3fig1}
\end{figure}

Figure \ref{ex3fig1} shows the behavior of  $u(\bm x,  \mu)$ along the line $x_2 = x_3 = 0$ for several values of $\mu$. As expected, the solution exhibits strong oscillations that depend on the parameter $\mu$. For lower values of $\mu$, the function displays more moderate oscillations, while for higher values, the frequency of oscillation increases significantly. This high-frequency behavior poses a challenge for the standard RB space, as accurately capturing such oscillations may require a large number of basis functions.

The parameter sample set $S_N$ is drawn from an extended Chebyshev distribution over the parameter domain. The extended Chebyshev distribution is chosen to ensure better coverage of the parameter space, particularly in capturing regions where the solution manifold may exhibit more complex behavior. To define the discrete solution manifold $\mathcal{M}_K$, we employ $K=100$ parameter points, which are sampled uniformly over the parameter domain. By discretizing the solution manifold in this way, we can assess the performance of the RB spaces to approximate the full solution space effectively. The combination of the extended Chebyshev distribution for $S_N$ and the uniform sampling for $\mathcal{M}_K$ ensures that both the training set  and the test set are  well distributed across the parameter domain.

\begin{figure}[htbp]
	\centering
	\begin{subfigure}[b]{0.32\textwidth}
		\centering		\includegraphics[width=\textwidth]{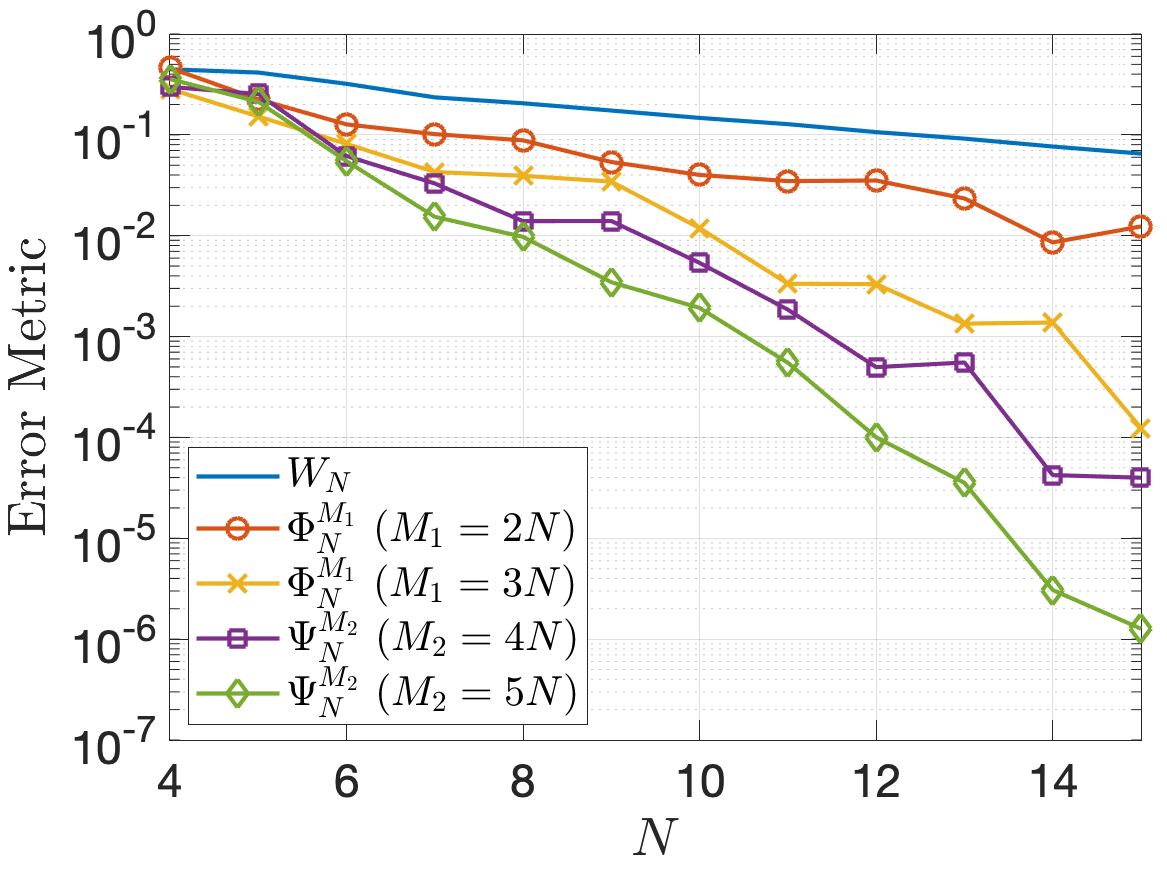}
		\caption{Hyperbolic tangent function}
	\end{subfigure}
	\begin{subfigure}[b]{0.32\textwidth}
		\centering		\includegraphics[width=\textwidth]{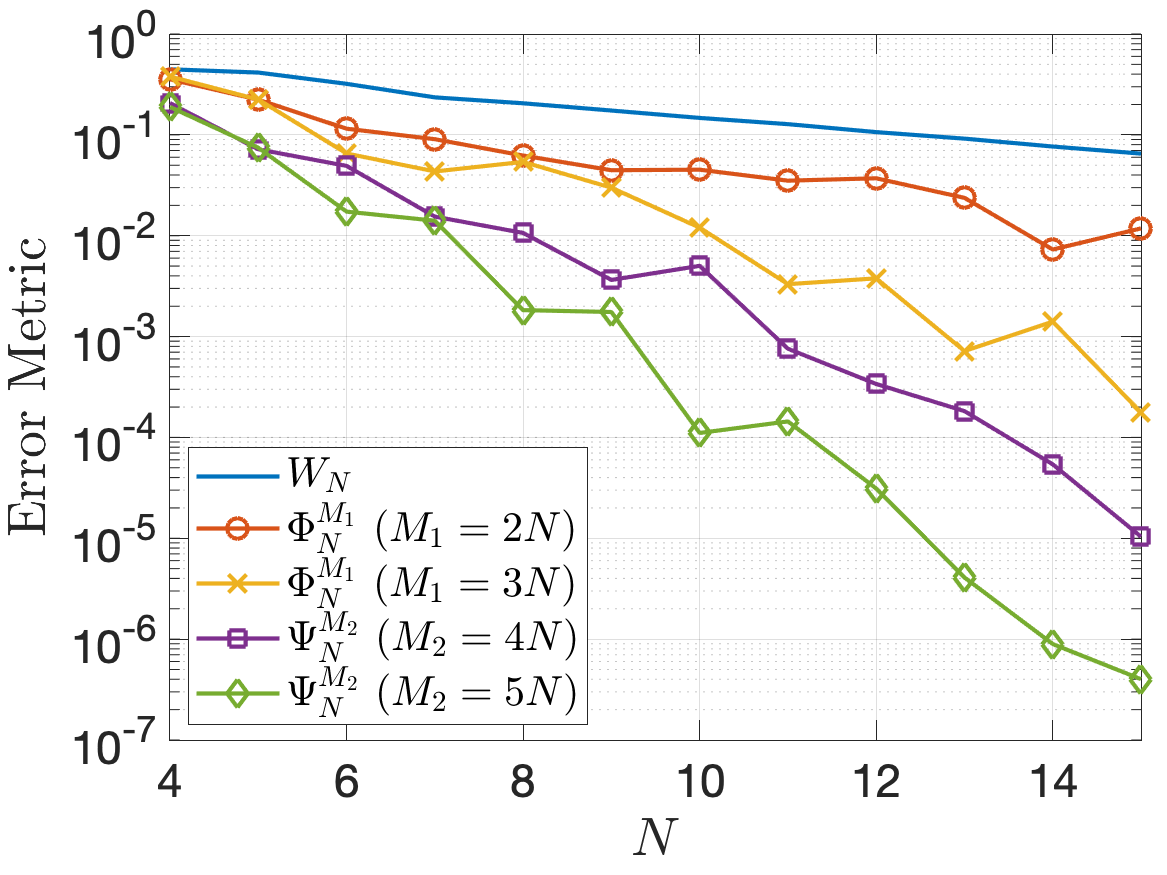}
		\caption{Sigmoid function}
	\end{subfigure}
 	\begin{subfigure}[b]{0.32\textwidth}
		\centering		\includegraphics[width=\textwidth]{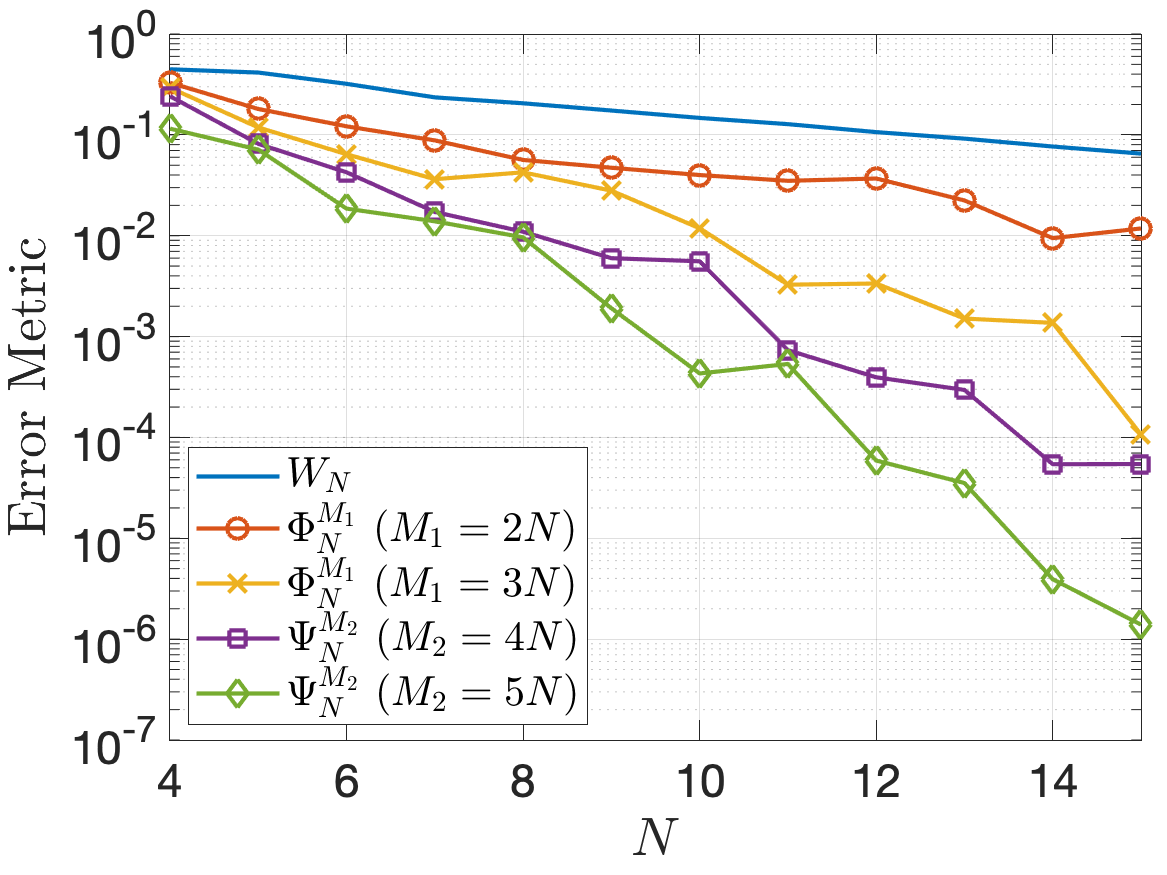}
		\caption{Arctangent function}
	\end{subfigure}	
	\begin{subfigure}[b]{0.32\textwidth}
		\centering		\includegraphics[width=\textwidth]{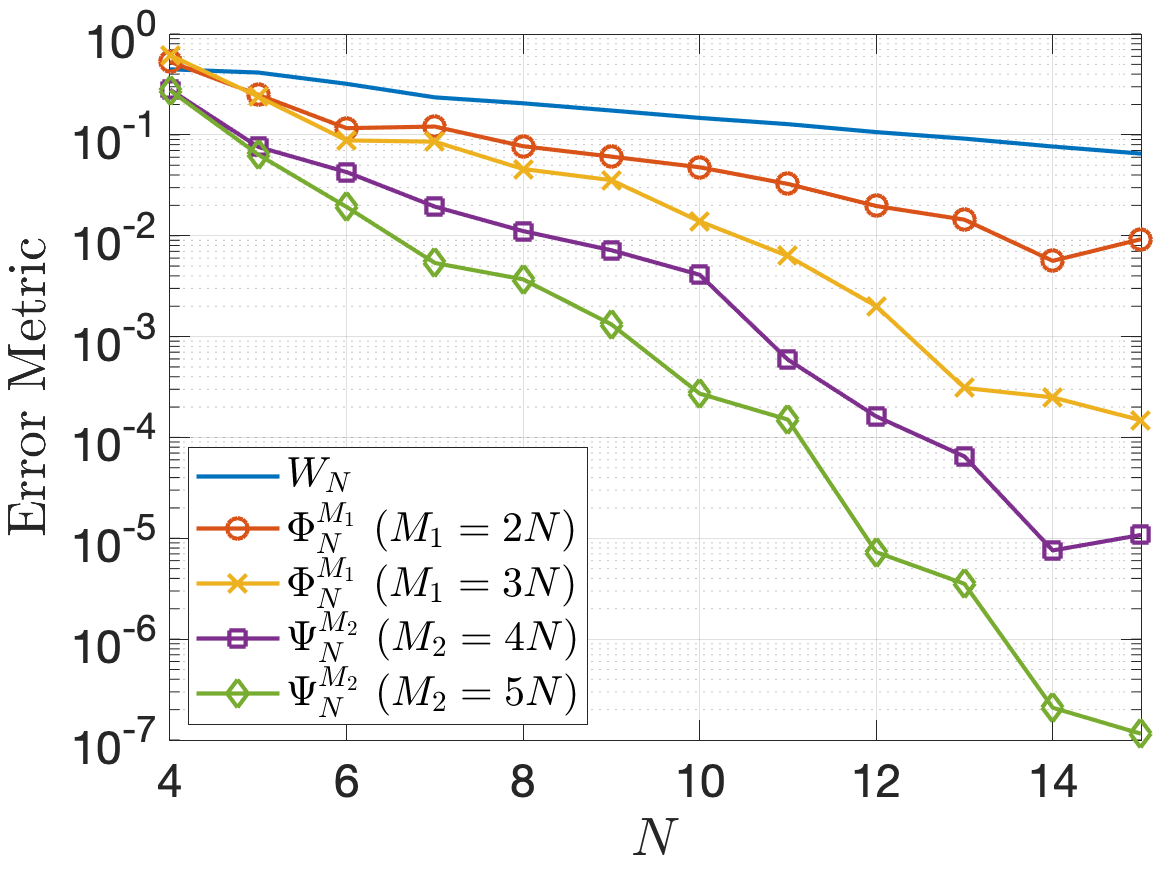}
		\caption{Softplus function}
	\end{subfigure}
  	\begin{subfigure}[b]{0.32\textwidth}
		\centering		\includegraphics[width=\textwidth]{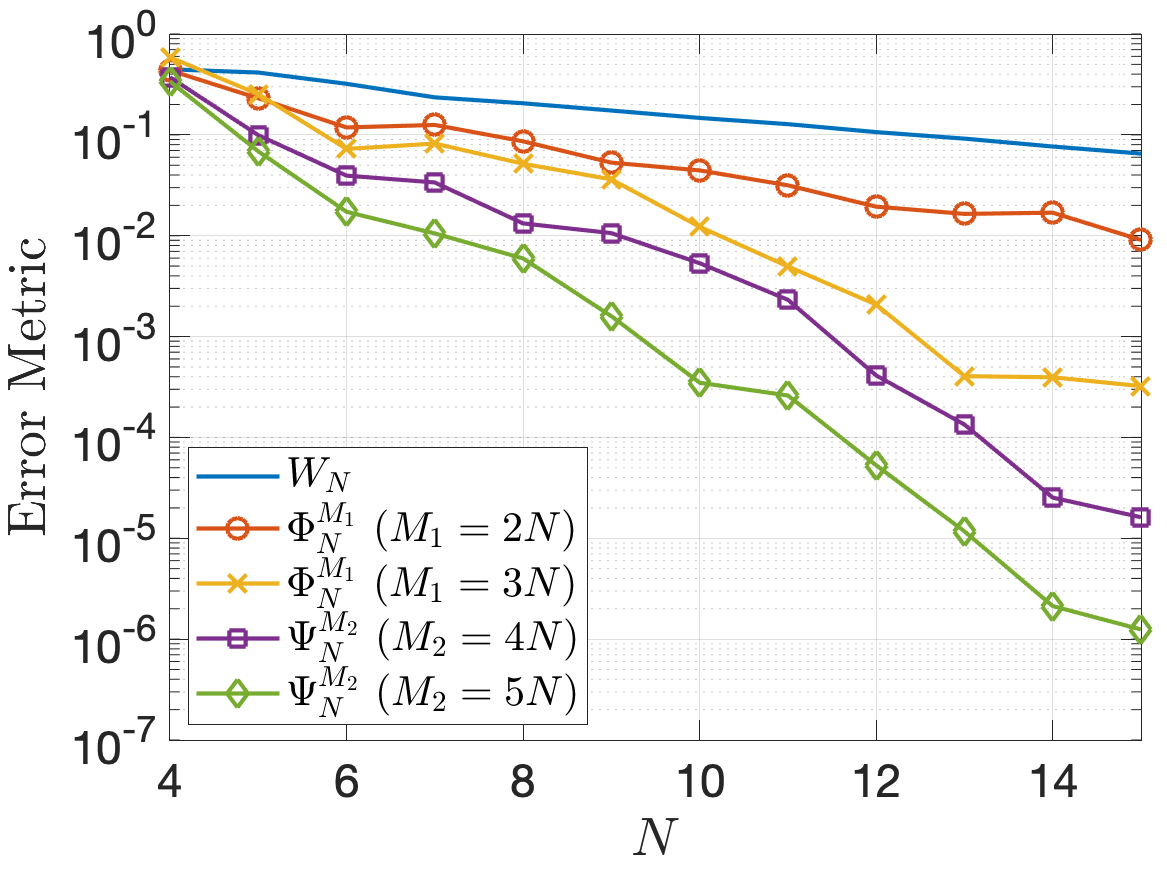}
		\caption{Exponential function}
	\end{subfigure}	
	\begin{subfigure}[b]{0.32\textwidth}
		\centering		\includegraphics[width=\textwidth]{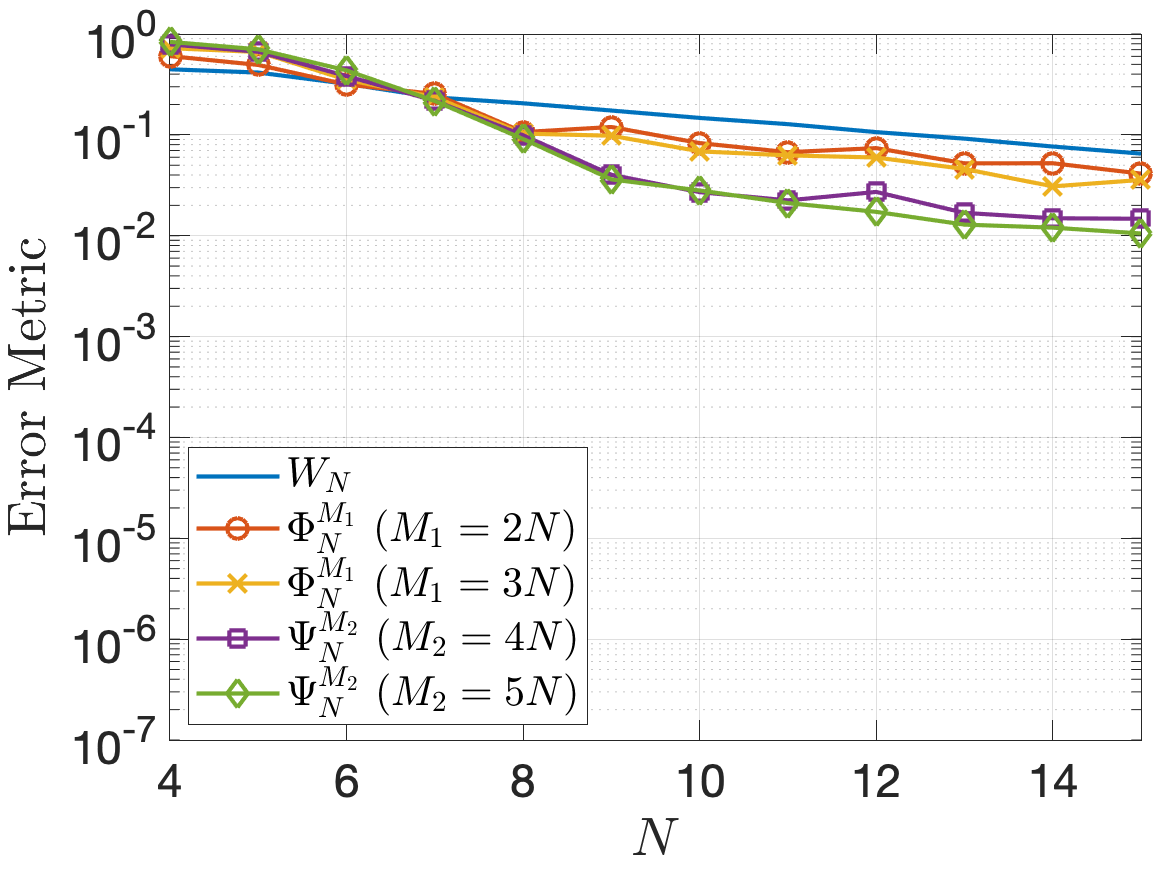}
		\caption{Quadratic function}
	\end{subfigure}
	\caption{The error metric for the standard RB space $W_N$ and the generative RB spaces $\Phi_N^{M_1}$ and $\Psi_N^{M_2}$ as a function of $N$ for $L= \min(5, N)$.}
	\label{ex3fig2}
\end{figure}

Figure \ref{ex3fig2} presents the error metric as a function of the number of basis functions $N$ for the standard and generative RB spaces. Across all subplots, the standard RB space consistently exhibits the highest error for all values of $N$, ranging from 0.45 to 0.065 as $N$ increases from 4 to 14. This slow error reduction reflects the difficulty in capturing the highly oscillatory behavior of the spherical Bessel function. The error decreases with increasing $N$, but the rate of reduction is slow, highlighting the inefficiency of the standard RB space when dealing with complex oscillatory functions. The generative RB spaces significantly outperform the standard RB space, with much lower errors across the range of $N$. For $N=14$, the error for $\Psi_N^{M_2}$ with $M_2 = 5N$ reaches values below $10^{-6}$ representing an improvement of four orders of magnitude compared to the standard RB space. This illustrates the effectiveness of the generative snapshot method in handling oscillatory solutions.

The softplus and exponential functions show the most rapid error convergence among the nonlinear functions used in the generative RB spaces.  They are highly effective at capturing the oscillatory behavior of the spherical Bessel function, making them particularly suited for complex problems with sharp transitions or high-frequency components. The hyperbolic tangent sigmoid, and arctangent functions also perform well and provide substantial error reduction compared to the standard RB space. In contrast, the quadratic function shows the slowest error reduction among these nonlinear functions. This suggests that the quadratic function is insufficient to capture the complexity of the spherical Bessel function.

In summary, the results demonstrate the clear superiority of generative RB spaces over the standard RB space when applied to the spherical Bessel function, which poses a challenging test case due to its oscillatory nature. Among the nonlinear functions considered, the softplus and exponential transformations offer the best performance, achieving rapid error reduction and high accuracy. The $\Psi_N^{M_2}$ spaces consistently outperform the $\Phi_N^{M_1}$ spaces. Therefore, higher-order transformations lead to richer snapshot sets that can more accurately represent the oscillatory behavior of the spherical Bessel function, resulting in faster convergence and lower error.

\section{Generative Reduced Basis Method}

In this section, we introduce a generative reduced basis (RB) method to construct reduced-order models (ROMs) for affine linear PDEs. The method leverages the generative snapshot technique to construct enriched RB spaces. The generative RB method not only improves accuracy, but also provides {\em a posterior} error estimates. These error estimates are used in a greedy algorithm to select parameter points to construct the generative RB spaces. 


\subsection{Finite Element Approximation}

The weak formulation for a parametrized linear PDE can be stated
as follows: given any $\bm \mu \in {\cal D} \subset \mathbb{R}^P$, we
evaluate $s^{\rm e}(\bm \mu) = \ell^O(u^{\rm e}(\bm \mu);\bm \mu)$, where $u^{\rm e}(\bm \mu) \in
X^{\rm e}$ is the solution of
\begin{equation}
a(u^{\rm e}(\bm \mu),v;\bm \mu) = \ell(v; \bm \mu), \quad \forall v \in X^{\rm e}.  
\label{eq:1}
\end{equation}
Here $\mathcal{D}$ is the parameter domain in which our $P$-tuple (input)
parameter $\bm \mu$ resides; $X^{\rm e}(\Omega)$ is an appropriate Hilbert space on $\Omega$;
$\Omega$ is a bounded domain in $\mathbb{R}^d$ with Lipschitz continuous boundary $\partial \Omega$; $\ell(\cdot;\bm \mu), \ \ell^O(\cdot;\bm \mu)$ are $X^{\rm
  e}$-continuous linear functionals;  and $a(\cdot,\cdot;\bm \mu)$ is a bilinear form of the parameterized PDE operator. In actual practice, we replace $u^{\rm e}(\bm \mu)$ with an approximate solution, $u(\bm \mu)$,
which resides in a finite element
approximation space $X \subset X^{\rm e}$ of {\em very} large dimension $\mathcal{N}$ and satisfies
\begin{equation}
a(u(\bm \mu),v;\bm \mu) = \ell(v), \quad \forall v \in X.  
\label{eq:1a}
\end{equation}
We then evaluate the finite element output as
\begin{equation}
s(\bm \mu) = \ell^O(u(\bm \mu)) \ .
\label{eq:1p}
\end{equation}
We assume that the FE solution and output are indistinguishable from the exact solution and output at the accuracy level of interest.

For simplicity of exposition, we assume that both $\ell$ and $\ell^O$ are independent of $\bm \mu$. Furthermore, we assume that the bilinear form $a$ may be expressed as an affine decomposition of the form
\begin{equation}
\label{affine}
a(w,v; \bm \mu) = \sum_{q=1}^{Q} \Theta^q(\bm \mu) a^q(w,v) 
\end{equation} 
where $a^q(\cdot, \cdot)$ are $\bm \mu$-independent bilinear forms, and $\Theta^q(\bm \mu)$ are $\bm \mu$-dependent functions. These assumptions can be relaxed by using the empirical interpolation method \cite{Barrault2004, Grepl2007} or its high-order variants \cite{Nguyen2023d,Nguyen2024}.

\subsection{Reduced Basis Approximation}

The parameter sample $S_N = \{\bm \mu_n\}_{n=1}^N$ is generated using a greedy sampling algorithm, which will be discussed in detail later. Once the parameter sample set $S_N$ is defined, we proceed to compute the FE solutions $\{u(\bm \mu_n)\}_{n=1}^N$ by solving the weak formulation (\ref{eq:1a}) for each parameter point in $S_N$. Rather than relying solely on these original solutions to construct the traditional RB space $W_N$, we apply the generative snapshot method, as outlined in the previous section, to construct the generative RB spaces 
\begin{equation}
\label{grbspaces}
\Phi_N^{M_1} = \mbox{span} \{\phi_m, 1 \le m \le M_1 \}, \quad \Psi_N^{M_2} = \mbox{span} \{\psi_m, 1 \le m \le M_2 \} .
\end{equation}
These generative RB spaces incorporate additional snapshots generated through nonlinear transformations to significantly expand the solution space, thereby improving the accuracy of the RB approximation.  Furthermore, we compute {\em a posteriori} error estimates by leveraging these enriched RB spaces, as discussed in the following section.

To develop the RB approximation of the FE formulation,  we perform a Galerkin projection of equation (\ref{eq:1a}) onto $\Phi_N^{M_1}$. Given a parameter point
$\bm \mu \in {\cal D}$, we evaluate $s_{M_1}(\bm 
 \mu) = \ell^O(u_{M_1}(\bm  \mu))$, where $u_{M_1}(\bm \mu) \in \Phi_N^{M_1}$ is the solution of
\begin{equation}
a(u_{M_1}(\bm  \mu),v; \mu) =  \ell(v), \quad
\forall v \in \Phi_N^{M_1}.   
\label{eq:6-24a}
\end{equation} 
We express $u_{M_1}(\bm  \mu) = \sum_{j=1}^{M_1} \alpha_{M_1,j}(\bm  \mu) \phi_j$ to evaluate the RB output  as 
\begin{equation}
\label{sNM1}
s_{M_1}(\bm \mu) = (\bm l_{M_1}^O)^T \bm \alpha_{M_1}(\bm \mu)    
\end{equation}
where $\bm \alpha_{M_1}(\bm \mu) \in \mathbb{R}^{M_1}$ as the solution of the following linear  system
\begin{equation}
\left(\sum_{q=1}^Q \Theta^q(\bm \mu) \bm A_{M_1}^q \right) \bm \alpha_{M_1}(\bm \mu)  = \bm l_{M_1} .
\label{eq:6-24b}
\end{equation} 
Here the matrices $\bm A_{M_1}^q \in \mathbb{R}^{M_1 \times M_1}, 1 \le q \le Q,$ and the vectors $\bm l_{M_1}, \bm l_{M_1}^0  \in \mathbb{R}^{M_1}$ are computed as 
\begin{equation}
A_{M_1, ij}^q  = a^q(\phi_i,\phi_j), \qquad  l_{M_1, i} =  \ell(\phi_i), \qquad  l^O_{M_1, i} =  \ell^O(\phi_i) , 
\label{eq:6-24f}
\end{equation} 
for $1 \leq i,j \leq M_1$. The RB system (\ref{eq:6-24b}) is derived from (\ref{eq:6-24a}) by choosing $v = \phi_i, 1 \le i \le M_1,$ and invoking the affine decomposition (\ref{affine}).

Since $\bm A_{M_1}^q, 1 \le q \le Q,$ and $\bm l_{M_1}, \bm l_{M_1}^O$ are independent of $\bm \mu$, they can be pre-computed and stored in the offline stage. In the online stage, for any $\bm \mu \in \mathcal{D}$, we solve the system (\ref{eq:6-24b}) for $\bm \alpha_{M_1}(\bm \mu)$ with $O(QM_1^2 + M_1^3)$ operation counts and evaluate $s_{M_1}(\bm \mu)$ from (\ref{sNM1}) with $O(M_1)$ operation counts.

\subsection{A Posteriori Error Estimation}

The error estimates provide a way to quantify the accuracy of the RB output and solution. The proposed error estimates rely on a more refined RB approximation, which is constructed using the larger generative RB space $\Psi_N^{M_2}$. Specifically, we evaluate the refined RB output  as 
\begin{equation}
\label{sNM2}
s_{M_2}(\bm \mu) = (\bm l_{M_2}^O)^T \bm \beta_{M_2}(\bm \mu)    
\end{equation}
where $\bm \beta_{M_2}(\bm \mu) \in \mathbb{R}^{M_2}$ as the solution of the following linear  system
\begin{equation}
\left(\sum_{q=1}^Q \Theta^q(\bm \mu) \bm A_{M_2}^q \right) \bm \beta_{M_2}(\bm \mu)  = \bm l_{M_2} .
\label{eq:6-24bw}
\end{equation} 
Here the matrices $\bm A_{M_2}^q \in \mathbb{R}^{M_2 \times M_2}, 1 \le q \le Q,$ and the vectors $\bm l_{M_2}, \bm l_{M_2}^0  \in \mathbb{R}^{M_2}$ are computed as 
\begin{equation}
A_{M_2, ij}^q  = a^q(\psi_i,\psi_j), \qquad  l_{M_1, i} =  \ell(\psi_i), \qquad  l^O_{M_2, i} =  \ell^O(\psi_i) , 
\label{eq:6-24fw}
\end{equation} 
for $1 \leq i,j \leq M_2$. 

The error estimate for $\epsilon_{M_1}^{s}(\bm \mu) = |s(\bm \mu) - s_{M_1}(\bm \mu)|$ is  evaluated as
\begin{equation}
\label{es1}
\bar{\epsilon}_{M_1M_2}^{s}(\bm \mu) =|s_{M_2}(\bm \mu) - s_{M_1}(\bm \mu)| ,
\end{equation}
where $s_{M_2}(\bm \mu)$ is the output computed using $\Psi_N^{M_2}$, and $s_{M_1}(\bm \mu)$  is the output obtained from the smaller RB space $\Phi_N^{M_1}$. Similarly, the estimate for the RB solution error  $\epsilon_{M_1}^{u}(\bm \mu) = \|u(\bm \mu) - u_{M_1}(\bm \mu)\|_X$ is computed as
\begin{equation}
\bar{\epsilon}_{M_1M_2}^{u}(\bm \mu) = \|u_{M_2}(\bm \mu) - u_{M_1}(\bm \mu)\|_X ,
\label{esu}
\end{equation}
where $u_{M_2}(\bm \mu)$ is the solution computed using $\Psi_N^{M_2}$, and $u_{M_1}(\bm \mu)$  is the solution from $\Phi_N^{M_1}$. By substituting  $u_{M_1}(\bm \mu) = \sum_{i=1}^{M_1} \alpha_{M_1,i} \phi_i$ and $u_{M_2}(\bm \mu) = \sum_{m=1}^{M_2} \beta_{M_2,m} \psi_m$ into (\ref{esu}), we obtain 
\begin{multline}
\label{uest}
\bar{\epsilon}_{M_1M_2}^{u}(\bm \mu) =  \left(\|u_{M_2}(\bm \mu) \|_X^2 + \|u_{M_1}(\bm \mu) \|_X^2 - 2 (u_{M_2}(\bm \mu), u_{M_1}(\bm \mu) )_X \right)^{1/2} \\
 =  \left( \left( \bm \alpha_{M_1}(\bm \mu) \right)^T \bm B_{M_1} \bm \alpha_{M_1}(\bm \mu) + \left( \bm \beta_{M_1}(\bm \mu) \right)^T \bm B_{M_2} \bm \beta_{M_1}(\bm \mu) \right. \\ \left. - 2 \left( \bm \alpha_{M_1}(\bm \mu) \right)^T \bm B_{M_1 M_2} \bm \beta_{M_2}(\bm \mu) \right)^{1/2}  
\end{multline}
where
\begin{equation}
\label{bq12}
B_{M_1, ij}  = (\phi_i, \phi_j)_X , \quad  B_{M_2, ml}  = (\psi_m, \psi_l)_X, \quad B_{M_1M_2, im}  = (\phi_i, \psi_m)_X      
\end{equation}
for $1 \le i,j \le M_1$ and $1 \le m,l \le M_2$. Computing $\bar{\epsilon}_{M_1M_2}^{u}(\bm \mu)$ from (\ref{uest}) requires only $O(M_2^2)$ operation counts.

\subsection{Offline-Online Decomposition}

Offline-online decomposition is a key strategy used in RB methods to efficiently solve parametrized PDEs. This decomposition separates the computationally demanding tasks, performed during the offline stage, from the real-time, low-cost computations required during the online stage. 


The offline stage of the generative RB approach is summarized in Algorithm 1. The offline stage is computationally intensive as it involves several key steps required to construct ROMs.  In the offline stage, we must compute $N$ FE solutions at selected parameter points to form the initial snapshot set. The generative snapshot method is then applied to generate additional snapshots and compress them to construct generative RB spaces. Moreover, the offline stage requires the computation of  inner products to form several parameter-independent quantities. The computational complexity of the offline stage is typically dominated by the number of FE solves. 

\begin{algorithm}
\begin{algorithmic}[1]
\REQUIRE{The parameter sample set $S_N = \{\bm \mu_n, 1 \le n \le N\}$.}
\ENSURE{$\bm l_{M_1}, \bm l_{M_1}^O, \bm A_{M_1}^q, \bm l_{M_2}, \bm l_{M_2}^O, \bm A_{M_2}^q, \bm B_{M_1}$, $\bm B_{M_2}$, $\bm B_{M_1M_2}$.}
\STATE{Solve the FE formulation (\ref{eq:1a}) for each $\bm \mu_n \in S_N$ to obtain $\{u(\bm \mu_n)\}_{n=1}^N$.}
\STATE{Construct the generative RB spaces defined in (\ref{grbspaces}) using the generative snapshot method.}
\STATE{Compute and store $\bm l_{M_1}, \bm l_{M_1}^O, \bm A_{M_1}^q$ from (\ref{eq:6-24f}), and $\bm l_{M_2}, \bm l_{M_2}^O, \bm A_{M_2}^q$ from (\ref{eq:6-24fw}).}
\STATE{Compute and store $\bm B_{M_1}$, $\bm B_{M_2}$, $\bm B_{M_1M_2}$ from (\ref{bq12}).}
\end{algorithmic}
\caption{Offline stage of the generative RB approach.}
\end{algorithm}

The online stage of the generative RB approach is summarized in Algorithm 2. The online stage computes the RB output $s_{M_1}(\bm \mu)$, and the error estimates $\bar{\epsilon}_{M_1M_2}^{s}(\bm \mu)$ and $\bar{\epsilon}_{M_1M_2}^{u}(\bm \mu)$ for any given input $\bm \mu \in \mathcal{D}$. The computational complexity to evaluate the RB output is $O(QM_1^2 + M_1^3)$. On the other hand, the computational complexity to evaluate error estimates is $O(QM_2^2 + M_2^3)$. Since $M_2 > M_1$, the overall computational complexity is $O(QM_2^2 + M_2^3)$. The computational cost of the online stage of the generative RB approach is higher than that of the standard RB approach, which is $O(QN^2 + N^3)$. Despite the higher complexity, the larger RB spaces and the error estimates ensure improved accuracy and reliability of the RB approximation, which justify the increased computational cost.

\begin{algorithm}
\begin{algorithmic}[1]
\REQUIRE{Parameter point $\bm \mu \in \mathcal{D}$.}
\ENSURE{RB output $s_{M_1}(\bm \mu)$, error estimates $\bar{\epsilon}_{M_1M_2}^{s}(\bm \mu)$ and $\bar{\epsilon}_{M_1M_2}^{u}(\bm \mu)$.}
\STATE{Form the RB system (\ref{eq:6-24b}) with $O(QM_1^2)$ operation counts.}
\STATE{Solve the RB system (\ref{eq:6-24b}) for $\bm \alpha_{M_1}(\bm \mu)$ and calculate $s_{M_1}(\bm \mu)$ from (\ref{sNM1}), which costs $O(M_1^3)$ operation counts.}
\STATE{Form the RB system (\ref{eq:6-24bw}) with $O(QM_2^2)$ operation counts.}
\STATE{Solve the RB system (\ref{eq:6-24bw}) for $\bm \beta_{M_2}(\bm \mu)$ and calculate $s_{M_2}(\bm \mu)$ from (\ref{sNM2}), which costs $O(M_2^3)$ operation counts.}
\STATE{Compute $\bar{\epsilon}_{M_1M_2}^{s}(\bm \mu)$ from (\ref{es1}) and $\bar{\epsilon}_{M_1M_2}^{u}(\bm \mu)$ from (\ref{uest}), which costs $O(M^2_2)$ operation counts.}
\end{algorithmic}
\caption{Online stage of the generative RB approach.}
\end{algorithm}

\subsection{Greedy Sampling Algorithm}

The greedy sampling algorithm is an iterative procedure employed in RB methods to select a set of parameter points from the parameter space that maximizes the accuracy of the ROM. The greedy algorithm constructs a reduced basis with the minimal number of basis functions required to achieve a desired level of accuracy. The greedy algorithm leverages the offline-online decomposition to efficiently construct the ROM during the offline phase and evaluate it  during the online phase. The  algorithm  relies on inexpensive and tight {\em a posteriori} error estimates to explore the parameter space effectively without requiring a large number of FOM solves. The greedy algorithm of the generative RB approach is summarized in Algorithm 3.

\begin{algorithm}
\begin{algorithmic}[1]
\REQUIRE{The initial parameter sample $S_N = \{\bm \mu_n\}_{n=1}^N$, the training set $\Xi_K^{\rm train} = \{{\bm \mu}_k^{\rm train}\}_{k=1}^K$, and error tolerance $\epsilon_{\rm tol}$.}
\ENSURE{$\bm l_{M_1}, \bm l_{M_1}^O, \bm A_{M_1}^q, \bm l_{M_2}, \bm l_{M_2}^O, \bm A_{M_2}^q, \bm B_{M_1}$, $\bm B_{M_2}$, $\bm B_{M_1M_2}$.}
\STATE{Perform Algorithm 1 to form and store $\bm l_{M_1}, \bm l_{M_1}^O, \bm A_{M_1}^q, \bm l_{M_2}, \bm l_{M_2}^O, \bm A_{M_2}^q, \bm B_{M_1}$, $\bm B_{M_2}$, $\bm B_{M_1M_2}$.}
\STATE{Perform Algorithm 2 to calculate $\bar{\epsilon}_{M_1M_2}^{s}({\bm \mu}_k^{\rm train}), 1 \le k \le K$.}
\STATE{Find the next parameter point as $\bm \mu_{N+1} = \arg \max_{1 \le k \le K} \epsilon_{M_1M_2}^{s}({\bm \mu}_k^{\rm train})$.}
\STATE{If $\bar{\epsilon}_{M_1M_2}^{s}(\bm \mu_{N+1}) \le {\epsilon}_{\rm tol}$ then stop. Otherwise update the parameter sample $S_N = S_N \cup {\bm \mu}_{N+1}$ and go back to Step 1.}
\end{algorithmic}
\caption{Greedy sampling of the generative RB approach.}
\end{algorithm}

\section{Numerical results}

In this section, we present numerical results from two affine linear PDEs to demonstrate the generative RB method and compare its performance to that of the standard RB method. To this end, we introduce the following maximum relative errors 
\begin{equation}
{\epsilon}_{\rm{max},{\rm rel}}^u =  \max_{\bm \mu \in \Xi^{\rm test}} \frac{\epsilon_{M_1}^u(\bm \mu)}{\|u(\bm \mu)\|_X} , \quad {\epsilon}_{\rm{max},{\rm rel}}^s = \max_{\bm \mu \in \Xi^{\rm test}} \frac{\epsilon_{M_1}^s (\bm \mu)}{|s(\bm \mu)|}    
\end{equation}  
and their estimates as 
\begin{equation}
\bar{{\epsilon}}_{\rm{max},{\rm rel}}^u =  \max_{\bm \mu \in \Xi^{\rm test}} \frac{\bar{\epsilon}_{M_1,M_2}^u(\bm \mu)}{\|u_{M_2}(\bm \mu)\|_X} , \quad \bar{\epsilon}_{\rm{max},{\rm rel}}^s = \max_{\bm \mu \in \Xi^{\rm test}} \frac{\bar{\epsilon}_{M_1,M_2}^s (\bm \mu)}{|s_{M_2}(\bm \mu)|}    ,
\end{equation}  
where $\Xi^{\rm test}$ is a test sample of $N^{\rm test}$ parameter points  distributed uniformly in the parameter domain. Similarly, the mean relative errors, ${\epsilon}_{\rm{mean},{\rm rel}}^u$ and ${\epsilon}_{\rm{mean},{\rm rel}}^s$, and their estimates, $\bar{\epsilon}_{\rm{mean},{\rm rel}}^u$ and $\bar{\epsilon}_{\rm{mean},{\rm rel}}^s$, are  defined.


\subsection{Convection-Diffusion Problem}

We consider a  linear convection-diffusion problem  
\begin{equation}
- \nabla^2u^{\rm e} + \nabla \cdot (\bm \mu u^{\rm e}) = 100, \quad \mbox{in } \Omega , 
\end{equation} 
with homogeneous Dirichlet condition $u^{\rm e}=0$ on  $\partial \Omega$. Here $\Omega = (0,1)^2$ is a unit square domain, while $\bm \mu = (\mu_1, \mu_2)$ is the parameter vector in a parameter domain $\mathcal{D} \equiv [0, 50] \times [0, 50]$. Let $X \in H_0^1(\Omega)$ be a finite element (FE) approximation space of dimension $\mathcal{N}$ with $X = \{v \in H_0^1(\Omega) : v|_K \in \mathcal{P}^3(T), \  \forall T \in \mathcal{T}_h \}$, where $\mathcal{P}^3(T)$ is a space of polynomials of degree $3$ on an element $T \in \mathcal{T}_h$ and $\mathcal{T}_h$ is a finite element grid of $32 \times 32$ quadrilaterals. The dimension of the FE space $X$ is $\mathcal{N} = 9409$. The FE approximation $u(\bm \mu) \in X$ of the exact solution $u^{\rm e}(\bm \mu)$ is the solution of
\begin{equation}
\int_\Omega \nabla u \cdot \nabla v -  \int_\Omega \bm \mu u \cdot \nabla v  = 100 \int_{\Omega}  v, \quad \forall  v
\in X \ .
\label{eq:7-6c}
\end{equation}
The output of interest is evaluated as $s(\bm \mu) = \ell^{O}(u(\bm \mu))$ with $\ell^{O}(v) \equiv \int_{\Omega} v$.  

\begin{figure}[htbp]
	\centering
	\begin{subfigure}[b]{0.49\textwidth}
		\centering		\includegraphics[width=\textwidth]{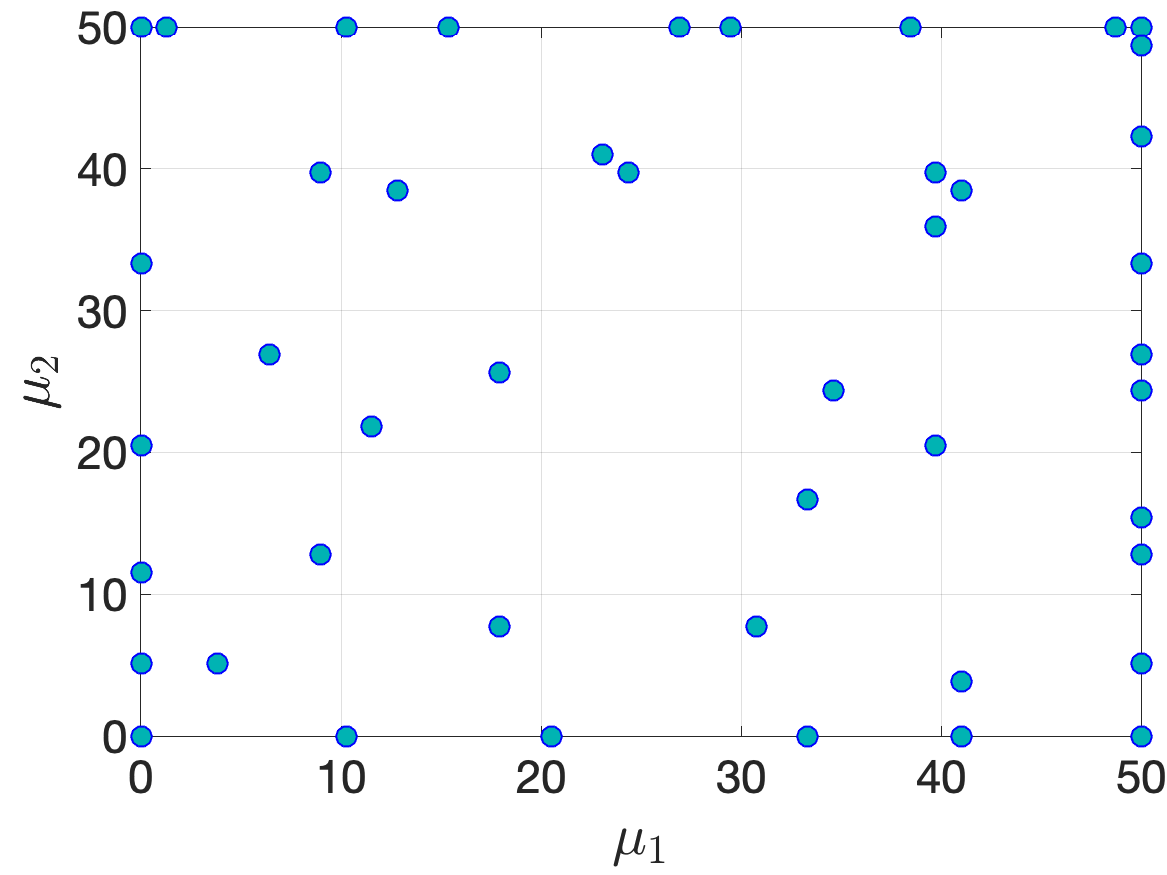}
		\caption{Parameter points in the sample $S_N$ for $N=45$.}
	\end{subfigure}
	\hfill
	\begin{subfigure}[b]{0.49\textwidth}
		\centering		\includegraphics[width=\textwidth]{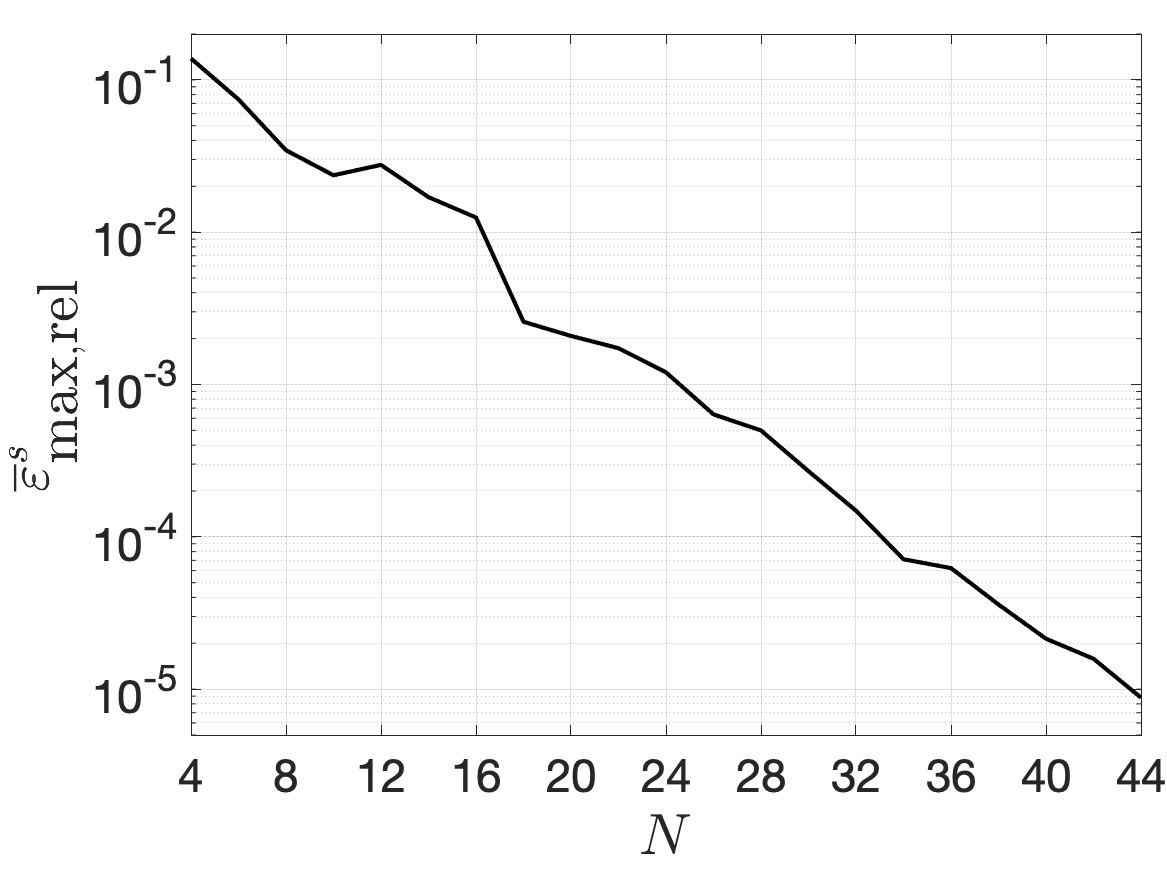}
		\caption{Maximum relative error estimate.}
	\end{subfigure}
	\caption{Greedy sampling algorithm results: (a) distribution of the parameter sample points selected from a grid of $40 \times 40$ uniformly distributed points in the parameter domain, and (b) the convergence of the maximum relative error estimate for the RB output.}
	\label{ex4fig1}
\end{figure}

Figure \ref{ex4fig1} shows the selected parameter points and the maximum relative error estimate as a function of $N$. The plot reveals that the majority of the points are concentrated along the boundary of the parameter domain. This distribution suggests that the regions near the boundary exhibit greater variability or complexity in the solution manifold, requiring more refined sampling. The greedy algorithm terminates at $N = 45$, where the maximum relative error estimate is less than $10^{-5}$.

\begin{figure}[h!]
	\centering
	\begin{subfigure}[b]{0.32\textwidth}
		\centering		\includegraphics[width=\textwidth]{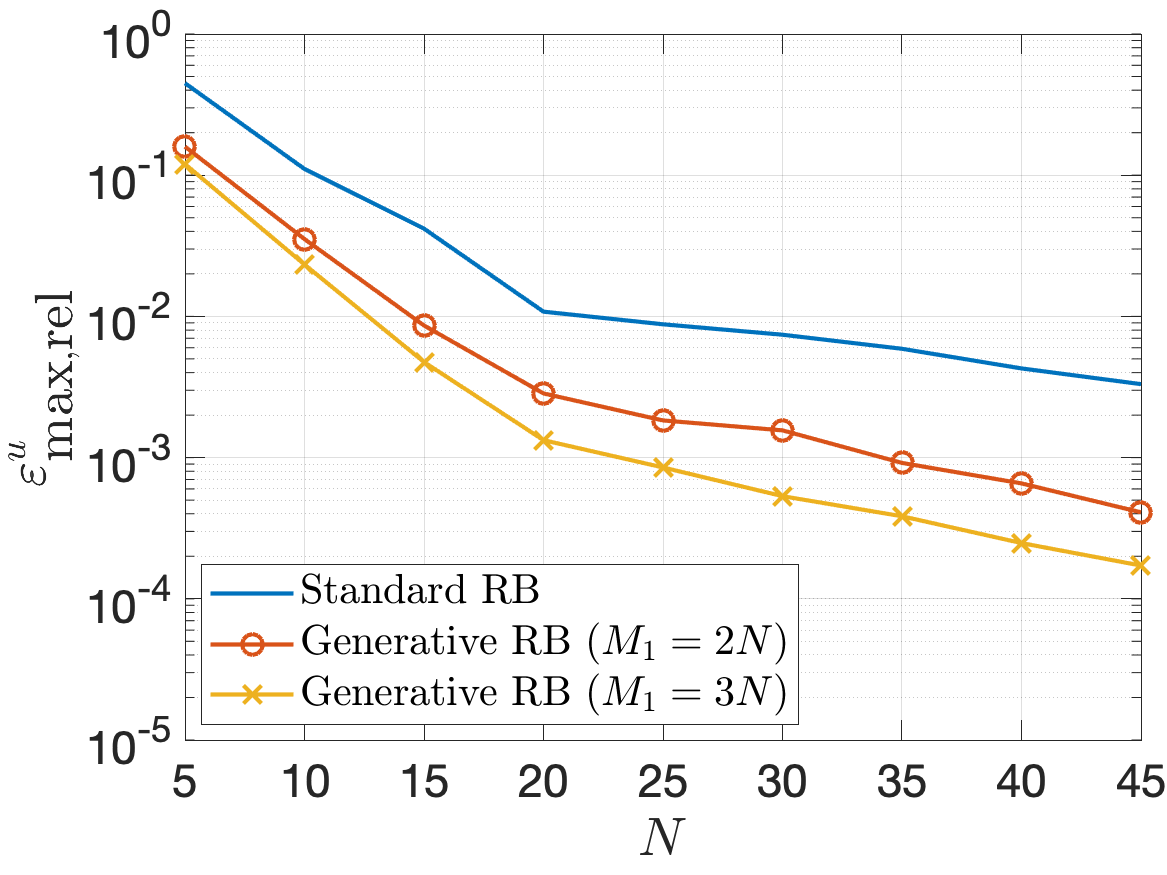}		
	\end{subfigure}	
	\begin{subfigure}[b]{0.32\textwidth}
		\centering		\includegraphics[width=\textwidth]{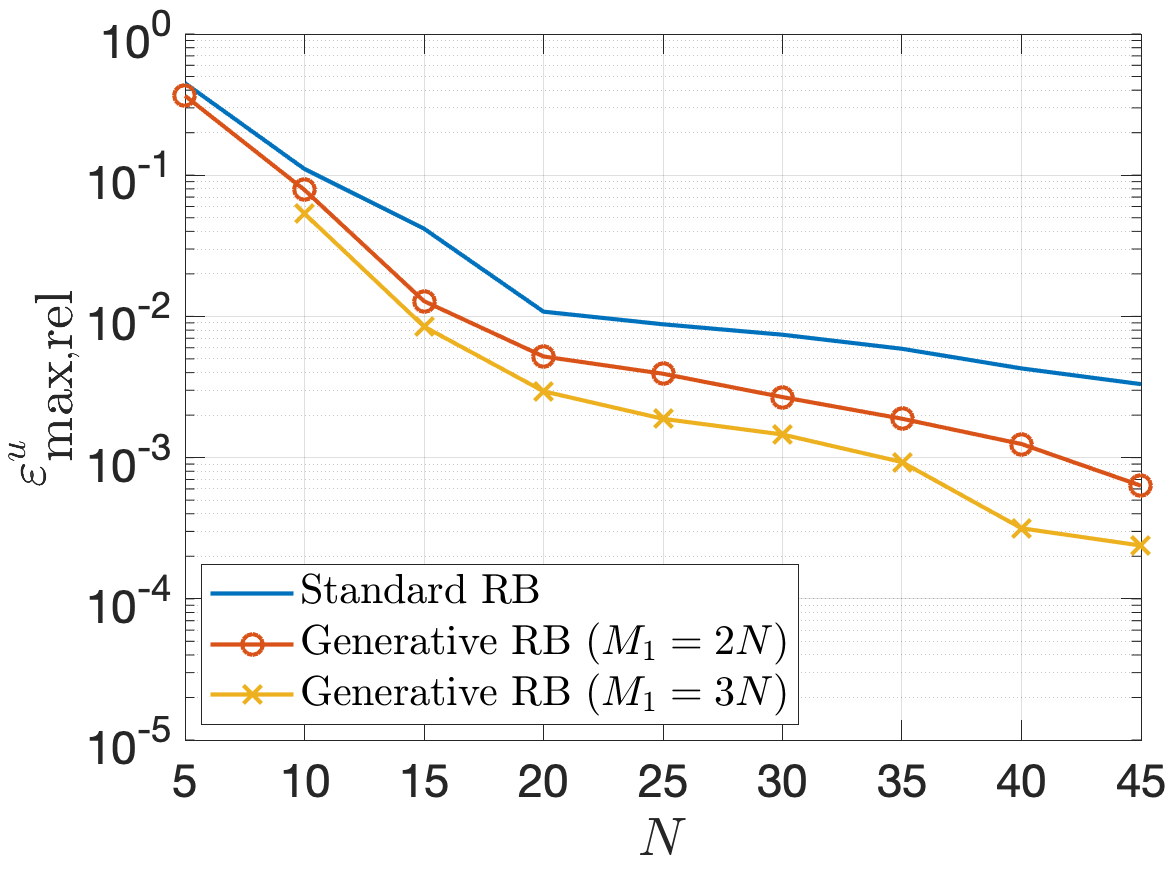}		
	\end{subfigure}        
 	\begin{subfigure}[b]{0.32\textwidth}
		\centering		\includegraphics[width=\textwidth]{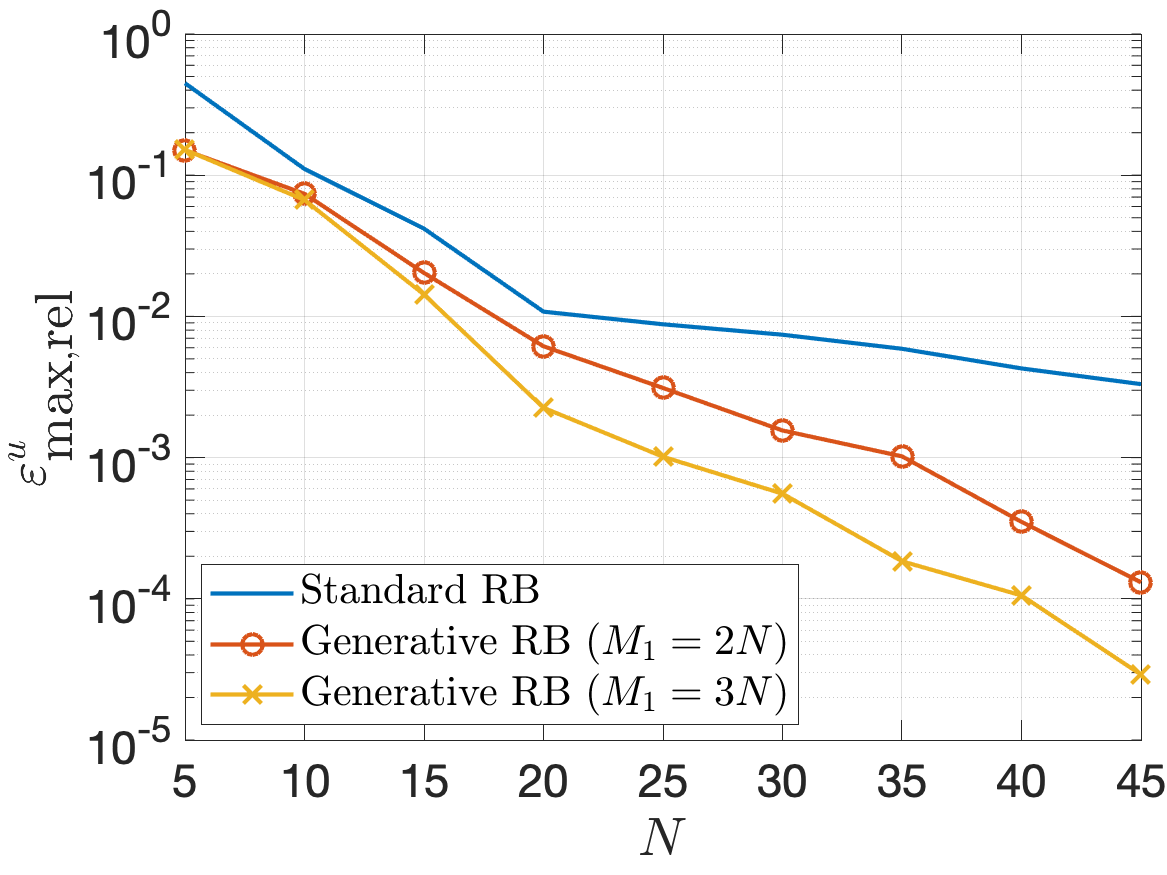}		
	\end{subfigure}
        \\
	\begin{subfigure}[b]{0.32\textwidth}
		\centering		\includegraphics[width=\textwidth]{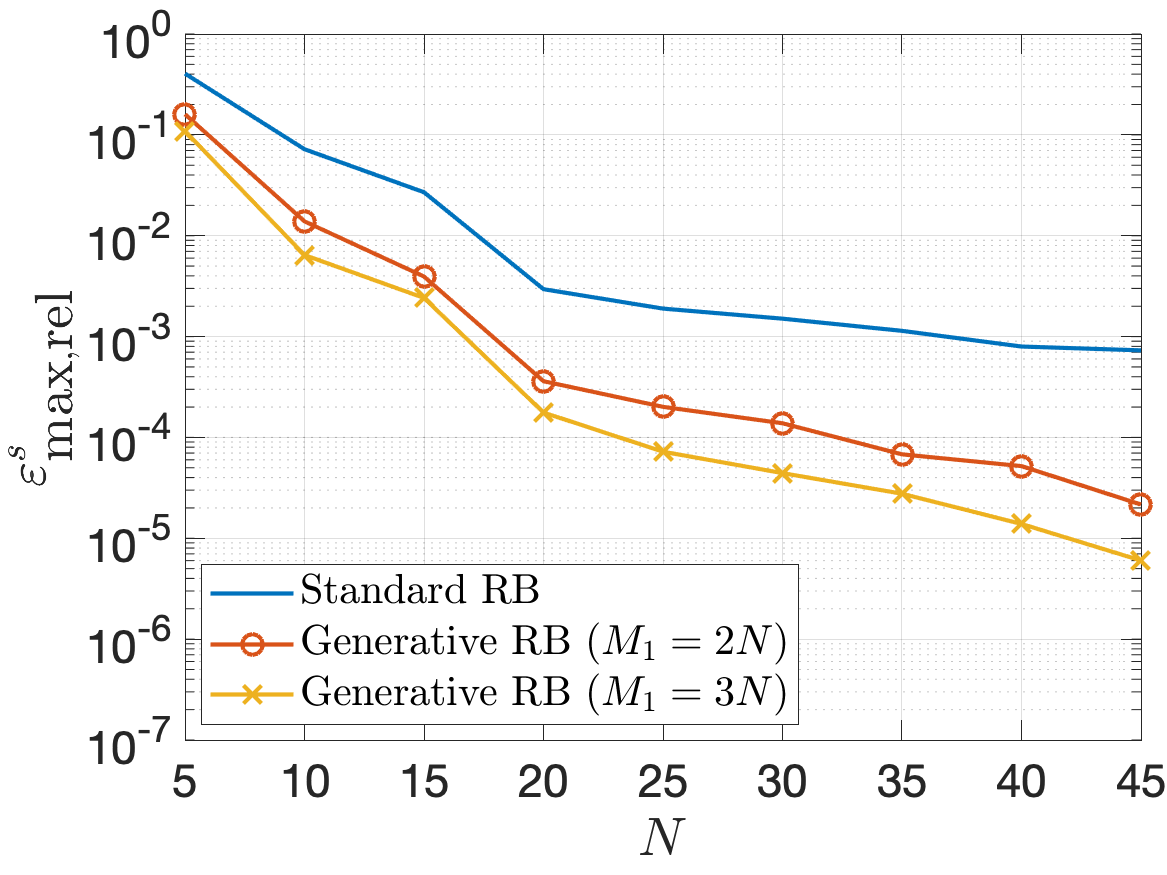}		
        \caption{Softplus function}
	\end{subfigure}         
  	\begin{subfigure}[b]{0.32\textwidth}
		\centering		\includegraphics[width=\textwidth]{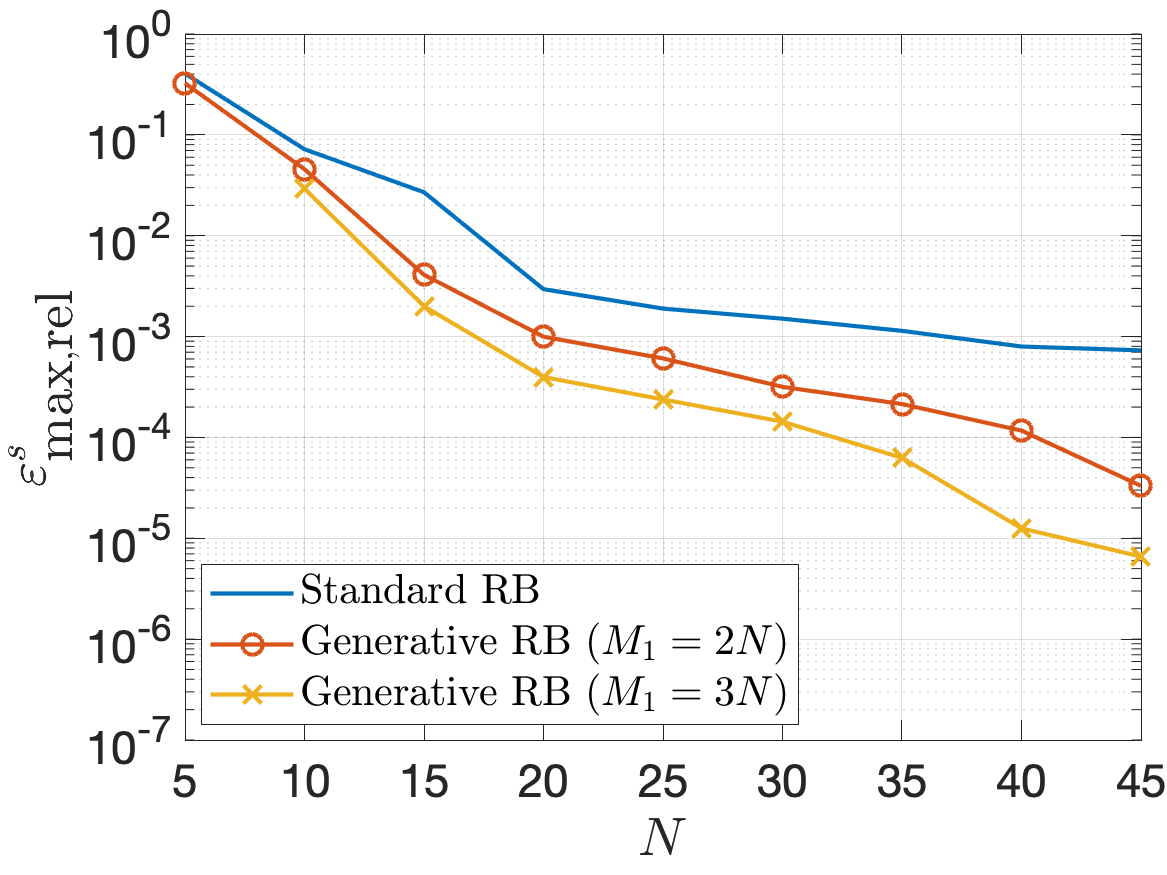}		
        \caption{Exponential function}
	\end{subfigure}
	\begin{subfigure}[b]{0.32\textwidth}
		\centering		\includegraphics[width=\textwidth]{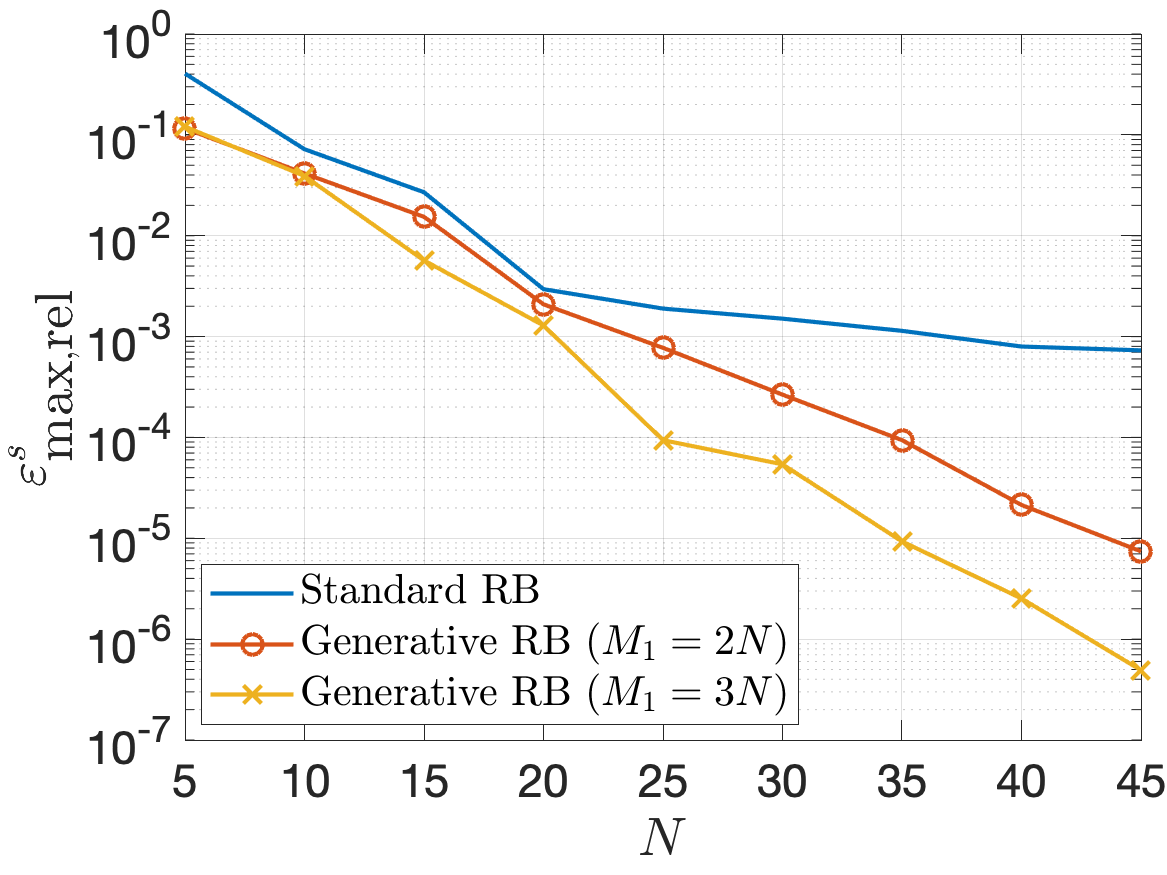}		
    \caption{Quadratic function}
	\end{subfigure}
	\caption{Convergence of the maximum relative error in the RB solution (top) and in the RB output (bottom) as a function of $N$.}
	\label{ex4fig2}
\end{figure}

Figure \ref{ex4fig2} presents the maximum relative errors, ${\epsilon}_{\rm{max},{\rm rel}}^u$ and ${\epsilon}_{\rm{max},{\rm rel}}^s$, as a function of $N$ for both the standard RB method and the generative RB method. The test set $\Xi^{\rm test}$ is a uniform grid of $N^{\rm test} = 30 \times 30$ points in the parameter domain. The generative RB method consistently shows lower errors compared to the standard RB method across all three activation functions considered. Among these, the softplus and exponential functions exhibit similar convergence behaviors, while the quadratic function achieves the fastest convergence rate and the smallest errors, particularly for larger values of $N$. In particular, at $N = 45$, the generative RB method with $M_1 = 3N$ produces output errors that are approximately 1000 times smaller than those of the standard RB method.

\begin{figure}[h!]
	\centering
	\begin{subfigure}[b]{0.32\textwidth}
		\centering		\includegraphics[width=\textwidth]{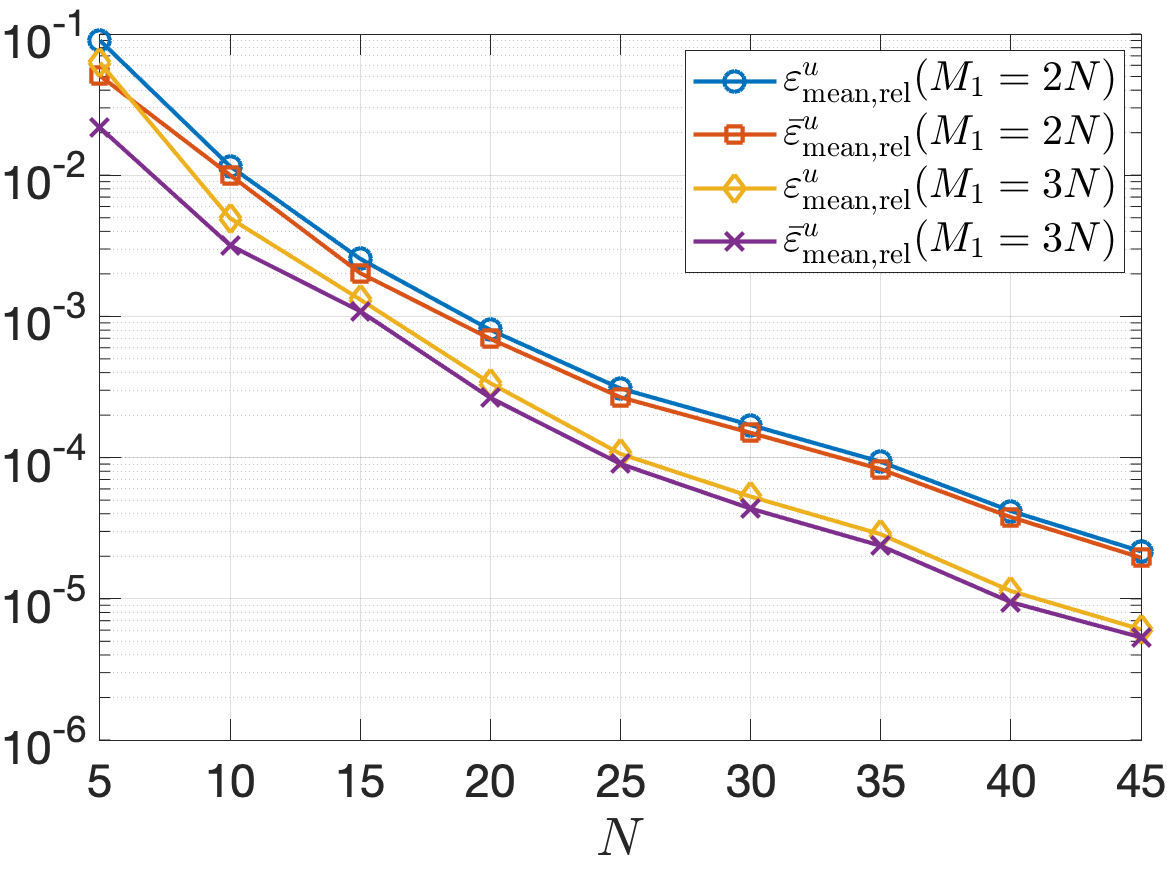}		
	\end{subfigure}	
	\begin{subfigure}[b]{0.32\textwidth}
		\centering		\includegraphics[width=\textwidth]{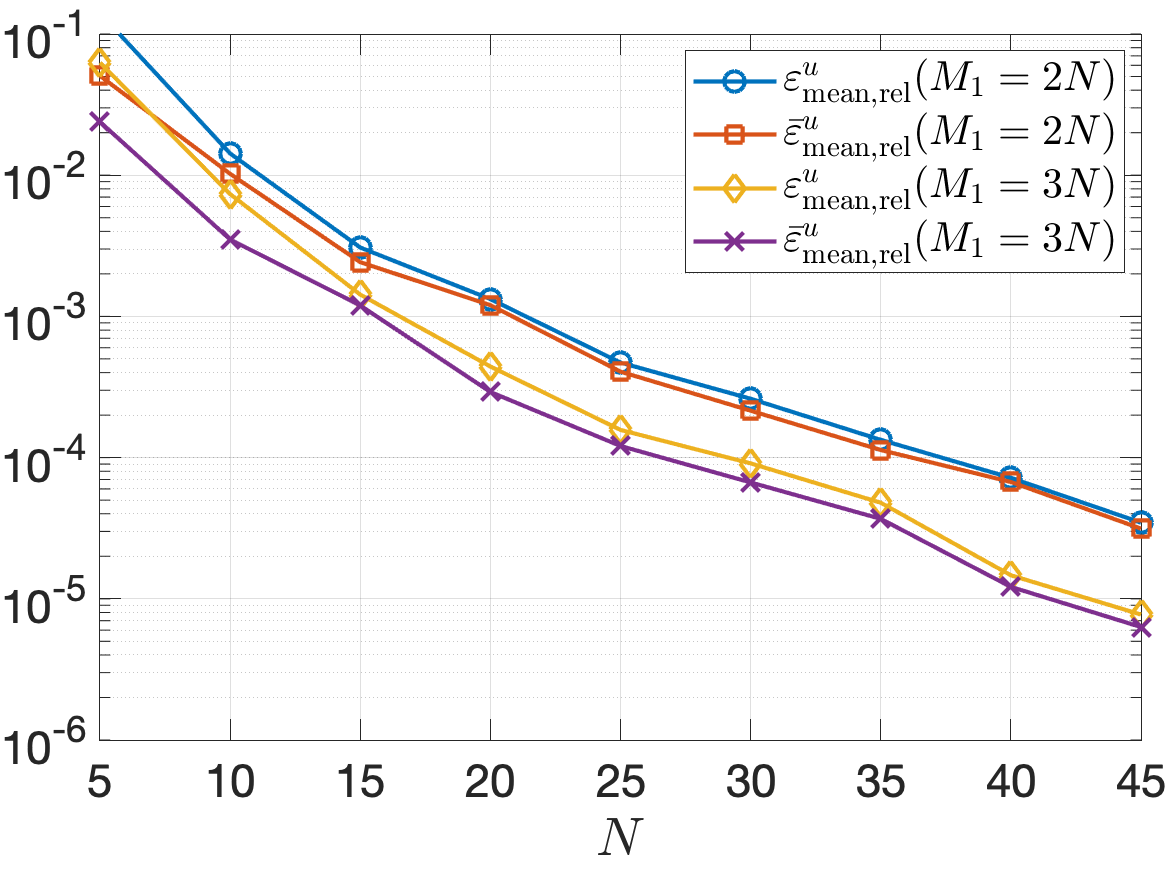}		
	\end{subfigure}        
 	\begin{subfigure}[b]{0.32\textwidth}
		\centering		\includegraphics[width=\textwidth]{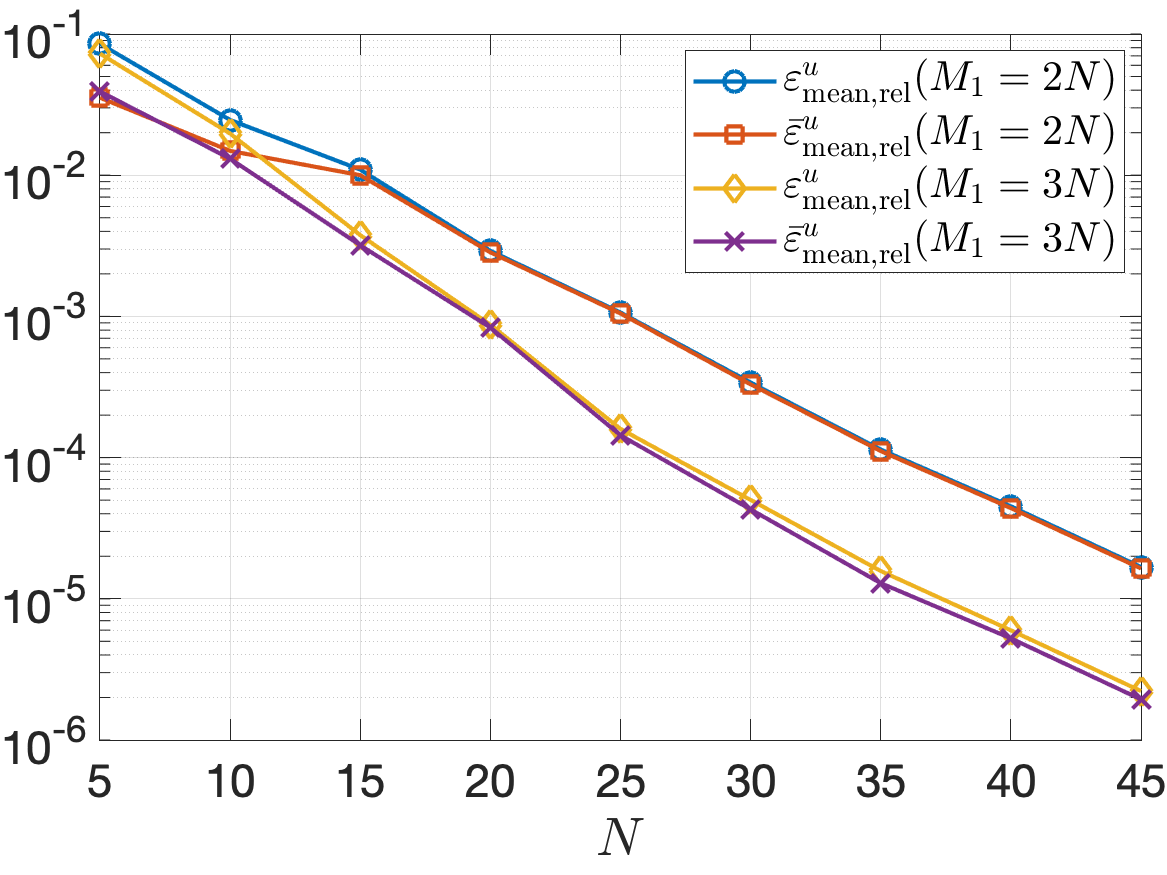}		
	\end{subfigure}
        \\
	\begin{subfigure}[b]{0.32\textwidth}
		\centering		\includegraphics[width=\textwidth]{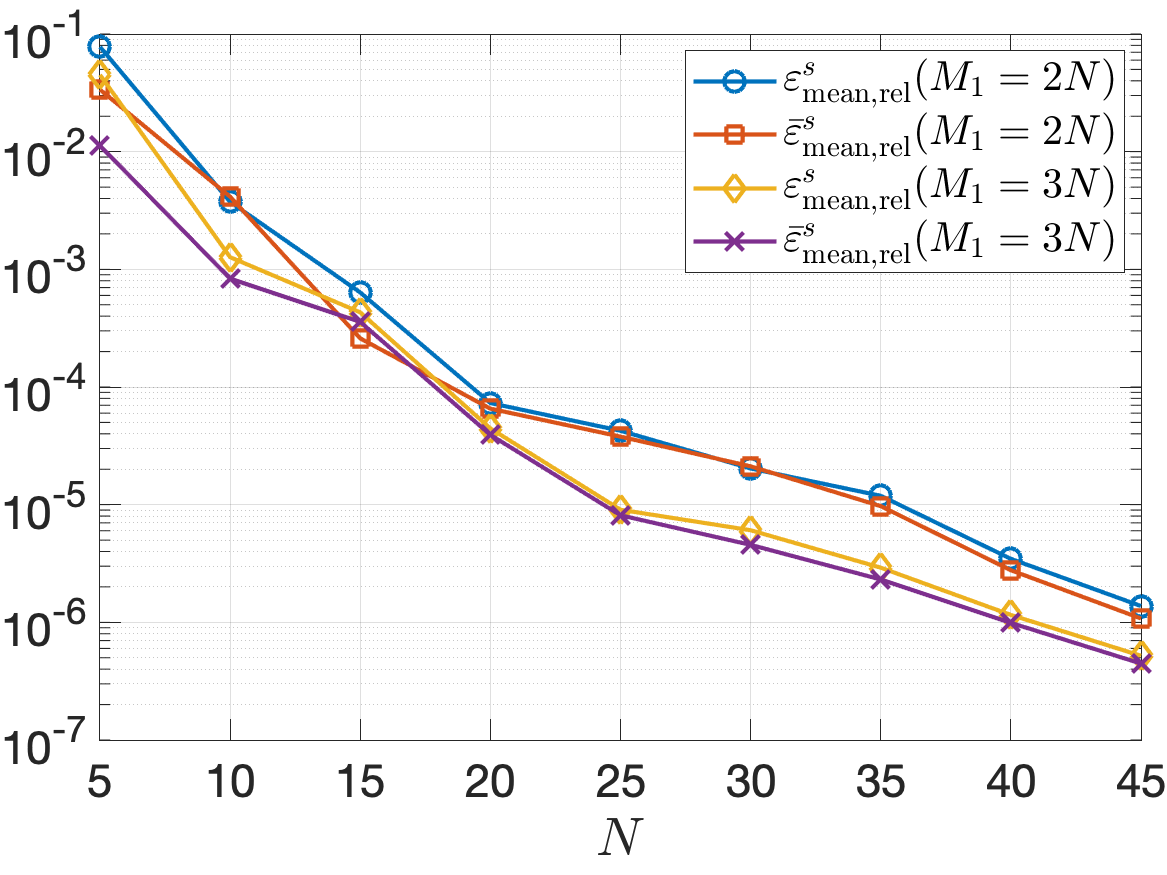}		
        \caption{Softplus function}
	\end{subfigure}         
  	\begin{subfigure}[b]{0.32\textwidth}
		\centering		\includegraphics[width=\textwidth]{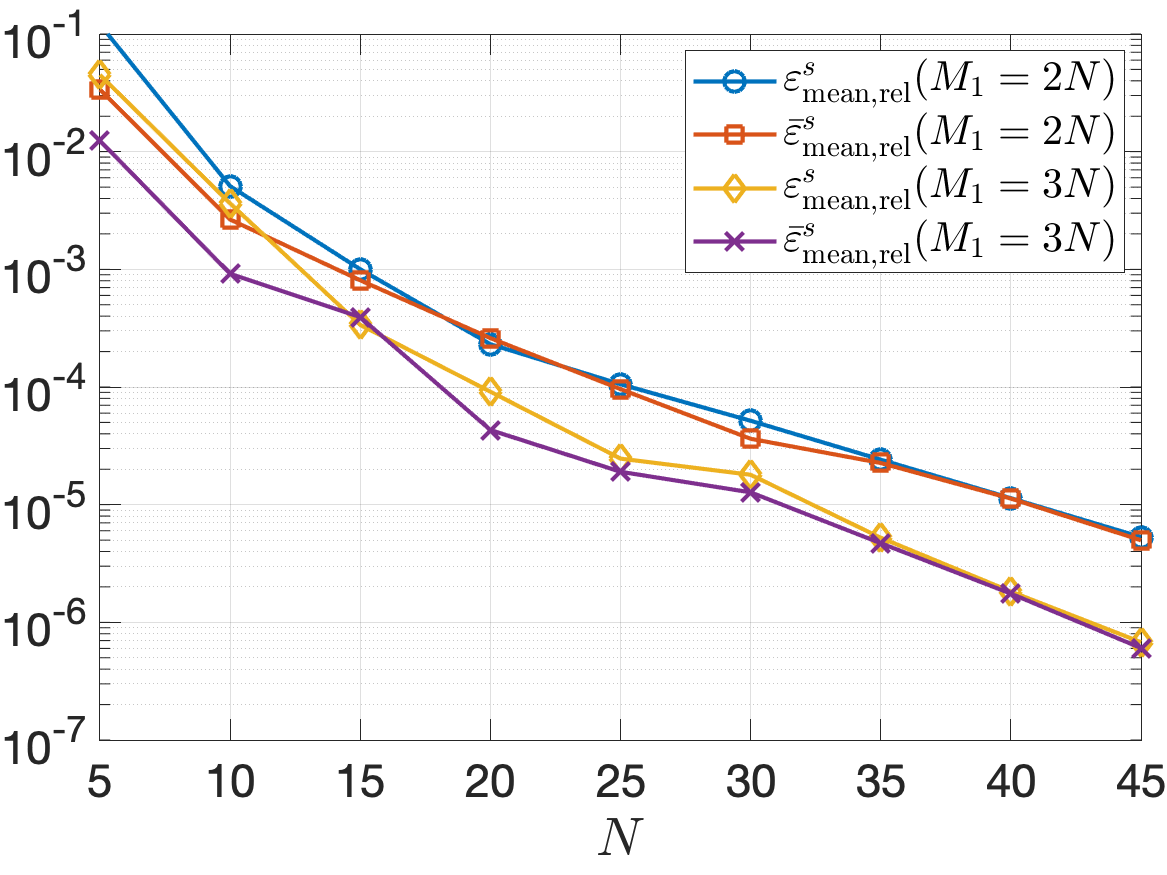}		
        \caption{Exponential function}
	\end{subfigure}
	\begin{subfigure}[b]{0.32\textwidth}
		\centering		\includegraphics[width=\textwidth]{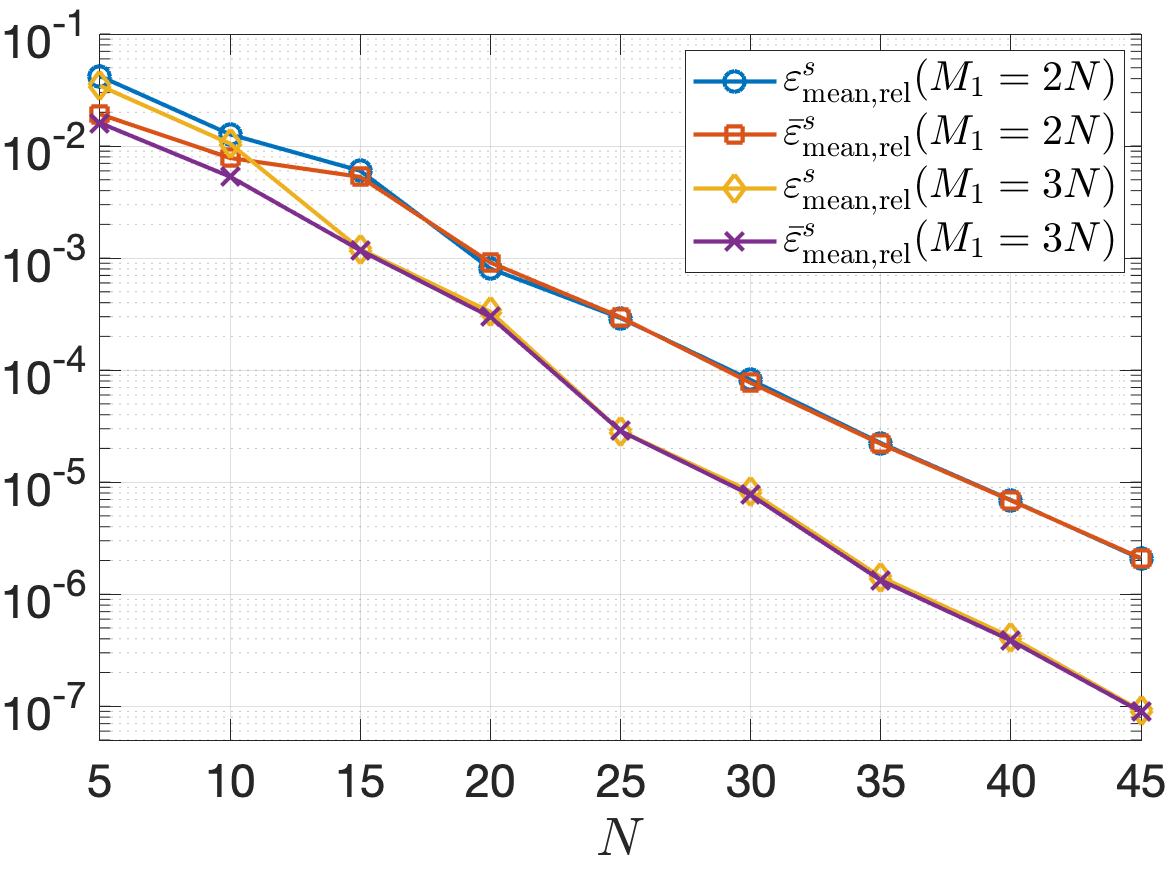}		
    \caption{Quadratic function}
	\end{subfigure}
	\caption{Convergence of the mean relative error and the mean error estimate in the RB solution (top) and in the RB output (bottom) as a function of $N$. The error estimates are calculated using $M_2 = M_1 + N$.}
	\label{ex4fig3}
\end{figure}

Figure \ref{ex4fig3} shows the mean relative errors, ${\epsilon}_{\rm{mean},{\rm rel}}^u$ and ${\epsilon}_{\rm{mean},{\rm rel}}^s$, together with their error estimates, $\bar{\epsilon}_{\rm{mean},{\rm rel}}^u$ and $\bar{\epsilon}_{\rm{mean},{\rm rel}}^s$, as functions of $N$. The error estimates are calculated using $M_2 = M_1 + N$. The error estimates closely match the true errors, indicating that the error estimates are very tight and provide a reliable measure of accuracy. Although error estimates are not provably rigorous -- evident in cases where they are slightly smaller than true errors -- their tightness serves as a strong indicator of the accuracy and reliability of the generative RB method. This tight correlation between estimated and true errors reinforces the efficacy of the generative RB method.


\subsection{Reaction-Diffusion Problem}

\begin{figure}[hthbp]
	\centering
	\begin{subfigure}[b]{0.285\textwidth}
		\centering		\includegraphics[width=\textwidth]{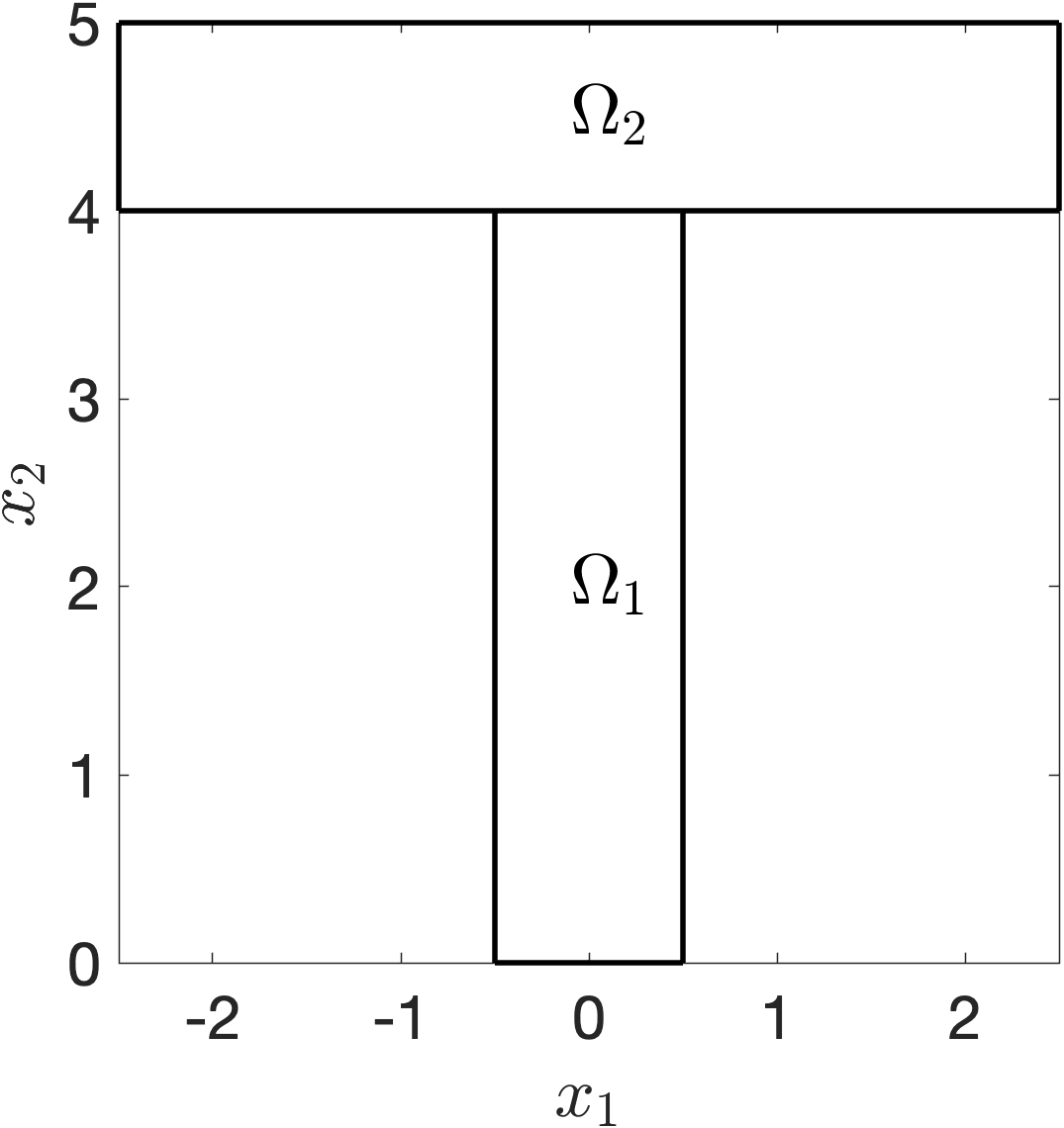}
		\caption{Physical domain}
	\end{subfigure}
	\hfill
	\begin{subfigure}[b]{0.33\textwidth}
		\centering
		\includegraphics[width=\textwidth]{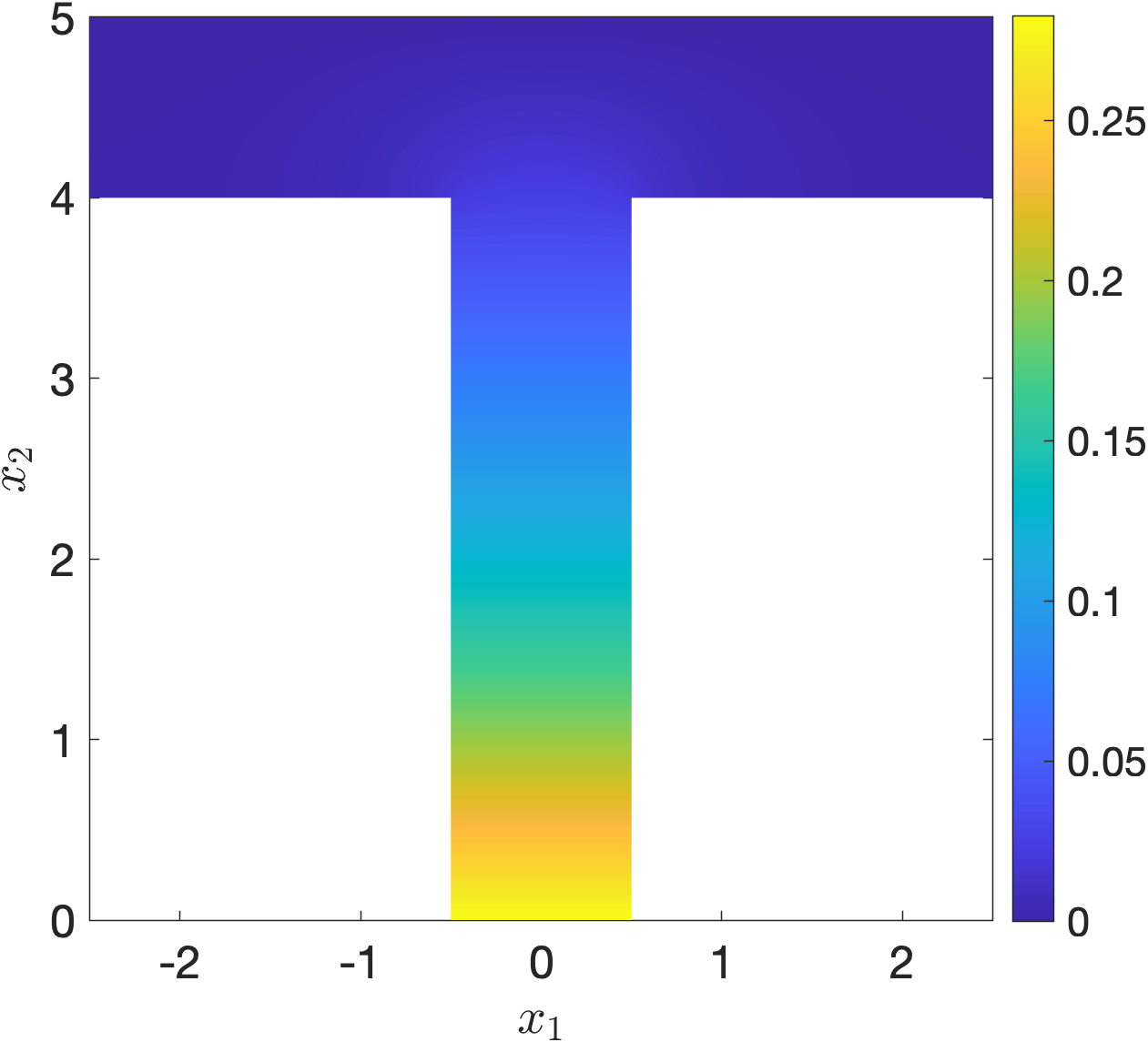}
		\caption{Solution for $\bm \mu = (10,10,0,0)$}
	\end{subfigure}
        \hfill
	\begin{subfigure}[b]{0.33\textwidth}
		\centering
		\includegraphics[width=\textwidth]{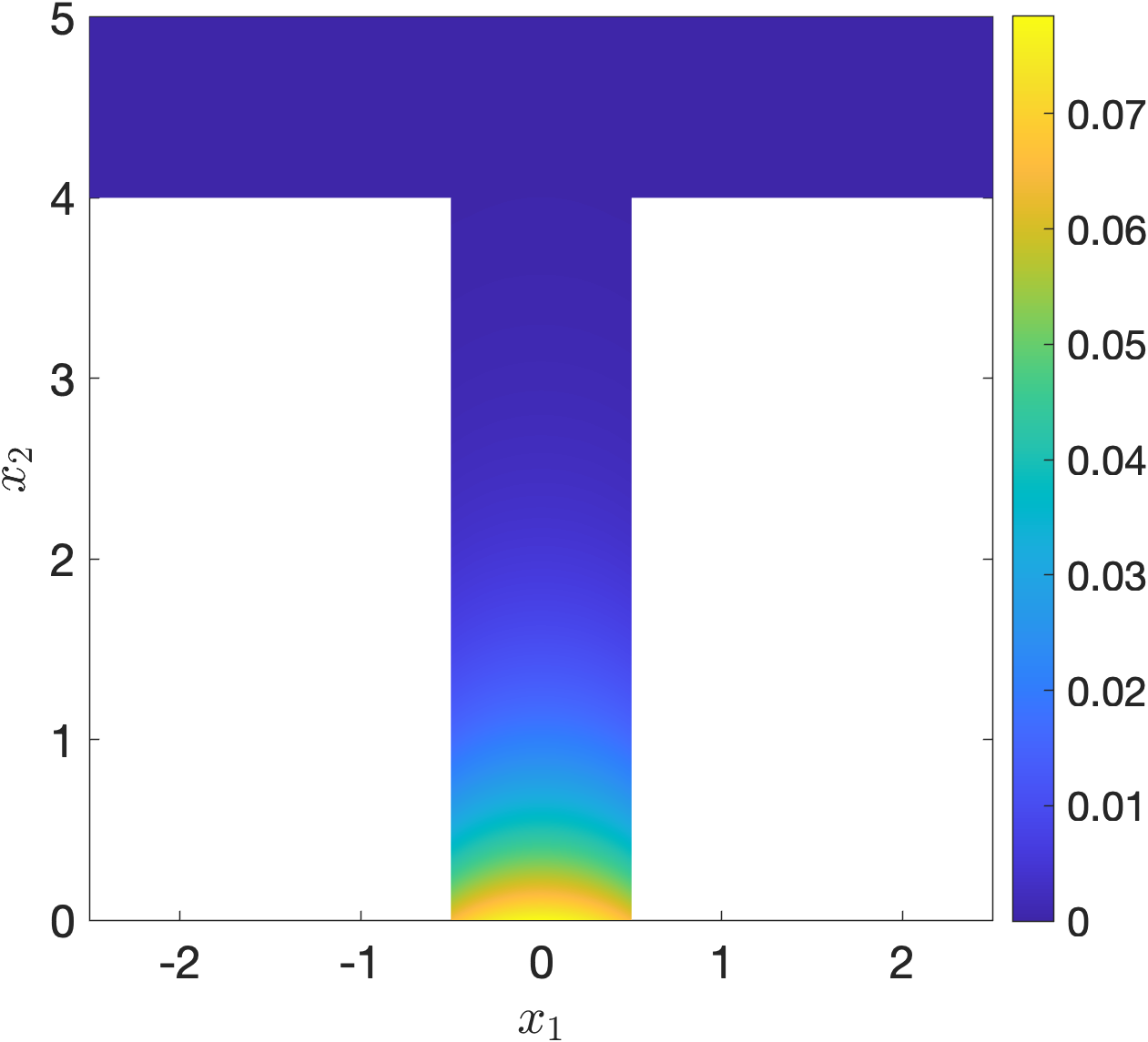}
		\caption{Solution for $\bm \mu = (10,10,10,10)$}
	\end{subfigure}
	\caption{Physical domain  and FE solutions at two parameter points for the linear reaction-diffusion problem.}
	\label{ex5fig0}
\end{figure}

We consider a linear reaction-diffusion problem 
\begin{equation}
-\nabla \cdot \left(\kappa(\bm x, \bm \mu) \nabla u \right) + u = 0 \quad \mbox{in } \Omega, 
\end{equation} 
 with  Dirichlet and Neumann boundary conditions
 \begin{equation}
 u = 0 \quad \mbox{ on } \Gamma_{2}, \quad \kappa(\bm x, \bm \mu)  \nabla u \cdot \bm n = 1 \quad \mbox{on } \Gamma_{\rm 1} ,
\end{equation} 
 and Robin boundary conditions
 \begin{equation}
 \kappa(\bm x, \bm \mu)  \nabla u \cdot \bm n + \mu_3 u = 0 \quad \mbox{on } \Gamma_{\rm 3}, \quad  \kappa(\bm x, \bm \mu)  \nabla u \cdot \bm n + \mu_4 u = 0 \quad \mbox{on } \Gamma_{\rm 4} .
\end{equation} 
Here $\Omega = \Omega_1 \cup \Omega_2$ is the $T$-shaped domain as shown in Figure \ref{ex5fig0}(a). $\Gamma_{\rm 1}$ is the bottom  boundary and $\Gamma_{\rm 2}$ is the top boundary, while $\Gamma_{\rm 3}$ and $\Gamma_{\rm 4}$ are the remaining part of the boundary of $\Omega_1$ and $\Omega_2$, respectively. The thermal conductivity has the form
\begin{equation}
\kappa(\bm x, \bm \mu)  = 
\left\{
\begin{array}{ll}
\mu_1, &  \bm x \in \Omega_1, \\ 
\mu_2, &  \bm x \in \Omega_2 .
\end{array}
\right. 
\end{equation}
The parameter domain is ${\cal D} \equiv [1, 10] \times [1,10] \times [0,10] \times [0, 10]$.  The output of interest is the average of the field variable over the
physical domain. For any given $\bm \mu \in
{\cal D}$, we evaluate $s(\bm \mu) = \int_\Omega u(\bm \mu)$, where  $u(\bm \mu) \in X \subset H_0^1(\Omega) \equiv \{v \in
H^1(\Omega) \mbox{ } | \mbox{ } v|_{\Gamma_1} = 0\}$ is the solution of
\begin{equation*}
\mu_1 \int_{\Omega_1} \nabla u \cdot \nabla v + \mu_2 \int_{\Omega_2} \nabla u \cdot \nabla v + \mu_3 \int_{\Gamma_3} uv + \mu_4 \int_{\Gamma_4} uv +\int_{\Omega} uv  = \int_{\Gamma_{\rm 2}} v, \quad \forall  v
\in X  .  
\label{eq:8-6}
\end{equation*}
The FE approximation space is $X = \{v \in H_0^1(\Omega) : v|_K \in \mathcal{P}^3(T), \  \forall T \in \mathcal{T}_h \}$, where $\mathcal{P}^3(T)$ is a space of polynomials of degree $3$ on an element $T \in \mathcal{T}_h$ and $\mathcal{T}_h$ is a mesh of $900$ quadrilaterals. The dimension of $X$ is $\mathcal{N} = 8401$. Figure \ref{ex5fig0} shows FE solutions at two different parameter points.

 Table \ref{tab:data} shows the parameter points and the maximum relative error estimate as a function of $N$ during the greedy cycle. The first four parameter points form the initial parameter set $S_N$. The greedy sampling is applied to a grid of $8 \times 8 \times 8 \times 8$ uniformly distributed parameter points and terminated at $N = 20$, where the maximum relative error estimate is less than $10^{-5}$. The error estimates are calculated using $M_1 = 2N$ and $M_2 = 4N$. All selected parameter points are on the boundary of the parameter domain, suggesting that the regions near the boundary exhibit greater variability or complexity in the solution manifold. 

\begin{table}[ht!]
\centering
\footnotesize
\begin{tabular}{|c|c|c|c|c|c|}
\hline
$N$ & $\mu_1$ & $\mu_2$ & $\mu_3$ & $\mu_4$ &  $\bar{\epsilon}_{\rm{max},{\rm rel}}^s$ \\ \hline
  1  &  1  &  1  &  0  &  0  &  --  \\  
 2  &  10  &  10  &  10  &  10  &  --  \\  
 3  &  10  &  1  &  0  &  10  &  --  \\  
 4  &  1  &  10  &  10  &  0  &  --  \\  
 5  &  10  &  10  &  0  &  0  &  1.06E-01  \\  
 6  &  10  &  1  &  0  &  0  &   6.99E-02  \\  
 7  &  3.5714  &  1  &  10  &  10  &  5.43E-03  \\  
 8  &  10  &  1  &  5.7143  &  0  &  5.61E-03  \\  
 9  &  1  &  1  &  5.7143  &  10  &  2.03E-03  \\  
 10  &  10  &  3.5714  &  0  &  0  &  1.01E-03  \\  
 11  &  10  &  1  &  1.4286  &  0  &  1.13E-03  \\  
 12  &  10  &  1  &  0  &  2.8571  &  8.77E-04  \\  
 13  &  10  &  4.8571  &  0  &  4.2857  & 3.59E-04  \\  
 14  &  6.1429  &  2.2857  &  0  &  10  &  1.78E-04  \\  
 15  &  3.5714  &  1  &  0  &  10  &  7.66E-05  \\  
 16  &  10  &  1  &  10  &  0  &  3.96E-05  \\  
 17  &  1  &  1  &  1.4286  &  7.1429  &  2.79E-05  \\  
 18  &  10  &  6.1429  &  2.8571  &  0  &  1.61E-05  \\  
 19  &  3.5714  &  10  &  0  &  0  &  1.62E-05  \\  
 20  &  10  &  1  &  0  &  7.1429  &  3.34E-06  \\     
 \hline 
\end{tabular}
\caption{Selected parameter points and the maximum relative errors during the greedy sampling cycle for the heat conduction problem. The greedy algorithm is applied to a grid of $8 \times 8 \times 8 \times 8$ uniformly distributed points in the parameter domain.}
\label{tab:data}
\end{table}

\begin{figure}[htbp]
	\centering
	\begin{subfigure}[b]{0.32\textwidth}
		\centering		\includegraphics[width=\textwidth]{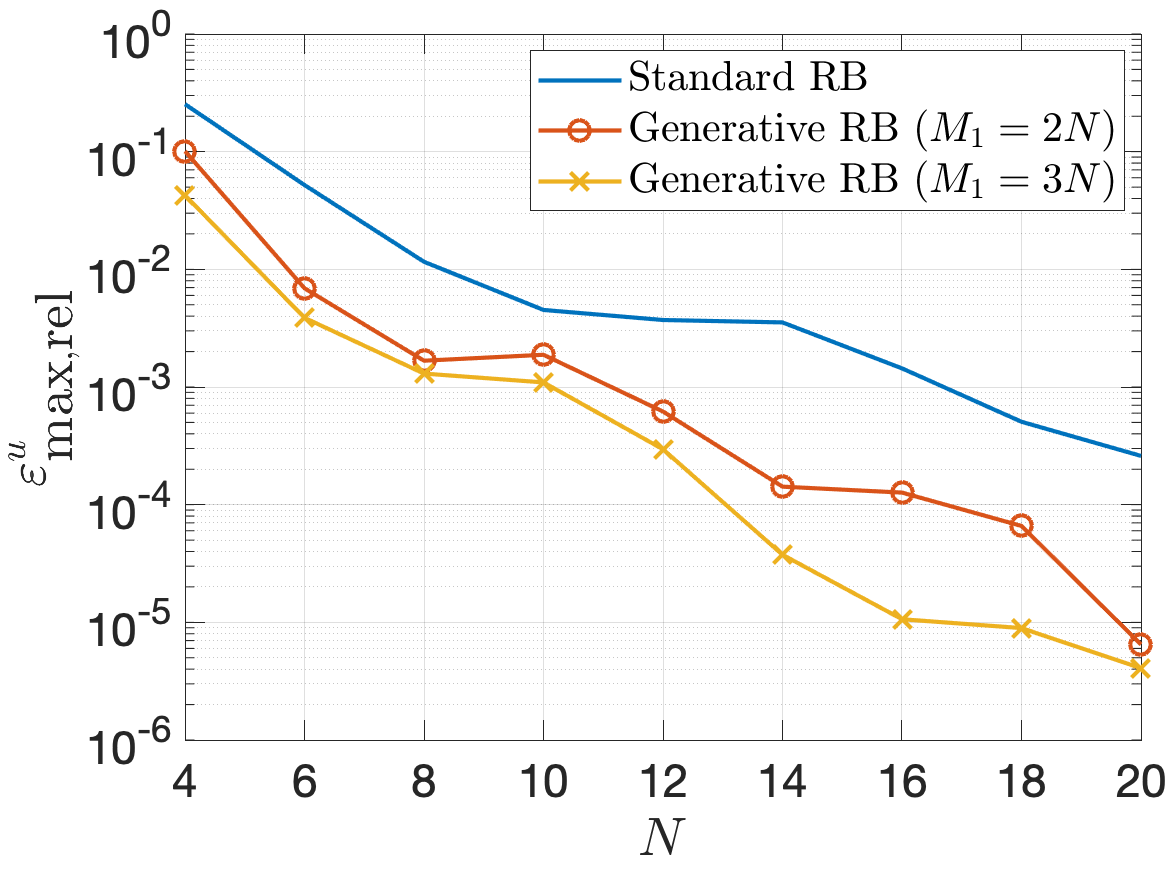}		
	\end{subfigure}	
	\begin{subfigure}[b]{0.32\textwidth}
		\centering		\includegraphics[width=\textwidth]{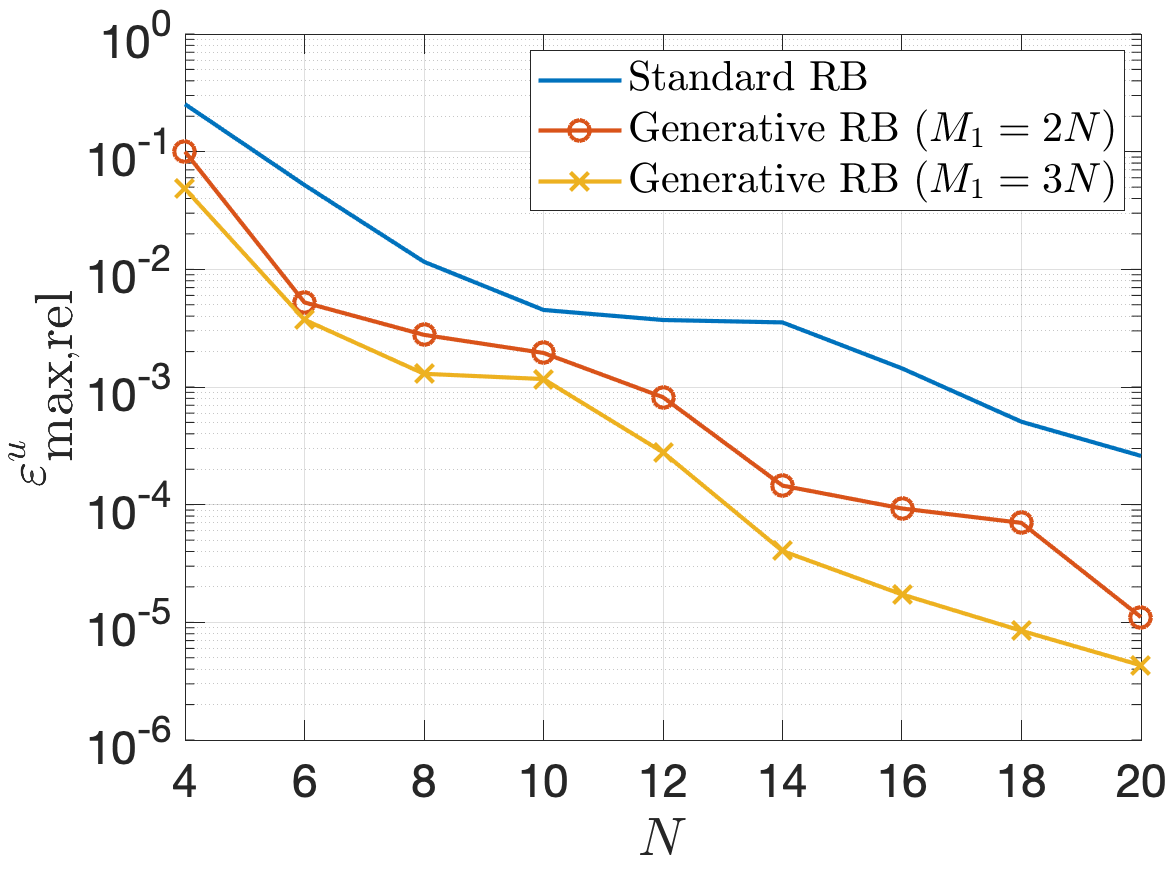}		
	\end{subfigure}        
 	\begin{subfigure}[b]{0.32\textwidth}
		\centering		\includegraphics[width=\textwidth]{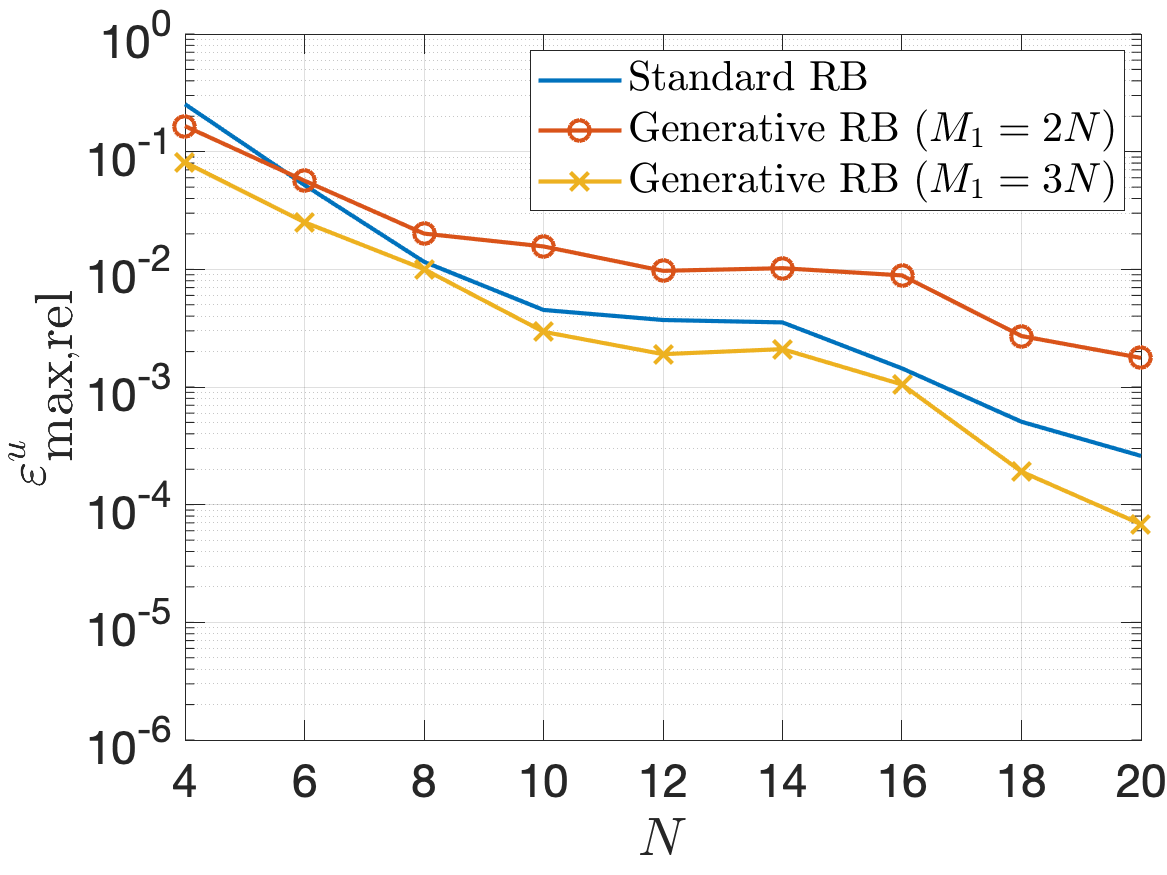}		
	\end{subfigure}
        \\
	\begin{subfigure}[b]{0.32\textwidth}
		\centering		\includegraphics[width=\textwidth]{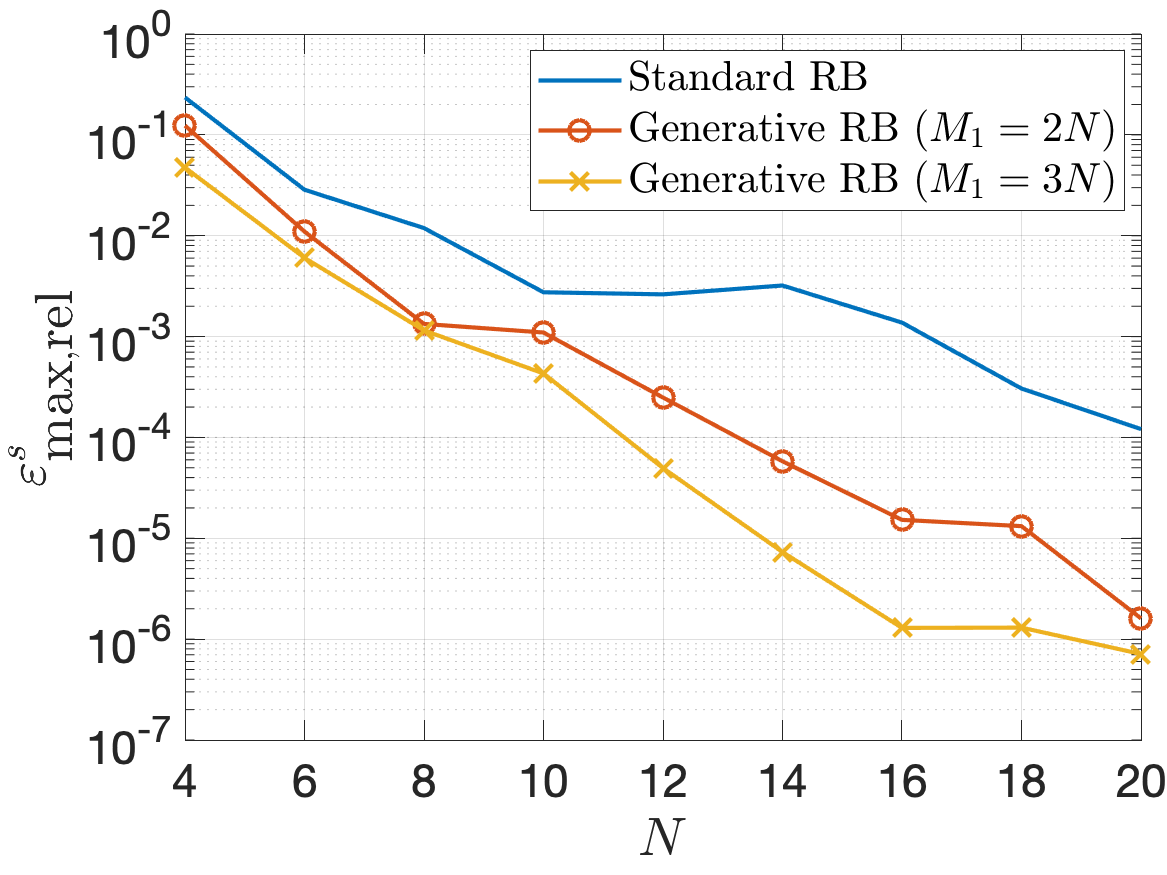}		
        \caption{Softplus function}
	\end{subfigure}         
  	\begin{subfigure}[b]{0.32\textwidth}
		\centering		\includegraphics[width=\textwidth]{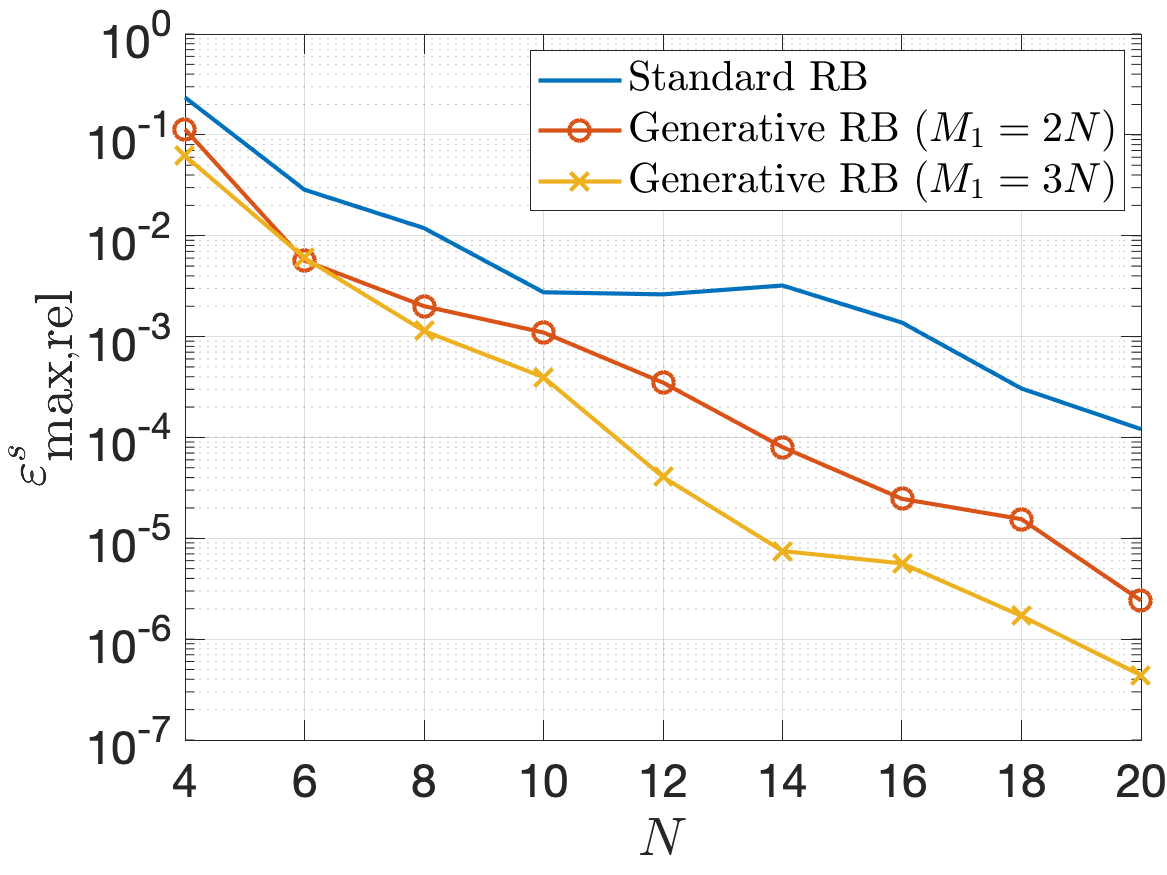}		
        \caption{Exponential function}
	\end{subfigure}
	\begin{subfigure}[b]{0.32\textwidth}
		\centering		\includegraphics[width=\textwidth]{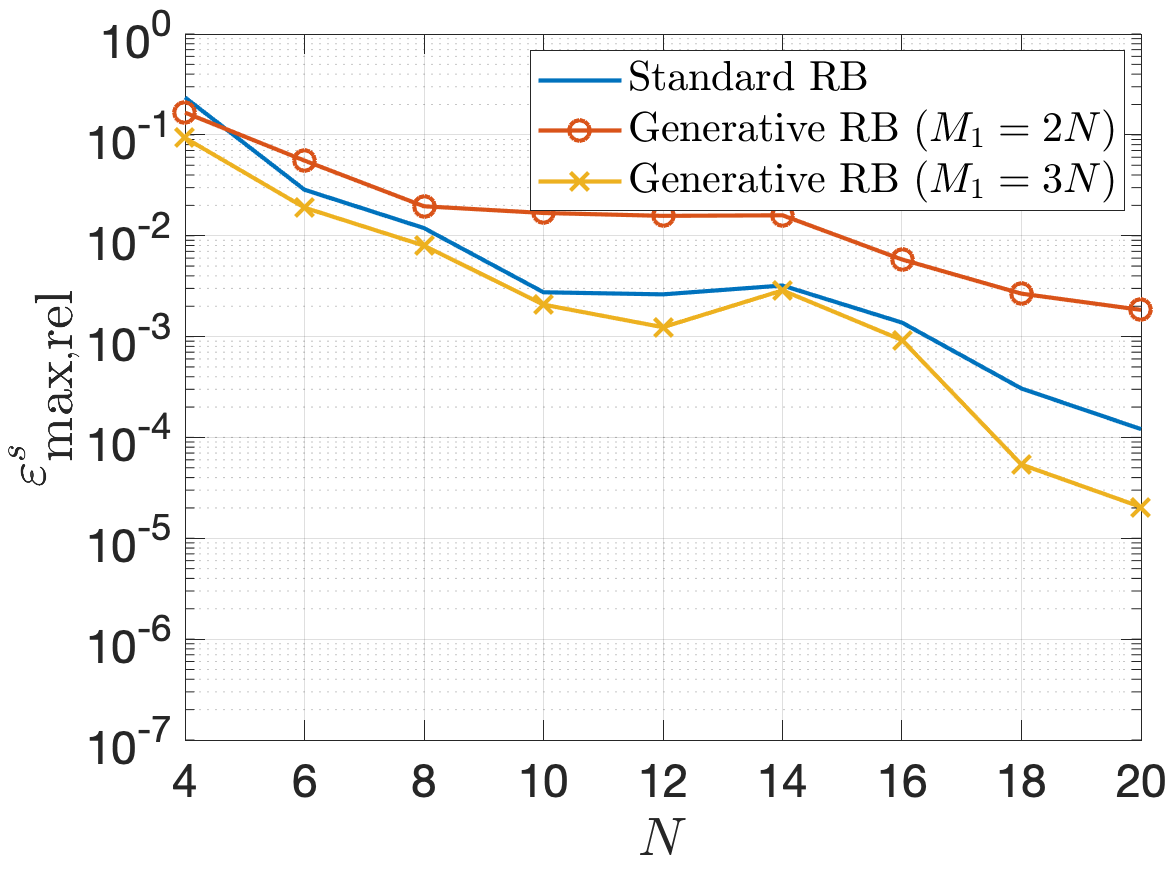}		
    \caption{Quadratic function}
	\end{subfigure}
	\caption{Convergence of the maximum relative error in the RB solution (top) and in the RB output (bottom) as a function of $N$.}
	\label{ex5fig1}
\end{figure}

Figure \ref{ex5fig1} presents the maximum relative errors, ${\epsilon}_{\rm{max},{\rm rel}}^u$ and ${\epsilon}_{\rm{max},{\rm rel}}^s$, as a function of $N$ for both the standard RB method and the generative RB method. The test set $\Xi^{\rm test}$ is a uniform grid of $N^{\rm test} = 6 \times 6 \times 6 \times 6$ points in the parameter domain. The generative RB method consistently shows lower errors compared to the standard RB method for the softplus and exponential functions. The softplus and exponential functions exhibit similar convergence behaviors. For the exponential function with $N = 20$ and $M_1 = 60$, the generative RB method produces output errors that are approximately 100 times smaller than those of the standard RB method. Unlike the previous example, the quadratic function performs poorly in this example.

\begin{figure}[h!]
	\centering
	\begin{subfigure}[b]{0.32\textwidth}
		\centering		\includegraphics[width=\textwidth]{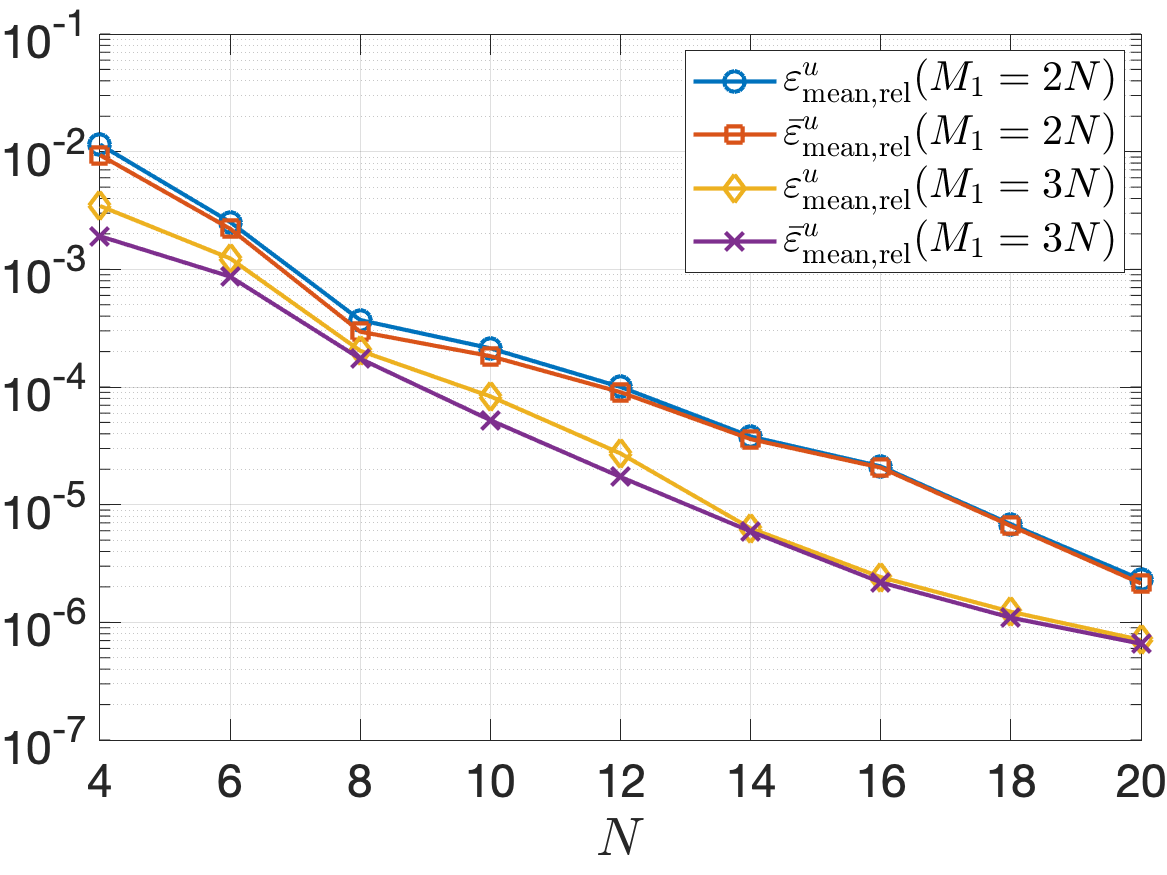}		
	\end{subfigure}	
	\begin{subfigure}[b]{0.32\textwidth}
		\centering		\includegraphics[width=\textwidth]{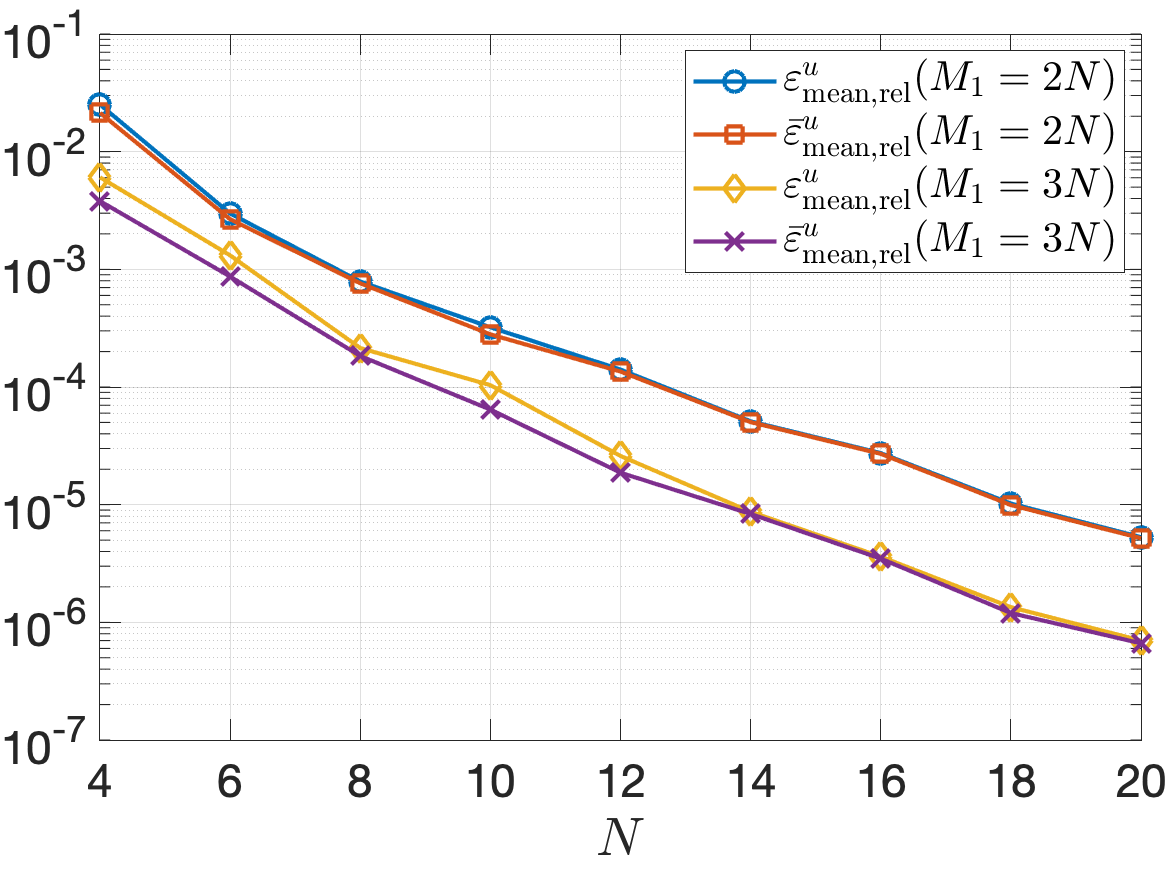}		
	\end{subfigure}        
 	\begin{subfigure}[b]{0.32\textwidth}
		\centering		\includegraphics[width=\textwidth]{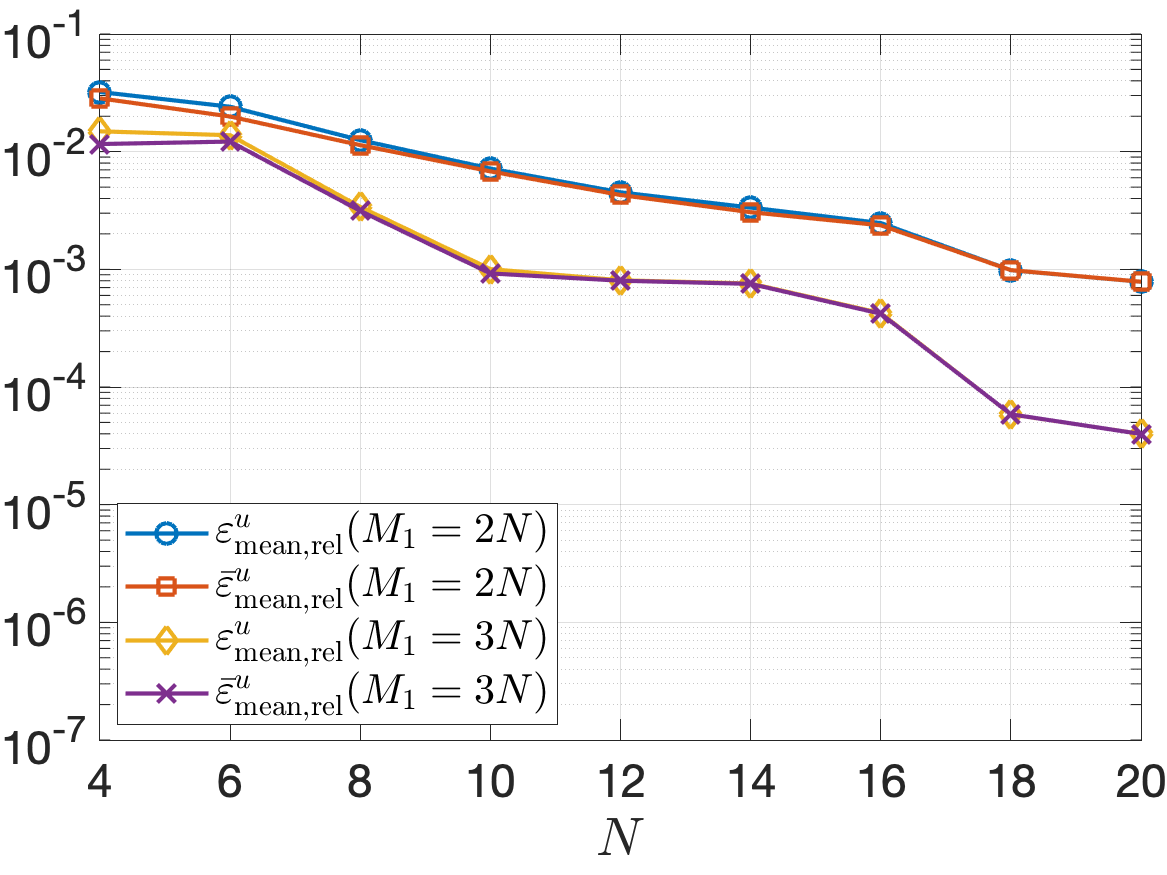}		
	\end{subfigure}
        \\
	\begin{subfigure}[b]{0.32\textwidth}
		\centering		\includegraphics[width=\textwidth]{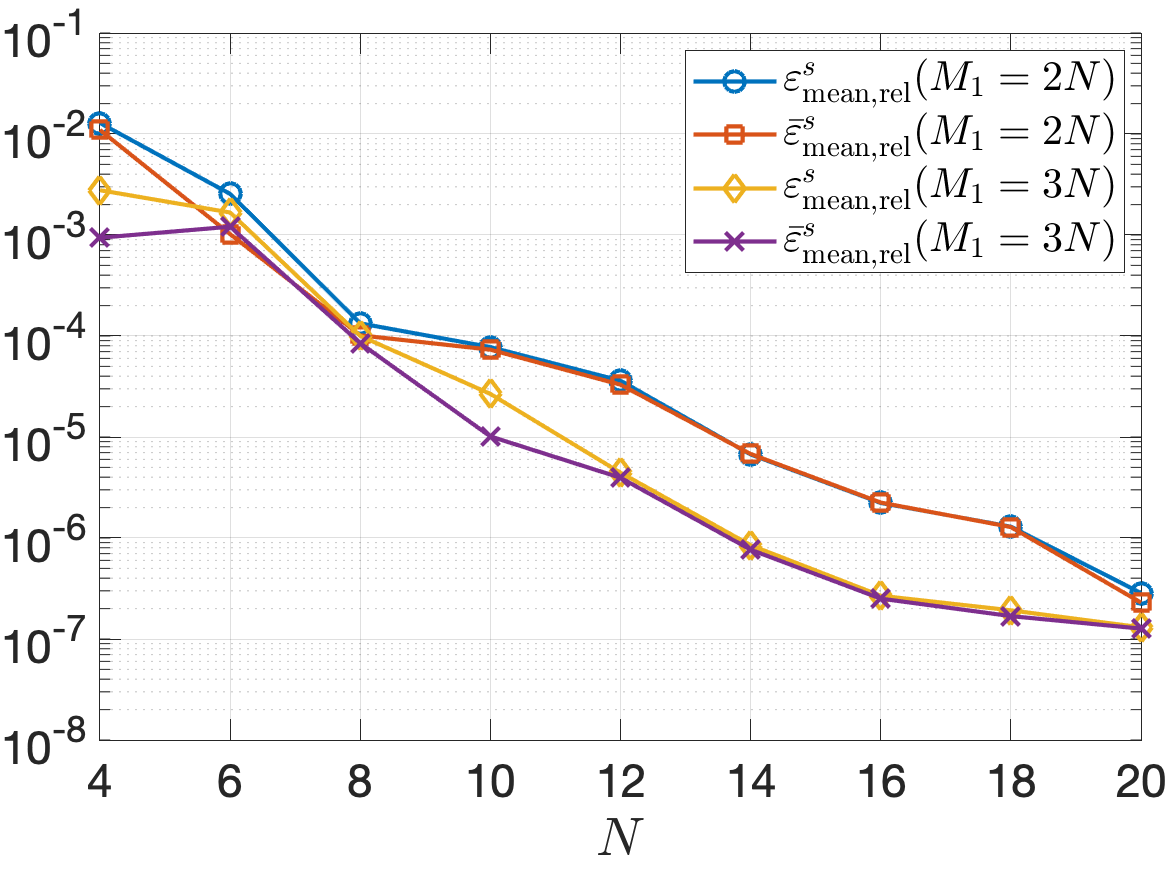}		
        \caption{Softplus function}
	\end{subfigure}         
  	\begin{subfigure}[b]{0.32\textwidth}
		\centering		\includegraphics[width=\textwidth]{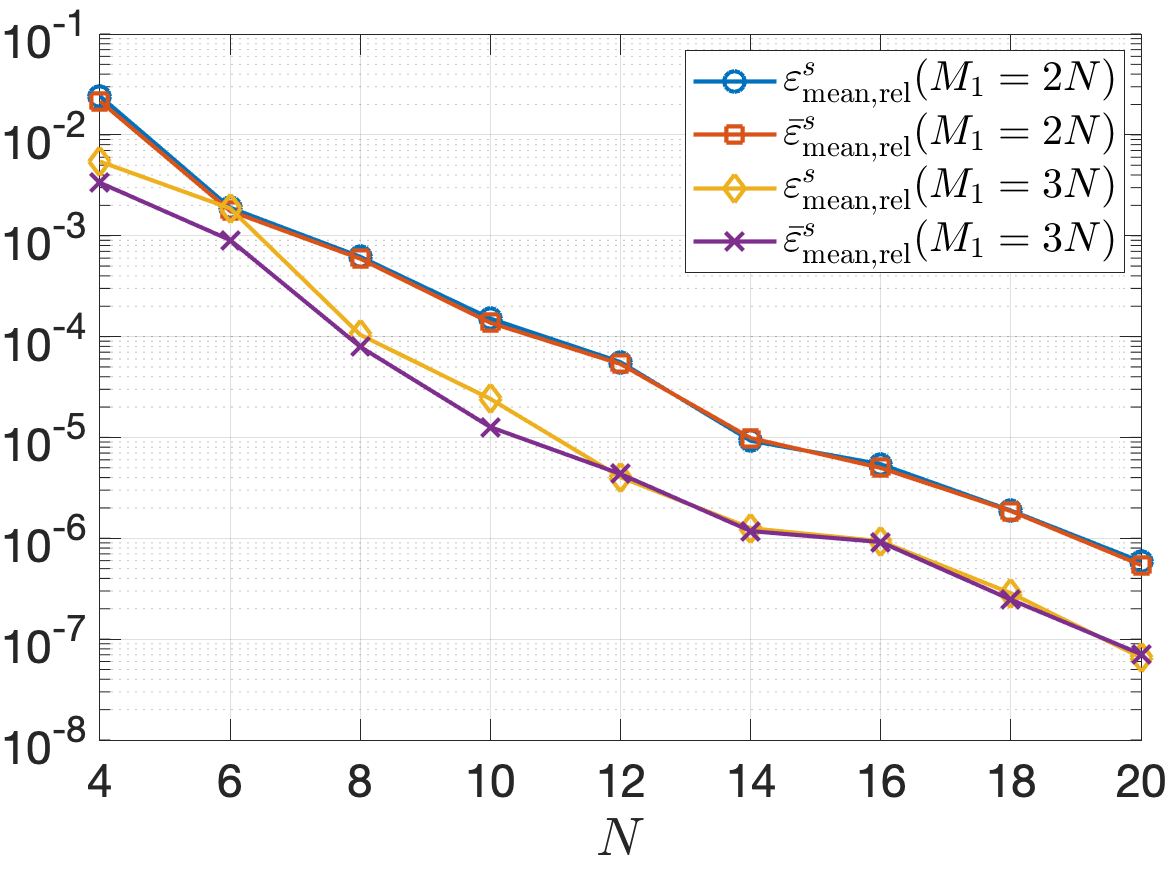}		
        \caption{Exponential function}
	\end{subfigure}
	\begin{subfigure}[b]{0.32\textwidth}
		\centering		\includegraphics[width=\textwidth]{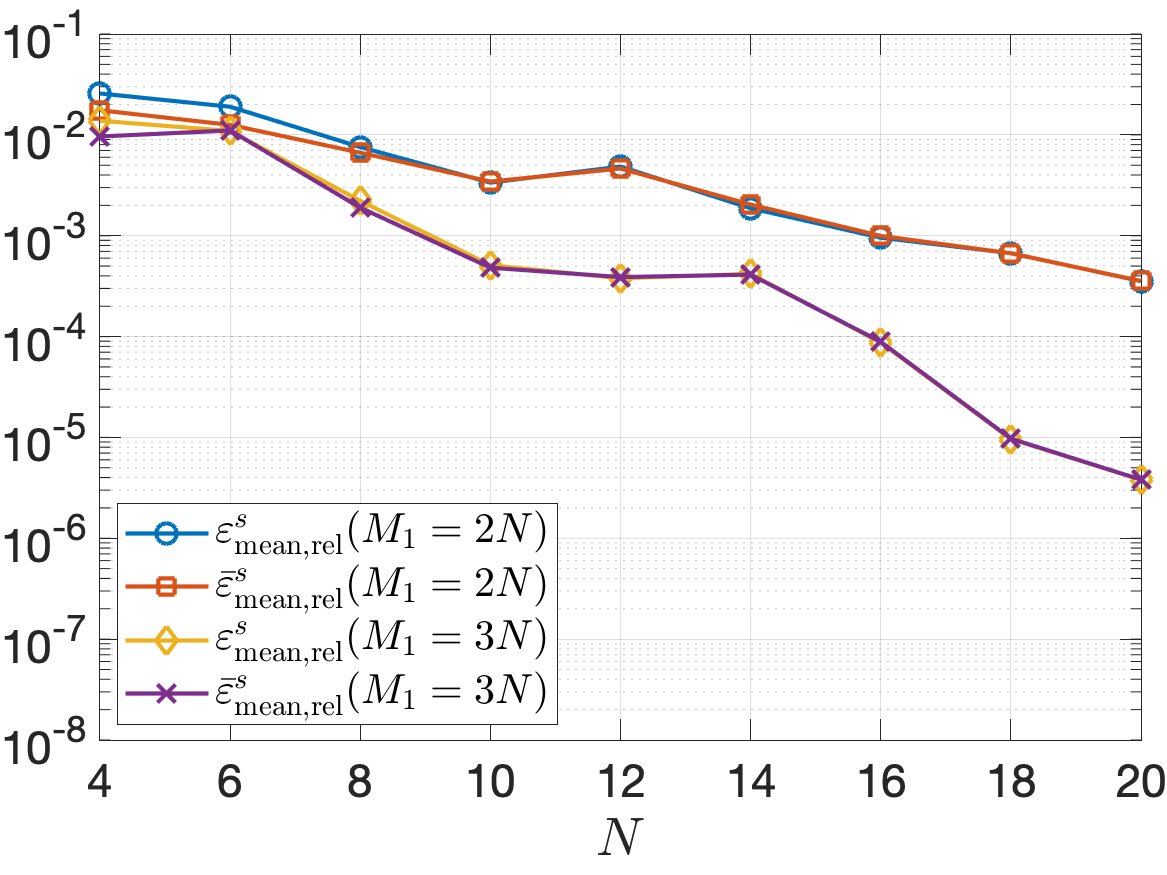}		
    \caption{Quadratic function}
	\end{subfigure}
	\caption{Convergence of the mean relative error and the mean error estimate in the RB solution (top) and in the RB output (bottom) as a function of $N$. The error estimates are calculated by using $M_2 = M_1 + N$.}
	\label{ex5fig2}
\end{figure}

Figure \ref{ex5fig2} shows the mean relative errors, ${\epsilon}_{\rm{mean},{\rm rel}}^u$ and ${\epsilon}_{\rm{mean},{\rm rel}}^s$, together with their error estimates, $\bar{\epsilon}_{\rm{mean},{\rm rel}}^u$ and $\bar{\epsilon}_{\rm{mean},{\rm rel}}^s$, as functions of $N$. The error estimates are calculated using $M_2 = M_1 + N$. The error estimates closely match the true errors and thus provide a reliable measure of accuracy. These results highlight the importance of selecting an appropriate activation function and demonstrate the potential of the generative RB method for achieving high accuracy for problems in high-dimensional parameter spaces.

\section{Conclusion}

In this paper, we introduced generative reduced basis (RB) methods to construct reduced-order models (ROMs) for parametrized partial differential equations. While the generative RB methods have demonstrated significant advantages compared to the standard RB method, there are several avenues for future work to further advance the approach.

Given that the choice of nonlinear activation function plays a critical role in the performance of the generative RB method, it becomes evident that optimizing the activation function is paramount for achieving optimal accuracy and efficiency. Future work should focus on developing systematic procedures to identify the most suitable activation function for a given problem. This could involve adaptive learning techniques that adjust the activation function based on error metrics, parameter space exploration, or the use of  optimization algorithms to tailor the activation function to the specific characteristics of the solution manifold. Such an approach would further enhance the generative RB method and make it even more effective for a wider range of complex high-dimensional problems.

One key area for future research is the extension of generative RB methods to non-affine and nonlinear problems, where traditional linear approximation techniques face greater challenges. Developing robust strategies to handle such complexities will broaden the applicability of generative RB methods to a wider range of physical and engineering problems.

Another direction is the integration of adaptive strategies within the generative framework. An adaptive refinement process could dynamically adjust the number and distribution of snapshots, further optimizing the construction of RB spaces based on the problem’s specific characteristics. This would enable more efficient ROMs, especially for problems where solution behavior varies significantly across the parameter space.

Additionally, machine learning techniques could be incorporated to enhance the scalability and generalization capabilities of generative RB methods. For example, data-driven ROMs could be trained on  smaller sets of solution snapshots by using the generative snapshot method. This would significantly reduce the offline computational cost and improve the overall efficiency of the method.

Lastly, there is potential for developing nonlinear manifold methods that go beyond linear approximation spaces. These methods could allow for the efficient approximation of the solution manifold by using nonlinear embeddings, particularly in cases where traditional RB methods exhibit slow convergence due to highly nonlinear features. However, a challenge associated with nonlinear manifold methods is the requirement for extensive snapshot sets to train the nonlinear mappings effectively. The generative snapshot method can be used to enrich the snapshot set and thus enhance the efficiency and accuracy of  nonlinear manifold methods.

\section*{Acknowledgements} \label{}
I would like to thank Professors Jaime Peraire, Anthony T. Patera, and Robert M. Freund at MIT for fruitful discussions. I gratefully acknowledge a Seed Grant from the MIT Portugal Program, the United States Department of Energy under contract DE-NA0003965, the National Science
Foundation under grant number NSF-PHY-2028125, and the Air Force Office of Scientific Research under Grant No. FA9550-22-1-0356 for supporting this work.

 \bibliographystyle{elsarticle-num} 
\bibliography{library.bib}




\end{document}